\documentclass[11pt,leqno]{article}

\usepackage{amsmath,amsfonts,amscd,amssymb,theorem}

\long\def\comment#1\endcomment{}

\makeatletter
\begingroup
\gdef\th@dotted{\normalfont\itshape
  \def\@begintheorem##1##2{%
        \item[\hskip\labelsep \theorem@headerfont ##1\ ##2.]}%
\def\@opargbegintheorem##1##2##3{%
   \item[\hskip\labelsep \theorem@headerfont ##1\ ##2\ (##3).]}}
\endgroup
\makeatother

\theoremstyle{dotted}

\newtheorem{theorem}{Theorem}[section]
\newtheorem{lemma}[theorem]{Lemma}
\newtheorem{conj}[theorem]{Conjecture}
\newtheorem{prop}[theorem]{Proposition}
\newtheorem{corr}[theorem]{Corollary}

\makeatletter
\begingroup
\gdef\th@upshape{\normalfont
  \def\@begintheorem##1##2{%
        \item[\hskip\labelsep \theorem@headerfont ##1\ ##2.]}%
\def\@opargbegintheorem##1##2##3{%
   \item[\hskip\labelsep \theorem@headerfont ##1\ ##2\ (##3).]}}
\endgroup
\makeatother

\theoremstyle{upshape}

\newtheorem{defn}[theorem]{Definition}
\newtheorem{remark}[theorem]{Remark}
\newtheorem{exa}[theorem]{Example}

\makeatletter
\renewcommand{\subsection}{\@startsection{subsection}{2}{0pt}{-3ex
plus -1ex minus -0.2ex}{-2mm plus -0pt minus
-2pt}{\normalfont\bfseries}} 
\renewcommand{\subsubsection}{\@startsection{subsubsection}{3}{0pt}{-3ex
plus -1ex minus -0.2ex}{-2mm plus -0pt minus
-2pt}{\normalfont\bfseries}} 
\makeatother

\makeatletter
\@addtoreset{equation}{section}
\makeatother

\newcommand{\cntrct}                
{\hspace{2pt}\raisebox{1pt}{\text{$\lrcorner$}}\hspace{2pt}}

\newcommand{\proof}[1][Proof.]{\smallskip\noindent{\em #1}}
\def\endproof{\hfill\ensuremath{\square}\par\medskip}

\def\eqref#1{\thetag{\ref{#1}}}

\let\latexref=\ref
\def\ref#1{{\normalfont{\latexref{#1}}}}

\newcommand{\wt}{\widetilde}
\newcommand{\wh}{\widehat}
\newcommand{\ol}{\overline}
\newcommand{\br}[1]{\langle #1 \rangle}

\newcommand{\ratto}{\dasharrow}        

\newcommand{\idot}{{\:\raisebox{1pt}{\text{\circle*{1.5}}}}}
\newcommand{\hdot}{{\:\raisebox{3pt}{\text{\circle*{1.5}}}}}
\newcommand{\Z}{{\mathbb Z}}
\newcommand{\N}{{\mathbb N}}
\newcommand{\RR}{{\mathbb R}}

\newcommand{\bC}{\overline{\C}}
\newcommand{\bGamma}{\overline{\Gamma}}

\newcommand{\SQ}{\phi}
\newcommand{\Sp}{\operatorname{\sf Sp}}

\newcommand{\eps}{\varepsilon}
\renewcommand{\phi}{\varphi}

\newcommand{\F}{{\sf F}}

\newcommand{\vH}{\check{H}}
\newcommand{\HH}{{\mathcal{H}}}

\def\dlim_#1{{\displaystyle\lim_{#1}}^\hdot}

\newcommand{\Hom}{\operatorname{Hom}}
\newcommand{\Ext}{\operatorname{Ext}}
\newcommand{\RHom}{\operatorname{RHom}}

\newcommand{\Ker}{\operatorname{Ker}}

\newcommand{\Fun}{\operatorname{Fun}}

\newcommand{\id}{\operatorname{\sf id}}
\newcommand{\Id}{\operatorname{\sf Id}}
\newcommand{\gr}{\operatorname{\sf gr}}

\newcommand{\A}{{\cal A}}
\newcommand{\D}{{\cal D}}
\newcommand{\DM}{{\cal D}{\cal M}}
\newcommand{\DS}{{\cal D}{\cal S}}
\newcommand{\DQ}{{\cal D}{\cal Q}}
\newcommand{\DT}{{\cal D}{\cal T}}

\newcommand{\C}{{\cal C}}
\newcommand{\B}{{\cal B}}

\newcommand{\Q}{{\cal Q}}
\newcommand{\T}{{\cal T}}

\newcommand{\R}{{\cal R}}

\newcommand{\I}{{\sf I}}

\newcommand{\Sets}{\operatorname{Sets}}

\newcommand{\Cat}{\operatorname{Cat}}
\newcommand{\Maps}{\operatorname{Maps}}

\newcommand{\Aut}{{\operatorname{Aut}}}
\newcommand{\End}{{\operatorname{End}}}

\newcommand{\Ind}{\operatorname{\sf Ind}}

\newcommand{\Sthom}{\operatorname{\sf StHom}}
\newcommand{\sthom}{\operatorname{\sf sthom}}
\newcommand{\SW}{\operatorname{\sf SW}}
\newcommand{\GSthom}{G\operatorname{\sf -StHom}}

\newcommand{\Add}{\operatorname{\sf Add}}

\newcommand{\amod}{{\text{\rm -mod}}}

\newcommand{\ppt}{{\sf pt}}

\newcommand{\wPsi}{\widetilde{\Psi}}
\newcommand{\wPhi}{\widetilde{\Phi}}
\newcommand{\wphi}{\widetilde{\phi}}
\newcommand{\hPhi}{\widehat{\Phi}}

\newcommand{\Ab}{\operatorname{Ab}}

\newcommand{\M}{\operatorname{\mathcal M}}
\newcommand{\Mm}{\operatorname{\sf M}}

\newcommand{\Restr}{\operatorname{\sf Restr}}
\newcommand{\Infl}{\operatorname{\sf Infl}}

\newcommand{\Ho}{\operatorname{\sf Ho}}
\newcommand{\Ass}{\operatorname{\sf Ass}}

\newcommand{\copr}{{\textstyle\coprod}}

\title{Derived Mackey functors}

\author{D. Kaledin\thanks{Partially supported by grant NSh-1987.2008.1}}

\date{{\em To Pierre Deligne, on the occasion of his 65th birthday}}

\begin{document}

\maketitle

\tableofcontents

\section*{Introduction.}

The notion of a ``Mackey functor'' associated to a finite group $G$
is a standard tool both in algebraic topology and in group
theory. It was originally introduced by Dress \cite{dress} and later
clarified by several people, in particular by Lindner \cite{lind};
the reader can find modern expositions in the topological context
e.g. in \cite{may1}, \cite{may2}, \cite{tD}, or a more algebraic
treatment in \cite{the}.

In this paper, we will be mostly concerned with applications to
algebraic topology. Of these, the main one is the following: the
category of Mackey functors is the natural target for equivariant
homology and cohomology.

Namely, assume given a CW complex $X$ equipped with a continuous
action of a finite group $G$. Then if the action is nice enough, the
cellular homology complex $C_\idot(X,\Z)$ inherits a $G$-action, so
that we can treat homology as a functor from $G$-equivariant CW
complexes to the derived category $\D(G,\Z)$ of representations of
the group $G$. However, this loses some essential information. For
example, for any subgroup $H \subset G$, the homotopy type of the
subspace $X^H \subset X$ of $H$-fixed points is a $G$-homotopy
invariant of $X$ in a suitable sense; but once we forget $X$ and
remember only the object $C_\idot(X,\Z) \in \D(G,\Z)$, there is no
way to recover the homology $H_\idot(X^H,\Z)$.

Thus, the Mackey functors: certain algebraic gadgets designed to
remember not only the homology $H_\idot(X,\Z)$ as a representation
of $G$, but also all the groups $H_\idot(X^H,\Z)$, $H \subset G$,
with whatever natural group action they possess, and some natural
maps between them. We recall the precise definitions in
Section~\ref{mack.rec}. For now, it suffices to say that in the
standard approach, Mackey functors form a tensor abelian category
$\M(G)$ such that, among other things,
\begin{enumerate}
\item for any $G$-equivariant CW complex $X$, we have natural
  homology objects $H_\idot^G(X,\Z) \in \M(G)$,
\item there is a forgetful exact tensor functor from $\M(G)$ to the
  category $\Z[G]\amod$ of representations of $G$ which recovers
  $H_\idot(X,\Z)$ with the natural $G$-action when applied to
  $H^G_\idot(X,\Z)$,
\item for any subgroup $H \subset G$, the homology $H_\idot(X^H,\Z)$
  can also be recovered from $H^G_\idot(X,\Z) \in \M(G)$,
\item $H^G_\idot(X,\Z)$ is compatible with stabilization and the
  tensor product, and extends to the ``genuine $G$-equivariant
  stable homotopy category'' of \cite{may1}, here denoted by
  $\Sthom(G)$.
\end{enumerate}
More precisely, for every subgroup $H \subset G$, one has an exact
functor from $\M(G)$ to the category of abelian groups which
associates an abelian group $M^H$ to every $M \in \M(G)$; then in
\thetag{iii}, there is a functorial isomorphism
$$
H_\idot(X^H,\Z) \cong H^G_\idot(X,\Z)^H.
$$
One can use the correspondence $M \mapsto M^H$ to visualize the
structure of the category $\M(G)$ in the following way.  For any
subgroup $H \subset G$, let $\M_H(G) \subset \M(G)$ be the full
subcategory spanned by such $M \in \M(G)$ that
\begin{itemize}
\item $M^{H'} = 0$ unless $H'$ contains a conjugate of $H$.
\end{itemize}
Then this is a Serre abelian subcategory, and the subcategories
$\M_H(G)$, $H \subset G$, form an increasing ``filtration by
support'' of the category $\M(G)$. The top associated graded
quotient of this filtration is equivalent to the category
$\Z[G]$-modules --- that is, we have
$$
\M(G)/\langle\M_H(G)\rangle_{\{e\} \neq H \subset G} \cong
\Z[G]\amod,
$$
where $\langle\M_H(G)\rangle_{\{e\} \neq H \subset G} \subset \M(G)$
is the Serre subcategory generated by $\M_H(G)$ for all subgroups $H
\subset G$ except for the trivial subgroup $\{e\} \subset G$. One
can also compute other quotients; for example, the smallest
subcategory $\M_G(G) \subset G$ corresponding to $G$ itself is
equivalent to the category $\Z\amod$ of abelian groups. More
generally, for any normal subgroup $N \subset G$ there exists a
fully faithful exact {\em inflation functor} $\Infl^N_G$ which gives
an equivalence
$$
\Infl^N_G:\M(G/N) \cong \M_N(G) \subset \M(G).
$$
Analogous structures also exist on the category $\Sthom(G)$. Namely,
for any $G$-spectrum $X \in \Sthom(G)$, one has the so-called {\em
Lewis-May fixed points spectrum} $X^H$, so that one can define the
subcategories $\Sthom_H(G)$ by \thetag{$\bullet$}. Then for a normal
subgroup $N \subset G$, one has a fully faithful embedding
$\Sthom(G/N) \cong \Sthom_N(G) \subset \Sthom(G)$. In particular,
the smallest subcategory $\Sthom_G(G) \subset \Sthom(G)$ is
equivalent to the non-equivariant stable homotopy category
$\Sthom$. Another important feature of the stable category is the
{\em geometric fixed points functor} $\Phi^H:\Sthom(G) \to \Sthom$;
on the level of Mackey functors, this corresponds to the projections
onto the associated graded quotients of the filtration by support.

\medskip

To a person trained in homological algebra, a natural next thing to
do is to consider the {\em derived} category $\D(\M(G))$ of the
abelian category $\M(G)$, and try to extend all of the above to the
``derived level'': one would like to have a natural equivariant
homology functor $C_\idot^G(-,\Z):\GSthom \to \D(\M(G))$, and one
would expect the category $\D(\M(G))$ to imitate the natural
structure of the category $\Sthom(G)$.

Unfortunately, and this came as a nasty surprise to the author, this
program does not work: the derived category $\D(\M(G))$ is not the
right thing to consider. The specific problem is that the inflation
functor $\Infl^N_G$ is not fully faithful on the level of derived
categories. Already in the case $G = \Z/p\Z$, $p$ prime, when the
only subgroups in $G$ are the trivial subgroup $\{e\} \subset G$ and
$G$ itself, we can consider the full triangulated subcategory
$\D_G(\M(G)) \subset \D(\M(G))$ spanned by such $M \in \D(\M(G))$
that $M^{\{e\}}=0$. Then while we do have the equivalence
$\D(\M(G))/\D_G(\M(G)) \cong \D(G,\Z)$, the functor
$$
\D(\Z\amod) \cong \D(\M_G(G)) \to \D_G(\M(G))
$$
is not an equivalence. The category $\D_G(\M(G))$ which ought to be
equivalent to the derived category of abelian groups is in fact
rather complicated and behaves badly. So, while one might be able to
construct a homology functor $\Sthom(G) \to \D(\M(G))$, it does not
seem to reflect the structure of $\Sthom(G)$ too closely, and in
particular, one cannot expect any reasonable compatibility with the
geometric fixed point functors $\Phi^H$.

\medskip

But fortunately, a moment's reflection on the definition of a Mackey
functor (to wit, the argument in the beginning paragraphs of
Section~\ref{mack.der}) shows why this is so, and in fact suggests
what the correct ``category of derived Mackey functors'' should
be. This is the subject of the present paper. For any finite group
$G$, we construct a tensor triangulated category $\DM(G)$ of
``derived Mackey functors'' which enjoys the following properties.
\begin{enumerate}
\item For any subgroup $H \subset G$ and any $M \in \DM(G)$, there
  exists a functorial ``fixed points'' object $M^H \in \D(\Z\amod)$.
\item For any subgroup $H \subset G$, let $\DM_H(G) \subset \DM(G)$
  be the full triangulated subcategory spanned by $M \in \DM(G)$
  satisfying \thetag{$\bullet$}. Then $\DM_H(G) \subset \DM(G)$ is
  admissible in the sense of \cite{boka} -- that is, the embedding
  functor $\DM_H(G) \hookrightarrow \DM(G)$ has a left and a
  right-adjoint -- and we have an equivalence
\begin{equation}\label{quo.eq}
\DM_H(G)/\langle\DM_{H'}(G)\rangle_{H \subset H' \subset G,H \neq H'}
\cong \D(W_H,\Z),
\end{equation}
where $W_H = N_H/H$ is the quotient of the normalizer $N_H \subset
G$ of $H \subset G$ by $H$ itself.
\item For any normal subgroup $N \subset G$, we have an equivalence
$$
\wh{\Infl}^N_G:\DM(G/N) \cong \DM_N(G) \subset \DM(G).
$$
\end{enumerate}
Informally speaking, the subcategories $\DM_H(G)$ form a filtration
of the category $\DM(G)$ indexed by the lattice of conjugacy classes
of subgroups in $G$; this filtration gives rise to a
``semiorthogonal decomposition'' in the sense of \cite{boka}, and
the graded pieces of the filtration are naturally identified with
$\D(W_H,\Z)$, $H \subset G$. We explicitly construct functors
$$
\wPhi^{[G/H]}:\DM(G) \to \D(W_H,\Z)
$$
that give the projections onto these graded pieces. We also compute
the gluing functors between $\D(W_H,\Z) \subset \DM(G)$; these are
naturally expressed in terms of a certain generalization of Tate
cohomology of finite groups.

The reader will notice that the properties of the category $\DM(G)$
are slightly stronger than what we have mentioned for the abelian
category $\M(G)$, in that we identify the associated graded pieces
of the filtration by support, and describe how these pieces are
glued. If one localizes the category $\M(G)$ by inverting the order
of the group $G$, then the corresponding statement for the localized
category $\wt{\M}(G)$ is a theorem of Thevenaz \cite{Th1} (in this
case, there is no gluing: the category $\M(G)$ is semisimple, and it
splits into a direct sum of the categories of representations of the
groups $W_H$, $H \subset G$). I do not know whether anything is
known for $\M(G)$ in the general non-semisimple case. I also do not
know whether analogous statements are knows for the category
$\Sthom(G)$ --- namely, whether the graded pieces of the filtration
by the subcategories $\Sthom_H(G) \subset \Sthom(G)$ have been
computed and/or whether the gluing functors have been identified.

Moreover, in this paper we construct a natural equivariant homology
functor $h^G$ from $G$-spectra to our derived Mackey functors such
that
\begin{enumerate}
\item the functor $h^G$ is tensor,
\item for any subgroup $H \subset G$, the functor $h^G$ sends the
 subcategory $\Sthom_H(G) \subset \Sthom(G)$ into $\DM_H(G) \subset
 \DM(G)$, and
\item for any subgroup $H \subset G$ and a $G$-spectrum $X$, the
underlying complex of $\wPhi^{[G/H]}(h^G(X)) \in \D(W_H,\Z)$ computes
the homology of the geometric fixed point spectrum $\Phi^H(X)$.
\end{enumerate}
Unfortunately, we can only do all of the above for {\em finite}
$G$-CW spectra $X$. I do not know whether the corresponding
statements are true for the whole category $\Sthom(G)$. I also
expect that for any $G$-spectrum $X$, $h^G(X)^H$ is the homology of
the Lewis-May fixed points spectrum $X^H$, but I do not presently
know how to prove it.

\medskip

The paper is organized as follows. In Section~\ref{pre.sec}, we give
some necessary standard facts from homological algebra and category
theory (this includes the notions related to semiorthogonal
decompositions of triangulated categories). In
Section~\ref{mack.rec}, we recall the usual definition of Mackey
functors and some of their basic properties. Then we give our
derived version. We actually give not one but two
definitions. First, we give a rather explicit definition using
$A_\infty$-categories and bar resolutions -- this is
Section~\ref{mack.der} (in fact, we work in a slightly larger
generality of a small category $\C$ which has fibered products --
for Mackey functors, this is the category of finite $G$-sets). Then
in Section~\ref{wald.sec}, we give a more invariant definition
somewhat in the spirit of Waldhausen's $S$-construction, and we show
that the two definitions are equivalent. In Section~\ref{pro.sec},
we show how to extend the basic properties of Mackey functors given
in Section~\ref{mack.rec} to the derived setting (in particular, we
construct the tensor product, the inflation functors
$\wh{\Infl}^N_G$ and the geometric fixed points functors
$\Phi^{[G/H]}$).

Formally, this is where the paper might have ended; however, neither
of our two equivalent definitions is suitable for computations. For
this reason, we give a third rather explicit description of the same
category, and this is the subject of the lengthy
Section~\ref{galois.sec}. Essentially, we do a sort of Koszul
duality --- we try to describe the same category $\DM(G)$ using the
geometric fixed points functors $\Phi^{[G/H]}$ as a ``fiber
functor''. This is surprisingly delicate; in particular, we have to
work with $A_\infty$-comodules over an $A_\infty$-coalgebra instead
of the more usual $A_\infty$-modules over an $A_\infty$-algebra.

One immediate advantage is that the functors $\Phi^{[G/H]}$ are
tensor, so that this new description is better compatible with the
tensor product in $\DM(G)$. However, the real reward comes in
Section~\ref{tate.sec}: we are able to prove the equivalences
\eqref{quo.eq}, thus describing the category of derived Mackey
functors as a successive extension of representation categories of
the subquotient groups $W_H$ of the group $G$, and we express the
gluing data between these representation categories in terms of a
generalization of Tate (co)homology. This ``generalized Tate
cohomology'' has the great advantage of being trivial in many cases;
even when it is not trivial, it is usually possible to compute it.

In the final Section~\ref{spectra}, we describe the relation between
the category of derived Mackey functors and the $G$-equivariant
stable homotopy category. Not being an expert on stable homotopy, I
have kept the exposition to an absolute minimum; however, I hope
that Section~\ref{spectra} does show that the category of derived
Mackey functors imitates the equivariant stable homotopy category in
a satisfactory way.

\medskip

It should be stressed that a topologist would learn almost nothing
from this paper: pretty much everything that we prove about derived
Mackey functors is well-known for equivariant spectra (possibly in a
different language). In a sense, the whole point of the paper is
that so much survives in the purely homological theory, which is
usually pretty trivial compared to the richness of the sphere
spectrum. On the other hand, things standard in algebraic topology
are not always well-known outside of it; a reader with a more
geometric and/or homological background can treat the paper as an
extended exercise in the theory of cohomological descent. As such,
it might even be useful, e.g. in the theory of Artin motives.

In the interest of full disclosure, I should mention that my
personal main reason for doing this research was its application to
the so-called ``cyclotomic spectra'', and a comparison theorem
between Topological Cyclic Homology, on one hand, and a syntomic
version of the Periodic Cyclic Homology, on the other hand. Needless
to say, in the end all of this had to be taken out and relegated to
a separate paper, which is ``in preparation''.

\subsection*{Acknowledgements.} It is a pleasure and an honor to
dedicate this paper to Pierre Deligne, on the occasion of his 65-th
birthday. The paper would be impossible without his elegant theory
of cohomological descent. I am very grateful to A. Beilinson,
R. Bezrukavnikov, M. Kontsevich, A. Kuznetsov, G. Merzon, D. Orlov,
L. Positselsky and V. Vologodsky for useful discussions. It was
A. Beilinson who originally asked whether one can do any gluing by
means of Tate homology, and attracted my attention to the problem of
descent for Artin motives. And last but foremost, this paper owes
its existence to conversations with L. Hesselholt, who introduced me
to Mackey functors and stable equivariant theory, and generally
encouraged me to think about the subject. A part of this work was
done during my stay at the University of Tokyo; the hospitality of
this wonderful place and my host Prof. Yu. Kawamata is gratefully
acknowledged. I am extremely grateful to the referee for many useful
suggestions, friendly criticism, and a lot of explanations about the
topological part of the story.

\section{Homological preliminaries.}\label{pre.sec}

\subsection{Generalities.}\label{gen.subs}
Throughout the paper we will fix once and for all an abelian
category $\Ab$; this is the category where our Mackey functors will
take values. We will assume that $\Ab$ is sufficiently nice, in that
it is a Grothendieck abelian category with enough projectives and
injectives. We will also assume that $\Ab$ is equipped with a
symmetric tensor product. In addition, we will assume fixed some
functorial DG enhancement for $\Ab$, so that for any $M,M' \in \Ab$
we have a complex
$$
\RHom^\hdot(M,M')
$$
of abelian groups whose homology computes $\Ext^\hdot(M,M')$, and
moreover, $\RHom^\hdot(-,-)$ is functorial in both arguments and
equipped with functorial and associative composition maps
$$
\RHom^\hdot(M,M') \otimes \RHom^\hdot(M',M'') \to \RHom^\hdot(M,M'')
$$
for any $M,M',M'' \in \Ab$. For example, $\Ab$ may be the category of
vector spaces over a field $k$, or the category of modules over a
commutative ring $R$, or indeed simply the category of abelian
groups. We will denote by $\D(\Ab)$ the (unbounded) derived category
of $\Ab$. We will tacitly assume that the DG enhancement is extended
to unbounded complexes, so that for any two complexes $M_\idot$,
$M_\idot'$ we have a complex of abelian groups
$$
\RHom^\hdot(M_\idot,M'_\idot)
$$
whose $0$-th homology classes are in one-to-one correspondence with
maps from $M_\idot$ to $M'_\idot$ in $\D(\Ab)$ (this is slightly
delicate -- to do it properly, one has to use either $h$-projective
or $h$-injective replacements in the sense of \cite{kel}, \cite{spa}).

Throughout the paper, we will work with many small
categories. Usually we will denote objects in a small category $\C$
by small roman letters; for any $c,c' \in \C$, we will denote by
$\C(c,c')$ the set of morphisms from $c$ to $c'$. For any small
category $\C$, we will denote by $\Fun(\C,\Ab)$ the category of all
functors from $\C$ to $\Ab$. This is also a Grothendieck abelian
category with enough projectives and injectives. We will denote its
derived category by $\D(\C,\Ab)$.

A functor $\gamma:\C \to \C'$ induces a restriction functor
$\gamma^*:\Fun(\C',\Ab) \to \Fun(\C,\Ab)$, and this has a
left-adjoint $\gamma_!:\Fun(\C,\Ab) \to \Fun(\C',\Ab)$ and a
right-adjoint $\gamma_*:\Fun(\C,\Ab) \to \Fun(\C',\Ab)$, known as
the left and right {\em Kan extensions}. By adjunction, $\gamma^*$
is exact, $\gamma_!$ is right-exact, and $\gamma_*$ is left-exact,
so that we have derived functors
$L^\hdot\gamma_!,R^\hdot\gamma_*:\D(\C,\Ab) \to \D(\C',\Ab)$.

In particular, if $\C' = \ppt$ is the point category, we have
$\Fun(\C',\Ab) = \Fun(\ppt,\Ab) \cong \Ab$, and if $\tau:\C \to
\ppt$ is the projection to the point, the Kan extensions $\tau_!$
and $\tau_*$ are given by the direct and inverse limit over the
category $\C$. Their derived functors are known as {\em homology}
$H_\idot(\C,-)$ and {\em cohomology} $H^\hdot(\C,-)$ of the category
$\C$.

On the other hand, if we choose an object $c \in \C$ and let
$i^{c}:\ppt \to \C$ be the functor which sends $\ppt$ to $c \in
\C$, then the Kan extension functors $i^{c}_!,i^{c}_*:\Ab \cong
\Fun(\ppt,\Ab) \to \Fun(\C,\Ab)$ are exact. To simplify notation, we
will denote $M_{c} = i^{c}_!(M)$, $M^{c} = i^{c}_*(M)$ for
any object $M \in \Ab$. Explicitly, the functors $M_{c},M^{c} \in
\Fun(\C,\Ab)$ are given by
\begin{equation}\label{repr.exp}
M_{c}(c') = \bigoplus_{\C(c,c')}M, \qquad
M^{c}(c') = \prod_{\C(c',c)}M,
\end{equation}
for any $c' \in \C$ (the sum and the product of copies of $M$
indexed by elements in the corresponding $\Hom$-sets). If $M$ is
projective, then $M_{c}$ is projective in $\Fun(\C,\Ab)$; if $M$
is injective, then $M^{c}$ is injective in $\Fun(\C,\Ab)$. For any
$M \in \C$, a map $f:c \to c'$ induces natural maps $M_{c'} \to
M_{c}$, $M^{c'} \to M^{c}$.

\subsection{Fibrations and base change.}\label{fib.subs}

In general, it is rather cumbersome to compute the Kan extensions
explicitly; however, there is one situation introduced in \cite{SGA}
where the computations are simplified. Namely, assume given a
functor $\gamma:\C' \to \C$ between small categories $\C$, $\C'$. A
morphism $f':c_0' \to c_1'$ in $\C'$ is called {\em Cartesian} with
respect to $\gamma$ if it has the following universal property:
\begin{itemize}
\item any morphism $f'':c_0'' \to c_1'$ such that $\gamma(f'') =
  \gamma(f')$ factorizes uniquely as $f'' = f' \circ i$ through a
  morphism $i:c_0'' \to c_0'$ such that $\gamma(i) = \id$.
\end{itemize}
The functor $\gamma$ is called a {\em fibration} if 
\begin{enumerate}
\item for any morphism $f:c_0 \to c_1$ in the category $\C$ and any
object $c_1' \in \C'$ with $\gamma(c_1') = c_1$ there exists a
Cartesian map $f':c_0' \to c_1'$ such that $\gamma(f') = f$,
\item and moreover, the composition of two Cartesian maps is
Cartesian.
\end{enumerate}

\begin{exa}\label{fib.fib}
Let $\C$ be a category with fibered products, let $\C'$ be the
category of diagrams $c_0 \to c_1$ in $\C$, and let $\gamma:\C' \to
\C$ be the functor which sends a diagram $c_0 \to c_1$ to
$c_1$. Then $\gamma$ is a fibration.
\end{exa}

For any object $c \in \C$, denote by $\C'_c$ the {\em fiber} of the
functor $\gamma$ over the object $c \in \C$ -- that is, the category
of objects $c' \in \C'$ such that $\gamma(c') = c$ and those
morphisms $i$ between them for which $\gamma(i)=\id$.  Note that by
the universal property of Cartesian maps, the map $f:c_0' \to c_1'$
in \thetag{i} is unique up to a unique isomorphism, so that, if
$\gamma$ is a fibration, setting $f^*(c_1') = c_0'$ defines a
functor $f^*:\C'_{c_1} \to \C'_{c_0}$ which we will call the {\em
transition functor} corresponding to $f$. The same universal property
provides a canonical map $(f \circ g)^* \cong f^* \circ g^*$ for any
two composable morphisms $f$, $g$ in $\C$, and the condition
\thetag{ii} insures that this is an isomorphism.  These isomorphisms
in turn satisfy a compatibility condition for composable triples,
and the whole thing has been axiomatized by Grothendieck under the
name of a ``pseudo-functor'' from $\C$ to the $2$-category $\Cat$ of
small categories; we refer the reader to \cite{SGA} for
details. Grothendieck also proved the inverse statement: every
fibration over $\C$ is uniquely defined by the corresponding
contravariant pseudo-functor from $\C$ to $\Cat$. Nowadays this is
usually called ``the Grothendieck construction''.

Assume given a fibration $\gamma:\C' \to \C$, another small category
$\C_1$, and a functor $\eta:\C_1 \to \C$. Define a small category
$\C_1'$ as a fibered product
$$
\begin{CD}
\C_1' @>{\eta'}>> \C'\\
@V{\gamma'}VV @VV{\gamma}V\\
\C_1 @>{\eta}>> \C.
\end{CD}
$$
Then we have a pair of adjoint {\em base change isomorphisms}
$$
\eta^* \circ R^\hdot\gamma_* \cong R^\hdot\gamma'_* \circ \eta^{'*},
\qquad
L^\hdot\eta'_! \circ \gamma^{'*} \cong \gamma^* \circ L^\hdot\eta_!.
$$
For the proof, see e.g. \cite[Lemma 1.7]{K}.

Dually, $\gamma:\C' \to \C$ is a {\em cofibration} if the
corresponding functor $\gamma^{opp}:\C^{opp'} \to \C^{opp}$ between
the opposite categories is a fibration. Under the Grothendieck
construction, cofibrations correspond to covariant pseudofunctors.
If $\gamma$ is a cofibration, we have a base change isomorphism
$$
\eta^* \circ L^\hdot\gamma_! \cong L^\hdot\gamma'_! \circ \eta^{'*}.
$$
In particular, for any functor $E \in \Fun(\C',\Ab)$, the value
$L^\hdot(E)(c)$ at an object $c \in \C$ is canonically given by
\begin{equation}\label{cofib.bc}
L^\hdot(E)(c) \cong H_\idot(\C'_c,E).
\end{equation}

\subsection{Bar-resolutions.}\label{bar.gene}
Another computational tool that we will need is the so-called {\em
bar-resolution}. Assume given a small category $\C$. Then every
functor $E \in \Fun(\C,\Ab)$ has a canonical resolution
$P_\idot(\C,E)$ with terms
$$
P_{i-1}(\C,E) = \bigoplus_{c_1 \to \dots \to c_i}E(c_1)_{c_i},
$$
where the sum is taken over all the diagrams $c_1 \to \dots \to c_i$
in the category $\C$, and the usual differential $\delta = d_1 - d_2
+ \dots \pm d_i$, where $d_l$ drops the object $c_l$ from the
diagram, and acts as the identity map if $1 < l < i$, as the natural
map $E(c_1)_{c_i} \to E(c_1)_{c_{i-1}}$ induced by $c_{i-1} \to c_i$
if $l = i$, and as the map $E(c_1)_{c_i} \to E(c_2)_{c_i}$ induced
by the map $E(c_1) \to E(c_2)$ if $l=1$. To see that this is indeed
a resolution, one evaluates $P_\idot(\C,E)$ at some object $c \in
\C$. By definition, the resulting complex is given by
$$
P_{i-1}(\C,E)(c) = \bigoplus_{c_1 \to \dots \to c_i \to c}E(c_1),
$$
where the sum is now over all the diagrams ending at $c \in \C$. If
one adds the term $E(c)$ in degree $-1$ corresponding to the diagram
consisting of $c$ itself, then the resulting complex is obviously
chain-homotopic to $0$ --- the contracting homotopy $h$ sends the
term corresponding to a diagram $c_1 \to \dots \to c_i \to c$ to the
term corresponding to $c_1 \to \dots \to c_i \to c \to c$, where the
last map $c \to c$ is the identity map.

Since all the objects $M_{c}$, $c \in \C$, $M \in \Ab$ are obviously
acyclic for the homology functor $H_\idot(\C,-)$, the bar resolution
can be used to compute the homology $H_\idot(\C,E)$. This results in
the {\em bar-complex} $C_\idot(\C,E)$ with terms
\begin{equation}\label{bar.def}
C_{i-1}(\C,E) = \bigoplus_{c_1 \to \dots \to c_i}E(c_1).
\end{equation}
This has the following standard properties.
\begin{enumerate}
\item The bar-complex $C_\idot(\C,E)$ is functorial with respect to
  $E$.
\item For any functor $\gamma:\C \to \C'$ and any $E:\C' \to \Ab$,
  there is a natural map $\gamma^E_*:C_\idot(\C,\gamma^*E) \to
  C_\idot(\C',E)$ which induces the natural adjunction map
  $H_\idot(\C,\gamma^*E) \to H_\idot(\C',E)$ on homology, and for
  any composable pair of functors $\gamma:\C \to \C'$, $\gamma':\C'
  \to \C''$, we have $\gamma^{\gamma^{'*}E}_* \circ \gamma^{'E}_* =
  (\gamma' \circ \gamma)^E_*$.
\item For any two categories $\C$, $\C'$ and functors $E:\C \to
  \Ab$, $E':\C' \to \Ab$, we have a K\"unneth-type quasiisomorphism
\begin{equation}\label{bar.prod}
C_\idot(\C,E) \otimes C_\idot(\C',E') \to C_\idot(\C \times \C',E
\boxtimes E'),
\end{equation}
and this is associative with respect to triple products.
\end{enumerate}
Of these, only \thetag{iii} is slightly non-obvious; the required
quasiisomorphism is given by the shuffle product.

In addition, the bar-resolution can be used to give a canonical DG
enhancement to the category $\Fun(\C,\Ab)$. Indeed, for any $E,E'
\in \Fun(\C,\Ab)$ we can compute $\Ext^\hdot(E,E')$ by the
bar-resolution; this results in a bicomplex
$\RHom^{\hdot,\hdot}(E,E') = \RHom^{\hdot,\hdot}_\C(E,E')$ with
terms
\begin{equation}\label{rhom.bar.eq}
\RHom^{i-1,\hdot}(E,E') = \prod_{c_1 \to \dots \to c_i}
\RHom^\hdot(E(c_1),E'(c_i)).
\end{equation}
We denote by $\RHom^\hdot(E,E')$ the total complex of this
bicomplex. Gives three objects $E,E',E'' \in \Fun(\C,\Ab)$, we have
a natural composition map $\RHom^\hdot(E,E') \otimes
\RHom^\hdot(E',E'') \to \RHom^\hdot(E,E'')$, and this is associative
in the obvious sense.

We will also need a slightly more refined version of the bar
resolution. For any $i$, the diagrams $c_1 \to \dots \to c_i$ and
isomorphisms between them form a groupoid $\C_i$. Denote by
$\sigma_i,\tau_i:\C_i \to \C$ the functors which send a diagram $c_1
\to \dots \to c_i$ to $c_1 \in \C$ resp. $c_i \in \C$. For any $E
\in \Fun(\C,\Ab)$, consider the complex
$$
\wt{P}_{i-1,\idot}(\C,E) = \tau_{i!}P_\idot(\C_i,\sigma_i^*E).
$$
Forgetting one vertex in a diagram gives a functor $\C_i \to
\C_{i-1}$, and this construction is strictly associative. Therefore
by the properties \thetag{i}, \thetag{ii} of the bar complex, we can
turn the collection $\wt{P}_{\idot,\idot}(\C,E)$ into a bicomplex,
with the second differential given by the same formula $\delta = d_1
- d_2 + \dots \pm d_i$ as in the case of $P_\idot(\C,E)$. The total
complex $\wt{P}_\idot(\C,E)$ is then also a resolution on the
functor $E$. To see this, one again evaluates at an object $c \in
\C$, and uses the same contracting homotopy $h$ as in the case of
$P_\idot(\C,E)$.

In the case of two functors $E,E' \in \Fun(\C,\Ab)$, we can use the
resolution $\wt{P}_\idot(\C,E)$ to compute $\RHom^\hdot(E,E')$;
this results in the triple complex with terms
\begin{equation}\label{rhom.bar.eq.1}
\RHom^{i,\hdot,\hdot}(E,E') =
C^\hdot(\C_i,\RHom^\hdot(\sigma_i^*E,\tau_i^*E')),
\end{equation}
a refinement of the double complex \eqref{rhom.bar.eq}. Its total
complex is functorially quasiisomorphic to $\RHom^\hdot(E,E')$.

\subsection{Semiorthogonal decompositions.}

We will also need some technology for working with triangulated
categories; the standard reference here is \cite{boka}.

In light of recent advances in axiomatic homotopy theory, it is
perhaps better to state explicitly that in this paper, our notion of
a {\em triangulated category} is the original notion of Verdier. A
full triangulated subcategory $\D \subset \D'$ in a triangulated
category $\D'$ is called {\em localizing} if the quotient $\D'/\D$
exists (in spite of the set-theoretic difficulties of the Verdier
construction). A full triangulated subcategory $\D \subset \D'$ is
called {\em left} resp. {\em right admissible} if the embedding
functor $\D \hookrightarrow \D'$ admits a left resp. right adjoint;
it is {\em admissible} if it is admissible both on the left and on
the right. The left orthogonal ${}^\perp\D \subset \D'$ consists of
objects $M \in \D'$ such that $\Hom(M,N) = 0$ for any $N \in
\D$. This is also a full triangulated subcategory in $\D'$, and it
known that $\D \subset \D'$ is left-admissible if and only if $\D$
and ${}^\perp\D$ generate the whole $\D'$. In this case, one says
that $\langle {}^\perp\D,\D \rangle$ is a {\em semi-orthogonal
  decomposition} of the triangulated category $\D'$. One shows that
$\D'$ is then generated by $\D$ and ${}^\perp\D$ in the following
strong sense: for any $M' \in \D'$, there exists a unique and
functorial distinguished triangle
\begin{equation}\label{tria}
\begin{CD}
{}^\perp M @>>> M' @>>> M @>>>
\end{CD}
\end{equation}
with $M \in \D$ and ${}^\perp M \in {}^\perp\D$. Moreover, the
category $\D \subset \D'$ is localizing, and we have a natural
identification ${}^\perp\D \cong \D'/\D$.

Analogously, the right orthogonal $\D^\perp \subset \D'$ consists of
objects $M \in \D'$ such that $\Hom(N,M)=0$ for any $N \in \D$, and
$\D \subset \D'$ is right-admissible if and only $\D'$ is generated
by $\D$ and $\D^\perp$; in this case, $\langle \D,\D^\perp \rangle$
is a semiorthogonal decomposition of the category $\D'$. We have the
following standard fact.

\begin{lemma}\label{orth}
Assume given a left-admissible triangulated subcategory $\D \subset
\D'$. Then the natural projection $\D^\perp \to \D'/\D \cong
{}^\perp\D$ is fully faithful, and it is an equivalence if and only
if $\D \subset \D'$ is right-admissible.\endproof
\end{lemma}

Given an admissible subcategory $\D \subset \D'$, one defines the
{\em gluing functor} $R:\D'/\D \to \D$ as the composition
$$
\begin{CD}
\D'/\D @>{\sim}>> \D^\perp \subset \D' @>>> \D,
\end{CD}
$$
where the second functor $\D' \to \D$ is left-adjoint to the
embedding $\D \subset \D'$. This is a triangulated functor defined
up to a canonical isomorphism. Objects $M' \in \D'$ are in natural
one-to-one correspondence with triples $\langle M^\perp,M,r \rangle$
of an object $M^\perp \in \D^\perp$, an object $M \in \D$, and a
gluing map $r:R(M^\perp) \to M[1]$.

We note that it is {\em not} possible to recover the category $\D'$
from $\D$, $\D'/\D$ and the gluing functor $R:\D'/\D \to \D$ (we can
recover objects, but not morphisms). However, there is the following
useful fact.

\begin{lemma}\label{orth.equi}
Assume given triangulated categories $\D'_1$, $\D'_2$ equipped with
left-admissible subcategories $\D_1 \subset \D'_1$, $\D_2 \subset
\D'_2$, and a triangulated functor $F:\D'_1 \to \D_2'$. Moreover,
assume that
\begin{enumerate}
\item $F$ sends $\D_1$ into $\D_2$, ${}^\perp\D_1$ into ${}^\perp\D_2$,
and the induced functors $F:\D_1 \to \D_2$, $F:\D_1'/\D_1 \to
\D_2'/\D_2$ are equivalences, and
\item $\D_1 \subset \D_1'$ is also right-admissible, and $F$ sends
  $\D_1^\perp$ into $\D_2^\perp \subset \D_2'$.
\end{enumerate}
Then $F$ is also an equivalence.
\end{lemma}

\proof{} Since $\D'_2$ is generated by $\D_2 = F(\D_1)$ and
${}^\perp\D_2 \cong F({}^\perp\D_1)$, the functor $F$ is essentially
surjective, and it suffices to prove that it is fully faithful ---
that is, for any $M,N \in \D_1'$, the map
$$
F:\Hom(M,N) \to \Hom(F(M),F(N))
$$
is an isomorphism. By \eqref{tria}, we may assume that $M$ lies
either in ${}^\perp\D_1$ or in $\D_1$. If $M \in {}^\perp\D_1$ and
$N \in \D_1$, then both sides are $0$. If $M,N \in {}^\perp\D_1$,
then the map is bijective by \thetag{i}. Thus we may assume $M \in
\D_1$. Decomposing $N$ by \eqref{tria} with respect to the
semiorthogonal decomposition $\langle\D_1,\D^\perp_1\rangle$, we see
that we may assume that either $N \in \D_1$ or $N \in
\D_1^\perp$. Then in the first case, the claim follows from
\thetag{i}, and in the second case, from \thetag{ii}.
\endproof

\subsection{$A_\infty$-structures.}\label{dgcat.subs}

To construct triangulated categories, we will use the machinery of
$A_\infty$-algebras and $A_\infty$-categories (this is very well
covered in the literature; a standard reference is, for example,
B. Keller's overview \cite{ke.a.inf}). We briefly recall the
relevant notions.

\subsubsection{Algebras and modules.}

An {\em $A_\infty$-algebra structure} on a graded free $\Z$-module
$A_\idot$ is given by a coderivation $\delta$ of the free non-unital
associative coalgebra $T_\idot(A_\idot[1])$ generated by $A_\idot$
shifted by $1$ such that $\delta^2=0$. Explicitly, the structure is
given by a collection of operations $b_n:A^{\hdot \otimes n} \to
A_\idot$, $n \geq 1$, and $\delta^2=0$ is equivalent to
\begin{equation}\label{a.inf.rel}
\sum_{i+j+l=n} b_{i+1+l} \circ (\id^{\otimes i} \otimes b_j \otimes
\id^{\otimes l}) = 0
\end{equation}
for any $n \geq 1$. For $n=1$, this reads as $b_1^2=0$, so that
$b_1$ is a differential which turns $A_\idot$ into a complex of
$\Z$-modules. After adding some signs depending on degrees of the
operands, the higher operations $b_n$, $n \geq 2$, can be arranged
together into an operad $\wt{\Ass}_\infty$ of complexes of
$\Z$-modules. Moreover, the action of symmetric groups on the
component complexes $\wt{\Ass}_\infty$ is irrelevant for the
definition of an $A_\infty$-algebra -- the operad $\wt{\Ass}_\infty$
is induced from an asymmetric operad $\Ass_\infty$ in the sense of
\cite{hinich}. The asymmetric operad $\Ass_\infty$ is equipped with
a canonical surjective augmentation quasiisomorphism $\Ass_\infty
\to \Ass$ onto the associative asymmetric operad $\Ass$. If one
forgets the differentials, $\Ass_\infty$ is the free asymmetric
operad generated by a single operation $b_n$ for each $n \geq
2$. Thus $\Ass_\infty$ is cofibrant in the natural closed model
structure on the category of asymmetric operads (see \cite{hinich}),
and the augmentation map $\Ass_\infty \to \Ass$ is a cofibrant
replacement for $\Ass$.

For any $A_\infty$-algebra $A_\idot$, the operation $m_2$ given by
\begin{equation}\label{m.2.b.2}
m_2(x,y) = (-1)^{\deg(x)}b_2(x,y)
\end{equation}
induces an associative multiplication on the homology groups
$H_\idot(A_\idot)$.  A {\em homological unit} in $A_\idot$ is an
element $1 \in A_0$ such that $b_1(1)=0$, and the cohomology class
of $1$ is the unit for the associative algebra $H_\idot(A_\idot)$. A
homological unit $1$ induces a contracting homotopy $h$ for the
differential $\delta$ on $T_\idot(A_\idot[1])$ by setting $h(a_1
\otimes \dots \otimes a_n) = 1 \otimes a_1 \otimes \dots \otimes
a_n$. We will assume that all $A_\infty$-algebras are equipped with
a homological unit.

An {\em $A_\infty$-morphism} $f$ between $A_\infty$-algebras
$A_\idot$, $A'_\idot$ is a DG coalgebra morphism
$T_\idot(A_\idot[1]) \to T_\idot(A'_\idot[1])$. Explicitly, $f$ is
given by a collection of maps $f_n:A^{\otimes n}_\idot \to
A'_\idot$, $n \geq 1$ such that
\begin{equation}\label{a.inf.map}
\sum_{i+j+l=n} f_{i+1+l} \circ (\id^{\otimes i} \otimes b_j \otimes
\id^{\otimes l}) = \sum_{i_1 + \dots + i_s} b_s \circ (f_{i_1}
\otimes \dots \otimes f_{i_s})
\end{equation}
for any $n \geq 1$. In particular, $f_1$ is a map of complexes, and
it induces an algebra map $H_\idot(f_1):H_\idot(A_\idot) \to
H_\idot(A'_\idot)$. The map $f$ is {\em unital} if so is
$H_\idot(f_1)$.

Every DG algebra is automatically an $A_\infty$-algebra (with
trivial $b_n$, $n \geq 3$). In particular, for a complex $M_\idot$
of objects in an abelian category $\Ab$ as in
Subsection~\ref{gen.subs}, $\End(M_\idot)$ is a DG algebra. The
structure of an {\em $A_\infty$ module over $A_\idot$} on the
complex $M_\idot$ is given by a unital $A_\infty$-morphism $A_\idot
\to \End(M_\idot)$. Explicitly, this given by a collection of maps
$$
b_n:A_\idot^{\otimes n-1} \otimes M_\idot \to M_\idot
$$
for all $n \geq 2$ satisfying \eqref{a.inf.rel} (where $b_1$ is the
differential on $M_\idot$). Equivalently, an $A_\infty$-module
structure on $M_\idot$ is given by a differential
$\delta$ on the cofree $T_\idot(A_\idot[1])$-comodule
$T_\idot(A_\idot[1]) \otimes M_\idot$ which turns it into a DG
comodule.

\subsubsection{Homotopy.}

The {\em homotopy category} $\Ho(A_\idot,\Ab)$ of $A_\infty$ modules
over $A_\idot$ is the full subcategory in the chain-homotopy
category of DG comodules over $T_\idot(A_\idot[1])$ spanned by
$A_\infty$-modules (objects are complexes $M_\idot$ equipped with an
$A_\infty$-module structure, maps are chain-homotopy classes of maps
between the corresponding DG comodules
$T_\idot(A_\idot[1])$. Explicitly, if we are given two complexes
$M_\idot$, $M'_\idot$ equipped with $A_\infty$-module structures
over $A_\idot$, then the graded group
$\Hom^\hdot_{A_\idot}(M_\idot,M'_\idot)$ of maps between the
corresponding $T_\idot(A_\idot[1])$-comodules can be canonically
written down as
\begin{equation}\label{alg.bar.hom}
\Hom^\hdot_{A_\idot}(M_\idot,M'_\idot)= \prod_{n \geq
  0}\Hom^{\hdot-n}(A_\idot^{\otimes n} \otimes M_\idot,M'_\idot),
\end{equation}
and it has a natural differential given by $d(a) = \delta \circ a -
a \circ \delta$. Maps in $\Ho(A_\idot,\Ab)$ are the degree-$0$
homology classes of this complex. The category $\Ho(A_\idot,\Ab)$ is
obviously triangulated. Inside it, we have the full triangulated
subcategory spanned by those $M_\idot$ which are acyclic as
complexes of objects in $\Ab$.

\begin{lemma}\label{a.inf.loc}
The subcategory of acyclic $A_\infty$-modules in $\Ho(A_\idot,\Ab)$
is localizing.
\end{lemma}

\proof{} As in the case of unbounded complexes of objects studied in
\cite{kel}, say that an $A_\infty$-module $M_\idot$ is {\em
$h$-injective} if it is right-orthogonal in $\Ho(A_\idot,\Ab)$ to
all acyclic $A_\infty$-modules. Then it suffices to prove that for
every $M_\idot \in \Ho(A_\infty,\Ab)$, there exists an $h$-injective
$\wt{M}_\idot \in \Ho(A_\idot,\Ab)$ equipped with a quasiisomorphism
$\wt{M}_\idot \to \wt{M}_\idot$. Choose a complex $\wt{M}_\idot$ of
objects in $\Ab$ which is $h$-injective and equipped with an
injective quasiisomorphism $M_\idot \to \wt{M}_\idot$. Then since we
assume that $A_\idot$ is a complexes of {\em free} $\Z$-modules, the
map $A_\idot^{\otimes n} \otimes M_\idot \to A^{\otimes n}_\idot
\otimes \wt{M}_\idot$ is an injective quasiisomorphism for any
$n$. Then we can solve the equations \eqref{a.inf.rel} by induction
on $n$, to obtain an $A_\infty$-module structure on $\wt{M}_\idot$
and an $A_\infty$-quasiisomorphism $M_\idot \to
\wt{M}_\idot$. Moreover, for any acyclic $A_\infty$-module
$N_\idot$, the terms in the complex
$\Hom_{A_\idot}(N_\idot,\wt{M}_\idot)$ can be rewritten as
\begin{equation}\label{a.to.the.left}
\Hom^\hdot(A_\idot^{\otimes n} \otimes N_\idot,M_\idot) \cong
\Hom^\hdot(A_\idot^{\otimes n},\Hom^\hdot(N_\idot,\wt{M}_\idot)),
\end{equation}
and since $A_\idot$ is a complex of free $\Z$-modules and
$\wt{M}_\idot$ is right-orthogonal to $N_\idot$, these are acyclic
complexes. Thus $\Hom_{A_\idot}(N_\idot,\wt{M}_\idot)$ is the
limit of an inverse system of acyclic complexes of abelian groups,
and the transition maps in this system are surjective.
Therefore the inverse limit is also acyclic, and
$\wt{M}_\idot$ is $h$-injective in $\Ho(A_\idot,\Ab)$.
\endproof

\begin{defn}
The {\em derived category} $\D(A_\idot,\Ab)$ is obtained by
localizing the category $\Ho(A_\idot,\Ab)$ with respect to the
subcategory of acyclic $A_\infty$-modules.
\end{defn}

We have the obvious forgetful functor $\D(A_\idot,\Ab) \to
\D(\Ab)$. It has both a left and a right-adjoint, the {\em free} and
the {\em cofree} module functors; they send a complex $M_\idot$ into
$A_\infty$-comodules given by $A_\idot \otimes M_\idot$,
resp. $\Hom_{\Z}(A_\idot,M_\idot)$, with an $A_\infty$-module
structure induced by the structure maps $b_n$ of the
$A_\infty$-algebra $A$. To see the adjunction, one notes that the
contracting homotopy $h$ given by the homological unit in $A_\idot$
induces a homotopy which contracts
$$
\Hom_{A_\idot}(A_\idot \otimes M_\idot,N_\idot)
$$
to $\Hom(M_\idot,N_\idot)$ for any $A_\infty$-module $N_\idot$, and
similarly for the cofree module $\Hom_\Z(A_\idot,N_\idot)$. In fact,
we have
\begin{equation}\label{fr.cofr}
\Hom^\hdot(A_\idot^{\otimes n} \otimes A_\idot \otimes M_\idot,N_\idot)
\cong \Hom^\hdot(A_\idot^{\otimes n},\Hom^\hdot_\Z(A_\idot,N_\idot))
\end{equation}
for any $n \geq 0$ and any two complexes $M_\idot$, $N_\idot$ in
$\Ab$, so that the complexes $\Hom^\hdot_{A_\idot}(-,-)$ one has to
contract are exactly the same in the free and in the cofree
case. The adjunction also holds on the level of homotopy
categories. In particular, for any $A_\infty$-module $M_\idot$, the
adjunction gives an $A_\infty$-map
$$
A_\idot \otimes M_\idot \to M_\idot.
$$
Iterating this construction, we obtain a version of the bar
resolution for $A_\infty$-modules; in effect, any $A_\infty$-module
$M_\idot$ is quasiisomorphic to the direct limit
$$
\lim_{\overset{n}{\to}} M^{(n)}_\idot
$$
so that the transition maps $M^{(n)}_\idot \to M^{(n+1)}_\idot$ are
injective, and their cokernels are free $A_\infty$-modules. Dually, we
have the cobar resolution, and every $M_\idot$ can represented as an
inverse limit of a system with surjective transition maps with
cofree kernels.

More generally, for any $A_\infty$-map $f:A_\idot \to A'_\idot$
between $A_\infty$-algebras, we have an obvious restriction functor
$f^*:\D(A'_\idot,\Ab) \to \D(A_\idot,\Ab)$. Replacing an
$A_\infty$-module with its free resp. cofree resolution, one easily
shows that $f^*$ has both a left-adjoint $f_!$ and a right-adjoint
$f_*$.

\begin{lemma}\label{a.inf.equi}
If the $A_\infty$-map $f:A_\infty \to A'_\infty$ is a
quasiisomorphism, then $f^*$ and $f_!$ are mutually inverse
equivalences of categories.
\end{lemma}

\proof{} Let $M_\idot$, $N_\idot$ be two $A_\infty$-modules over
$A'_\infty$. Rewrite the terms of the complex \eqref{alg.bar.hom} as
in \eqref{a.to.the.left}. Then for any $n \geq 0$, the natural map
$$
f^*:\Hom^\hdot(A_\idot^{\otimes n},\Hom^\hdot(N_\idot,\wt{M}_\idot))
\to \Hom^\hdot(A_\idot^{'\otimes n},\Hom^\hdot(N_\idot,\wt{M}_\idot))
$$
is a quasiisomorphism. This implies that $f^*$ is fully
faithful. To see the it is essentially surjective, note that it
commutes with filtered direct limits, and for any complex $M_\idot$,
$f^*(A'_\idot \otimes M_\idot)$ is quasiisomorphic to the free
$A_\infty$-module $A_\idot \otimes M_\idot$.
\endproof

\subsubsection{Coalgebras.}\label{coa.subsu}

An {\em $A_\infty$-coalgebra} is an $A_\infty$-algebra in the
category opposite to that of abelian groups. Explicitly, an
$A_\infty$-coalgebra structure on a graded $\Z$-module $A_\idot$ is
given by a collection of operations
$$
b_n:A_\idot \to A_\idot^{\otimes n}, \qquad n \geq 1,
$$
subject to relations \eqref{a.inf.rel} (where the composition
$\circ$ should be understood in the reverse order). Equivalently,
this structure is encoded by a collection of maps
$$
A_\idot \otimes \Ass_\infty \to A_\idot^{\otimes n},
$$
where $\Ass_\infty$ is the asymmetric $A_\infty$-operad. As in the
algebra case, we will only consider $A_\infty$-coalgebras which are
free as $\Z$-modules. For any $A_\infty$-coalgebra $A_\idot$, $b_1$
is a differential, that is, $b_1^2=0$, so that $A_\idot$ becomes a
complex of free abelian groups. The dual complex $(A_\idot)^* =
\Hom_\Z(A_\idot,\Z)$ is naturally an $A_\infty$-algebra. A {\em
  homological counit} for an $A_\infty$-coalgebra $A_\idot$ is a map
$1:A_0 \to \Z$ such that $1 \circ b_1 = 0$, and $1 \in (A_0)^*$ is a
homological unit for $(A_\idot)^*$. We will only consider counital
coalgebras.

An {\em $A_\infty$-morphism} $f:A'_\idot \to A_\idot$ of
$A_\infty$-coalgebras is given by a collection of maps
$$
f_n:A'_\idot \to A_\idot^{\otimes n}, \qquad n \geq 1
$$
satisfying \eqref{a.inf.map}. Such a map is a quasiisomorphism if so
is its component $f_1$. The dual maps $f_n^*$ give an $A_\infty$-map
$f^*:(A_\idot)^* \to (A'_\idot)^*$; $f$ is {\em counital} if $f^*$
is unital. As in the algebra case, we will only consider counital
$A_\infty$-maps.

An {\em $A_\infty$-comodule} over an $A_\infty$-coalgebra $A_\idot$
in an abelian category $\Ab$ is given by a graded object $M_\idot$
in $\Ab$ together with maps
\begin{equation}\label{como.str}
b_n:M_\idot \to M_\idot \otimes A_\idot^{\otimes n-1},\qquad n \geq
1,
\end{equation}
again subject to \eqref{a.inf.rel}. Again, $b_1$ is a differential
on $M_\idot$. Moreover, $M_\idot$ is automatically an
$A_\infty$-module over $(A_\idot)^*$, and $H_\idot(M_\idot)$ is a
module over the cohomology algebra $H^\hdot((A_\idot)^*)$. If this
module is unital, the $A_\infty$-comodule $M_\idot$ is called
counital. We will only consider counital comodules.

Equivalently, an $A_\infty$-structure on $A_\idot$ can be described
as a square-zero derivation of the completed tensor algebra
$\wh{T}^\hdot(A_\idot[-1])$, and an $A_\infty$-comodule structure on
a complex $M_\idot$ is the same as a differential on the completed
tensor product
\begin{equation}\label{wh.M}
\wt{T}^\hdot(A_\idot[-1]) \wh{\otimes} M_\idot =
\lim_{\overset{n}{\gets}} T^\hdot(A_\idot[-1])/T^{\geq
    n}(A_\idot[-1]) \otimes M_\idot
\end{equation}
which turns it into a DG $\wh{T}^\hdot(A_\idot[-1])$-module. The
homotopy category $\Ho(A_\idot,\Ab)$ of $A_\infty$-comodules over
$A_\idot$ in $\Ab$ is the full subcategory in the chain-homotopy
category of topological DG $\wh{T}^\hdot(A_\idot[-1])$-modules
spanned by DG modules of the form \eqref{wh.M}. Given two
$A_\infty$-comodules $M_\idot$, $M'_\idot$, we define a complex
$\Hom^\hdot_{A_\idot}(M_\idot,M'_\idot)$ as
$$
\prod_{n \geq 0}\Hom^\hdot(M_\idot,A^{\otimes n} \otimes M'_\idot),
$$
with the differential dual to that of \eqref{alg.bar.hom}. Then the
space of maps from $M_\idot$ to $M'_\idot$ in $\Ho(A_\idot,\Ab)$ is
the $0$-th homology group of this complex. An object $M_\idot \in
\Ho(A_\idot,\Ab)$ is called {\em acyclic} if it is acyclic as a
complex in $\Ab$.

\begin{lemma}\label{weib}
  Assume that $\Ab$ is the category of modules over a ring $R$. Then
  the subcategory of acyclic complexes in $\Ho(A_\idot,\Ab)$ is
  localizing.
\end{lemma}

\proof{} We adopt the method of \cite[Proposition 10.4.4]{We}; I am
grateful to the referee for suggesting this reference.

As in \cite[Proposition 10.4.4]{We}, it suffices to prove that for
any $A_\infty$-comodule $M_\idot \in \Ho(A_\idot,\Ab)$, there exists
a set $A$ of quasiisomorphisms $r_\alpha:M^\alpha_\idot \to
M_\idot$, $\alpha \in A$, such that for any quasiisomorphism
$r:M'_\idot \to M_\idot$, one of the maps $r_\alpha$ factors through
$r$. Let $\kappa$ be an infinite cardinal larger than the
cardinality of $\bigoplus M_r$. There is at most a set of
quasiisomorphisms $r_\alpha:M^\alpha_\idot \to M_\idot$ such that
the cardinality of $\bigoplus M^\alpha_i$ is at most $\kappa$; take
them all. Assume given a quasiisomorphism $r:M'_\idot \to M_\idot$.

Every element $m \in M'_i$, $i \in \Z$, lies in an at most countable
$A_\infty$-subcomodule $M^m_\idot \subset M'_i$ -- indeed, we can
take the abelian subgroup generated by $m$, add to it all the
left-hand sides of the structure maps $b_n$ of \eqref{como.str}, and
repeat the procedure by induction. Therefor there exists an
$A_\infty$-subcomodule $M^{(1)}_\idot \subset M'_\idot$ of
cardinality at most $\kappa$ such that the map
$$
r^{(1)}_i:H_i(M^{(1)}_\idot) \to H_i(M_\idot)
$$
induced by $r$ is surjective for every $i \in \Z$. Repeating the
procedure, we obtain a system of subcomodules $M^{(n)}_\idot \subset
M'_\idot$, $n \geq 1$, such that $M^{(n+1)}_\idot$ contains
$M^{(n)}_\idot$, the map
$$
r^{(n)}_i:H_i(M^{(n)}_\idot) \to H_i(M_\idot)
$$
is still surjective for any $n \geq 2$, $i \in \Z$, and the natural
map
$$
H_i(M^{(n)}_\idot) \to H_i(M^{(n+1)}_\idot)
$$
induced by the embedding $M^{(n)}_\idot \to M^{(n+1)}_\idot$
annihilates $\Ker r^{(n)}_i$ for any $i \in \Z$. Let $M''_\idot =
\bigcup M^{(n)}_\idot \subset M'_\idot$. Then the natural map
$$
r'': M''_\idot \to M_\idot
$$
is a quasiisomorphism by construction, and it is of the form
$r^\alpha$ for some $\alpha$ in the indexing set $A$.
\endproof 

\begin{remark}
  The proof can probably be modified so that it only requires our
  original assumptions on $\Ab$, but I haven't pursued it for lack
  of interesting examples.
\end{remark}

\begin{defn}
The {\em derived category} $\D(A_\idot,\Ab)$ of $A_\infty$-comodules
over $A_\idot$ in $\Ab$ is the quotient of the homotopy category
$\Ho(A_\idot,\Ab)$ by the subcategory of acyclic objects.
\end{defn}

As in the algebra case, an $A_\infty$-comodules $M_\idot \in
\Ho(A_\idot,\Ab)$ is called {\em $h$-injective} if it
right-orthogonal to all acyclic objects. We have an obvious
forgetful functor $\Ho(A_\idot,\Ab) \to \Ho(\Ab)$ onto the
chain-homotopy category of complexes in $\Ab$, and it has a
right-adjoint which sends $M_\idot \in \Ho(\Ab)$ into the cofree
comodule $A_\idot \otimes M_\idot$, with the comodule structure maps
$b_n$ given by the structure maps of $A_\infty$. As in the algebra
case, the counit on $A_\idot$ induces a homotopy which contracts the
complex
$$
\Hom^\hdot_{A_\idot}(N_\idot,A_\idot \otimes M_\idot)
$$
onto $\Hom^\hdot(N_\idot,M_\idot)$ for any $N_\idot,M_\idot \in
\Ho(A_\idot,\Ab)$, and this gives the adjunction. In particular, if
$M_\idot \in \D(\Ab)$ is $h$-injective, the cofree
$A_\infty$-comodule $A_\idot \otimes M_\idot$ is $h$-injective in
$\Ho(A_\idot,\Ab)$. One is tempted now to dualize the
bar-construction and obtain an $h$-injective replacement for any
$M_\idot \Ho(A_\idot,\Ab)$. However, {\em this does not work}. The
reason is the following: to be $h$-injective, the cobar resolution
of an $A_\infty$-comodule $M_\idot$ has to be a projective limit of
$h$-injective comodules. Thus as a graded object in $\Ab$, it is of
the form
$$
A^\hdot \otimes \prod_{n \geq 0} (A^{\otimes n}_\idot \otimes
M_\idot).
$$
However, this is different from
$$
\prod_{n \geq 0} (A_\idot \otimes A^{\otimes n}_\idot \otimes M_\idot),
$$
and it is the latter, not the former which can be contracted onto
$M_\idot$ by the homological counit of $A_\infty$. This is why the
existence of $\D(A_\idot,\Ab)$ has to be proved by an indirect
method. And even when $\Ab$ is as in Lemma~\ref{weib}, it is not
clear at present whether for an arbitrary $A_\idot$, any $M_\idot
\in \Ho(A_\idot,\Ab)$ is quasiisomorphic to an $h$-injective
$M'_\idot$.

\begin{remark}
  If $A$ is simply a coalgebra, not an $A_\infty$ or DG coalgebra,
  then more is known, since a very comprehensive study of the
  homological properties of unbounded complexes of comodules has
  been done recently by L. Positelski \cite{lenya}. In particular,
  it has been proved in \cite{lenya} that $h$-injective replacements
  do exist. However, the proof is very indirect, and it is not clear
  at all whether it can be generalized to the $A_\infty$-case.
\end{remark}

As in the algebra case, an $A_\infty$-map $f:A_\idot \to A'_\idot$
induces a natural corestriction functor $f^*:\D(A_\idot,\Ab) \to
\D(A'_\idot,\Ab)$, assuming that $\Ab$ is as in Lemma~\ref{weib} so
that both categories are well-defined. If $A'_\idot = \Z$ and $f$ is
the counit, this is the forgetful functor $\D(A_\idot,\Ab) \to
\D(\Ab)$, and it has a right-adjoint given by the cofree comodule
construction. In general, it is not clear whether $f^*$ admits a
right-adjoint (and the situation with the left-adjoint is even
worse).

Let us list some other things that do not work for coalgebras.
\begin{enumerate}
\item There are no free comodules, only the cofree ones
  (\eqref{fr.cofr} does not work, since it would involve double
  dualization).
\item The proof of Lemma~\ref{a.inf.loc} breaks down because an
  analog of \eqref{a.to.the.left} is not an isomorphism.
\item For those interested in such things, the category of
  $A_\infty$-comodules for a general $A_\idot$ does not admit a
  closed model structure (or at least, none such is known).
\item Finally, Lemma~\ref{a.inf.equi} fails. Not only does its proof
  break down, the statement itself is false. In fact, one case where
  this happens will be the main subject of Section~\ref{galois.sec}.
\end{enumerate}

\subsubsection{Categories.}

Informally, a (small) $A_\infty$-category is an $A_\infty$-algebra
``with many objects''. To keep track of the combinatorics, it is
convenient to use the approach of \cite{lef}. For any set $S$,
denote by $\Z_S\amod$ the category of $S$-graded $\Z$-modules. Fix a
set $S$, and consider the category $\Z_{S \times S}\amod$. Equip it
with a tensor product by setting
$$
(A' \otimes A'')_{s',s''} = \bigoplus_{s \in S}A'_{s',s} \otimes
A''_{s,s''}
$$
for any $A',A'' \in \Z_{S \times S}\amod$. This is not symmetric, but
neither are the operads $\Ass$ and $\Ass_\infty$, so that speaking
about associative and $A_\infty$-algebras in $\Z_{S \times S}\amod$
makes perfect sense. For any small additive category $\B$ with the
set of objects $S$, the sum
$$
\B^S = \bigoplus_{s,s'}\B(s,s')
$$
is an associative algebra in $\Z_{S \times S}\amod$. Then an {\em
  $A_\infty$-category} $\B_\idot$ consists of
\begin{enumerate}
\item a small graded additive category $\overline{\B}_\idot$, with a
  certain set of objects $S$, and
\item an $A_\infty$-algebra $\B_\idot$ in $\Z_{S \times S}\amod$
  equipped with an isomorphism
$$
H_\idot(\B_\idot) \cong \overline{\B}_\idot^S
$$
of graded associative algebras in $\Z_{S \times S}\amod$.
\end{enumerate}
Our reason for making this slightly convoluted definition is that it
automatically takes care of the units. As in the algebra case, we
will assume that $\B_\idot(b,b')$ is a complex of {\em free}
$\Z$-modules for any two objects $b,b' \in \B_\idot$.

For any map of sets $f:S \to S'$, we have an obvious pseudotensor
restriction functor $f^*:\Z_{S' \times S'}\amod \to \Z_{S \times
  S}\amod$; an {\em $A_\infty$-functor} between
$A_\infty$-categories $\B_\idot$, $\B'_\idot$ joists of a functor
$f:H_\idot(\B_\idot) \to H_\idot(\B'_\idot)$ and an
$A_\infty$-morphism $\B_\idot \to f^*\B'_\idot$.

For any abelian category $\Ab$, the category $\Ab_S$ of $S$-graded
objects in $\Ab$ is naturally a module category over the tensor
category $\Z_{S \times S}\amod$. An {\em $A_\infty$-functor} from an
$A_\infty$-category $\B_\idot$ to $\Ab$ is then an $A_\infty$-module
over $\B_\idot$ in $\C_S$ such that units act by identity maps on
the corresponding $\overline{\B}_\idot$-module
$H_\idot(M_\idot)$. As in the algebra case, we have the homotopy
category $\Ho(\B_\idot,\Ab)$ and the derived category
$\D(\B_\idot,\Ab)$ of $A_\infty$-functors from $\B_\idot$ to
$\Ab$. The role of free and cofree modules is played by
representable and corepresentable functors: for any object $b \in
\B_\idot$ and any $M_\idot \in \Ho(\Ab)$, these are given by
$$
M_\idot^b(b') = M_\idot \otimes \B_\idot(b,b'), \qquad M_{\idot
  b}(b') = \Hom_{\Z}(\B_\idot(b',b),M_\idot)
$$
for any $b' \in \B_\idot$. We also have the bar and cobar
resolution, so that the derived category $\D(\B_\idot,\Ab)$
is generated by representable resp. corepresentable functors in the
same sense as the category of $A_\infty$-modules is generated by
free resp. cofree modules. For any $A_\infty$-functor $f:\B_\idot
\to \B'_\idot$ between two $A_\infty$-categories, we have the
restriction functor $f^*:\D(\B'_\idot,\Ab) \to \D(\B_\idot,\Ab)$ and
its two adjoints
$$
f_!:\D(\B_\idot,\Ab) \to \D(\B'_\idot,\Ab), \qquad
f_*:\D(\B_\idot,\Ab) \to \D(\B'_\idot,\Ab).
$$
An $A_\infty$-functor $f:\B_\idot \to \B'_\idot$ is a {\em
  quasiequivalence} if the corresponding functor $H_\idot(\B_\idot)
\to H_\idot(\B'_\idot)$ is an equivalence, and the natural map
$$
f_1:\B_\idot(b,b') \to \B'_\idot(f(b),f(b'))
$$
is a quasiisomorphism for any $b,b' \in \B_\idot$. It is not
difficult to generalize Lemma~\ref{a.inf.equi} and show that for a
quasiequivalence $f$, the functors $f^*$ and $f_!$ are mutually
inverse equivalences of the derived categories.

We note that any small category $\C$ defines an additive category
$\Z[\C]$ with the same objects, and morphisms given by
$\Z[\C](c,c') = \Z[\C(c,c')]$, where $\Z[S]$ for a set $S$ means the
free abelian group generated by $S$. Then $\Z[\C]$ can be treated as
an $A_\infty$-category (placed in homological degree $0$), and we of
course have $\D(\C,\Ab) \cong \D(\Z[\C],\Ab)$.

\medskip

We will also need a version of this for coalgebras. In fact, it will
more convenient to use a slightly more refined notion. Assume given
a small category $\C$, with the set of objects $\C^{ob}$ and the set
of morphisms $\C^{mor}$. Introduce a tensor product on the category
$\Z_{\C^{mor}}\amod$ by setting
$$
(A' \otimes A'')_f = \bigoplus_{f',f'' \in \C^{mor},f' \circ f'' =
  f}A'_{f'} \otimes A''_{f''}
$$
for any $A',A'' \in \Z_{\C^{mor}}\amod$ and $f \in \C^{mor}$. Then
$\Z_{\C^{mor}}\amod$ is a monoidal category, and for any abelian
category $\Ab$ as in Subsection~\ref{gen.subs}, $\Ab_{\C^{ob}}$
is a module category over $\Z_{\C^{ob}}\amod$.

\begin{defn}\label{DG.coa}
A {\em $\C$-graded $A_\infty$-coalgebra} $\A_\idot$ is an
$A_\infty$-coalgebra in the monoidal category
$\Z_{\C^{mor}}\amod$. An {\em $A_\infty$-comodule} in $\Ab$ over
$A_\idot$ is an $A_\infty$-module over $A_\idot$ in the module
category $\Ab_{\C^{mor}}$.
\end{defn}

As in the $A_\infty$-category case, we will always assume that our
graded coalgebras consist of free $\Z$-modules; we will also assume
that coalgebras and comodules are (homologically) counital in the
same sense as in the non-graded case. Explicitly, a $\C$-graded
$A_\infty$-coalgebra $A_\idot$ is given by a collection
$\A_\idot(f)$ of complexes of free $\Z$-modules numbered by
morphisms $f$ of the category $\C$, together with a homological
counit map $\A_\idot(\id_c) \to \Z$ for any object $c \in \C$ and a
comultiplication map
$$
b_n:\A_\idot(f_1 \circ \dots \circ f_n) \to \A_\idot(f_0) \otimes
\dots \otimes A_\idot(f_n)
$$
for any $n \geq 2$ and any $n$-tuple of composable maps
$f_1,\dots,f_n$ in $\C$, subject to \eqref{a.inf.rel}.  An
$\Ab$-valued $\A_\idot$-comodule $E_\idot$ in the category $\Ab$ is
a collection of complexes $E_\idot(c)$ in $\Ab$, one for each object
$c \in \C$, together with maps
$$
b_n:E_\idot(c) \to \A_\idot(f_1) \otimes \dots \otimes A_\idot(f_n)
\otimes E_\idot(c')
$$
for any $n \geq 1$ and any $n$-tuple of composable maps
$f_0,\dots,f_n$, $f_0 \circ \dots \circ f_n:c \to c'$ in $\C$, again
subject to \eqref{a.inf.rel}.

As in the non-graded case, any $\C$-graded $A_\infty$-coalgebra
$A_\idot$ produces the homotopy and the derived categories of
$A_\idot$-comodules, denoted $\Ho(A_\idot,\Ab)$
resp. $\D(A_\idot,\Ab)$. For any object $c \in \C$ and any complex
$M_\idot$ in $\Ab$, we have the corepresentable $A_\idot$-comodule
$M_\idot^c$ given by
\begin{equation}\label{corep.coa}
M_\idot^c(c') = \prod_{f:c' \to c}A_\idot(f) \otimes M_\idot.
\end{equation}
For any $\C$-graded $A_\infty$-map $f:A_\idot \to A'_\idot$, we have
the restriction functor $f^*:\D(A_\idot,\Ab) \to
\D(A'_\idot,\Ab)$. As in the coalgebra case, it does not have to be
an equivalence even if $f$ is a quasiisomorphism in a suitable
sense. On the other hand, assume given another small category $\C'$
and a functor $\rho:\C' \to \C$. Define a $\C'$-graded
$A_\infty$-coalgebra $\rho^*A_\idot$ by setting
$$
\rho^*A_\idot(f) = A_\idot(\rho(f))
$$
for any morphism $f$ in $\C'$, with the same structure maps
$b_n$. We then have an obvious pullback functor
$$
\rho^*:\D(A_\idot,\Ab) \to \D(\rho^*A_\idot,\Ab).
$$
I do not know under what assumptions, if any, either of the functors
$\rho^*$, $f^*$ has a right adjoint.

\subsection{$2$-categories.}\label{2cat.subs}

To produce $A_\infty$-categories, we will use $2$-categories; we end
the preliminaries with a brief sketch of the corresponding
construction.

Assume given a small monoidal category $\C$, and consider the bar
complex $C_\idot(\C,\Z)$. If $\C$ is strictly associative, then the
tensor product functor $m:\C \times \C \to \C$ induces an
associative DG algebra structure on $C_\idot(\C,\Z)$ (apply the
properties \thetag{ii} and \thetag{iii} of
Subsection~\ref{bar.gene}). More generally, if we also have an
object $T \in \Fun(\C,\Z)$ and a map
\begin{equation}\label{T.t.map}
T \boxtimes T \to m^*T
\end{equation}
which is associative on triple products, then $C_\idot(T)$ becomes
an associative DG algebra (plug in the property \thetag{i}).

However, monoidal categories in nature are usually associative only
up to an isomorphism -- there is an associativity isomorphism $m
\circ m_{12} \cong m \circ m_{23}$ satisfying the pentagon equation.
We observe that in this case $C_\idot(\C,\Z)$ is no longer a DG
algebra, but it has an $A_\infty$-algebra structure.

Indeed, for any $n \geq 2$, let $I_n$ be the groupoid whose objects
are all possible $n$-ary operations obtained from a single binary
operation, and which has exactly one morphism between every two
objects. Then $I_n$, $n \geq 2$ form an asymmetric operad of
categories, and any weakly associative monoidal category $\C$ is an
algebra over this operad: we have natural functors
$$
I_n \times \C^n \to \C
$$
for every $n$. The bar complexes $C_\idot(I_n,\Z)$ form an
asymmetric operad of complexes of $\Z$-modules, $C_\idot(\C,\Z)$ is
an algebra over this operad, and the operad itself is a resolution
of the trivial asymmetric operad $\Ass$. The asymmetric operad
$\Ass_\infty$ is another such resolution, and it is
cofibrant. Therefore the augmentation map $\Ass_\infty \to \Ass$
factors through a map $\Ass_\infty \to C_\idot(I_\idot,\Z)$. Fixing
this map once and for all, we turn the bar complex $C_\idot(\C,\Z)$
for any monoidal category $\C$ into an $A_\infty$-algebra.

Analogously, a (weakly) monoidal functor between monoidal categories
$\C$, $\C'$ induces an $A_\infty$-map between the
$A_\infty$-algebras $C_\idot(\C,\Z)$, $C_\idot(\C',\Z)$.

The same construction obviously works for bar complex
$C_\idot(\C,T)$ with coefficients, where $T \in \Fun(T,\Z)$ is
equipped with an associative map \eqref{T.t.map}.

Moreover, if we have a $2$-category $\Q$ with a certain set of
objects $\{c\}$ and categories of $1$-morphisms $\Q(c,c')$, $c,c'
\in \C$, then the same construction produces an $A_\infty$-category
with the same objects, and with the bar complexes
$C_\idot(\Q(c,c'),\Z)$ as complexes of morphisms. This is also
functorial with respect to $2$-functors, and has an obvious version
with coefficients.

\section{Recollection on Mackey functors.}\label{mack.rec}

This ends the preliminaries. For the convenience of the reader, we
start the paper itself by briefly recalling the definitions and
known facts about Mackey functors (we more-or-less follow the
expositions in \cite{may2} and \cite{the}).

\subsection{Definitions.}\label{m.defn.subs}

Assume given a group $G$, and let $\Gamma_G$ be the category of
finite sets equipped with a $G$-action. This category obviously has
pullbacks. Define a bigger category $Q\Gamma_G$ as follows: objects
are the same as in $\Gamma_G$, maps from $S_1$ to $S_2$ are
isomorphism classes of diagrams $S_1 \gets S_1' \to S_2$,
composition of $S_1 \gets S_1' \to S_2$ and $S_2 \gets S_2' \to S_3$
is given by the diagram $S_1 \gets S_1' \times_{S_2} S_2' \to
S_3$. The category $Q\Gamma_G$ is self-dual, $Q\Gamma_G \cong
Q\Gamma_G^{opp}$. Moreover, disjoint unions of sets give finite
coproducts both in the category $\Gamma_G$ and in the category
$Q\Gamma_G$. Every $G$-finite set $S \in \Gamma_S$ can be
canonically decomposed into such a disjoint union
\begin{equation}\label{unio}
S = \coprod_{p \in S/G} S_i
\end{equation}
of subsets $S_p$ on which $G$ acts transitively; we call them {\em
$G$-orbits}. This decomposition is valid both in $\Gamma_G$ and in
$Q\Gamma_G$.

\begin{defn}\label{mackey}
A {\em $G$-Mackey functor} $M$ is a functor $M:Q\Gamma_G \to \Ab$ to
the category $\Ab$ of abelian groups which is additive in the
following sense: for any $S \in Q\Gamma_S$, the natural map
$$
\bigoplus_{p \in S/G}M(S_p) \to M(S)
$$
induced by the decomposition \eqref{unio} is an isomorphism.
\end{defn}

Mackey functors obviously form an abelian category which we denote
by $\M(G,\Ab)$, or simply by $\M(G)$. By definition, $\M(G,\Ab)$ is
a full subcategory in $\Fun(Q\Gamma_G,\Ab)$, and one checks easily
that the embedding $\M(G,\Ab) \to \Fun(Q\Gamma_G,\Ab)$ admits a left
adjoint, which we call {\em additivization} and denote by
$\Add:\Fun(Q\Gamma_G,\Ab) \to \M(G,\Ab)$ (in fact, $\Add$ is also
right-adjoint to the embedding).

For any cofinite subgroup $H \in G$, the value $M([G/H])$ of a
$G$-Mackey functor $M$ on the $G$-orbit $[G/H]$ is usually denoted
by $M^H$. By the additivity property, $M^H$ for all $H \subset G$
completely define $M$. Explicitly, a $G$-Mackey functor $M$ is given
by
\begin{enumerate}
\item an abelian group $M^H$ for any cofinite subgroup $H \in G$,
and
\item two maps $f_*:M^{H_1} \to M^{H_2}$, $f^*:M^{H_2} \to M^{H_1}$
  for any two cofinite subgroups $H_1, H_2 \subset G$ and a
  $G$-equivariant map $f:[G/H_1] \to [G/H_2]$,
\end{enumerate}
such that $f^* \circ g^* = (g \circ f)^*$ and $g_* \circ f_* = (g
\circ f)_*$ for any two composable maps $f:[G/H_1] \to [G/H_2]$,
$g:[G/H_2] \to [G/H_3]$, and for any two maps $f:[G/H_1] \to [G/H]$,
$g:[G/H_2] \to [G/H]$, we have
\begin{equation}\label{d.cos}
g^* \circ f_* = \sum_{p \in S/G}f_{p*} \circ
g_p^*,
\end{equation}
where we let $S = [G/H_1] \times [G/H_2]$, take its decomposition
\eqref{unio}, and let $g_p:S_p \to [G/H_1]$, $f_p:S_p \to [G/H_2]$ be
the natural projections. We note that since $S/G = G\backslash(G
\times G)/(H_1 \times H_2) \cong H_1\backslash G/H_2$, the
components $S_p$ correspond to double cosets $H_1gH_2 \subset G$;
for this reason, \eqref{d.cos} is known as the {\em double coset
formula}. This is the original definition of Mackey functors
introduced by Dress \cite{dress}; the version with the category
$Q\Gamma_G$ is due to Lindner \cite{lind}. The collection $\langle
f_*,f^*\rangle$ without the condition \eqref{d.cos} is sometimes
called a {\em bifunctor} (from the category of finite $G$-orbits to
$\Ab$).

\begin{exa}\label{r.exa}
Representation ring: setting $[G/H] \mapsto R_H$, the representation
ring of the group $H$, defines a Mackey functor, with $f^*$ given by
restriction and $f_*$ given by induction. This is the origin of the
notion and the name: the double coset formula for $R_H$ was found by
Mackey.
\end{exa}

\begin{exa}\label{h.exa}
Cohomology: setting $[G/H] \mapsto H^\hdot(H,\Z)$, $f^*$ given by
restriction, $f_*$ given by corestriction, defines a (graded) Mackey
functor.
\end{exa}

To obtain a third useful definition of Mackey functors, one
considers an {\em additive} category $\B^G$ defined as follows:
objects of $\B^G$ are finite $G$-orbits $[G/H]$, and the set of maps
from $S_1$ to $S_2$ is the free abelian group generated by
isomorphism classes of diagrams $S_1 \gets S \to S_2$, where $S$ is
another finite $G$-orbit. Composition $g \circ f$ of two maps
$f:[G/H_1] \to [G/H_2]$, $g:[G/H_2] \to [G/H_3]$ represented by
diagrams $[G/H_1] \gets S_f \to [G/H_2]$, $[G/H_2] \gets S_g \to
[G/H_3]$ is given by
$$
g \circ f = \sum_p (g \circ f)_p,
$$
where we consider the decomposition \eqref{unio} of the fibered
product $S = S_f \times_{[G/H_2]} S_g$, and let $(g \circ f)_p$ be the
map represented by the diagram $[G/H_1] \gets S_p \to [G/H_3]$. Then a
$G$-Mackey functor $M$ is obviously the same thing as an additive
functor $\B^G \to \Ab$. The category $\B^G$ can be described more
explicitly in terms of the so-called ``Burnside rings''.

\begin{defn}
The {\em Burnside ring} $\A^G$ of a group $G$ is the abelian group
generated by isomorphism classes $[S]$ of objects $S \in \Gamma_S$
in the category $\Gamma_G$, modulo the relations $[S_1] + [S_2] =
[S_1 \copr S_2]$, and with the product given by $[S_1] \cdot [S_2]
= [S_1 \times S_2]$.
\end{defn}

Then the endomorphism ring of the trivial $G$-orbit $[G/G] \in \B^G$
obviously coincides with the Burnside ring $\A^G$. And more
generally, given two finite orbits $S_1,S_2 \in \B^G$, we have
\begin{equation}\label{b.g.hom}
\B^G(S_1,S_2) = \bigoplus_{p \in (S_1 \times S_2)/G} \A^{H_p},
\end{equation}
where $S_p = [G/H_p]$ are the components in the decomposition
\eqref{unio} of the product $S = S_1 \times S_2$.

\begin{remark}
Normally the definitions in this subsection are given for a {\em
finite} group $G$; however, everything works in a slightly wider
generality, and this is sometimes useful. Of course, the category
$\M(G)$ as defined here only depends on the profinite completion of
the group $G$. There is also a version of Mackey functors for
topological groups such as finite-dimensional Lie groups, see
e.g. \cite{tD}; however, this is beyond the scope of the present
paper.
\end{remark}

\subsection{Functoriality.}\label{mack.func.subs}

For any cofinite subgroup $H \subset G$ of a group $G$ and a finite
$H$-set $S$, the product $S \times_H G = (S \times G)/H$ is
naturally a finite $G$-set, with the $G$-action through the second
factor. This defines a functor $\gamma_H^G:\Gamma_H \to
\Gamma_Q$. In fact, if we denote $S = [G/H]$, then $\gamma_H^G$ is an
equivalence between $\Gamma_H$ and the category $\Gamma_G/S$ of
finite $G$-sets equipped with a map to $S$. The functor $\gamma_H^G$
obviously commutes with fibered products, thus extends to a functor
$\gamma_H^G:Q\Gamma_H \to Q\Gamma_G$. It also commutes with disjoint
unions, so that we can define an exact functor
$$
\Restr^G_H:\M(G) \to \M(H)
$$
which sends $M:Q\Gamma_H \to \Ab$ to $\gamma^H_G \circ M:Q\Gamma_G
\to \Ab$. For any $G$-Mackey functor $M \in \M(G)$, the $H$-Mackey
functor $\Restr^G_H(M)$ is called the {\em restriction} of $M$ to $H
\subset G$.

Assume that a subgroup $H \subset G$ is normal, and let $N=G/H$ be
the quotient. Then any $N$-set is also a $G$-set, so that we have an
obvious full embedding $\Gamma_N \to \Gamma_G$ which induces a full
embedding $Q\Gamma_N \to Q\Gamma_G$ compatible with disjoint
unions. This induces a functor $\Psi^H:\M(G) \to \M(N)$ (usually
$\Psi^H(M)$ is denoted simply by $M^H$, but this might cause
confusion). However, we also have a functor $\Infl^N_G:\M(N) \to
\M(G)$ called {\em inflation} and given by ``extension by $0$'': we
set
$$
\Infl^N_G(M)^K =
\begin{cases}
M^{K/H}, &\quad H \subset K,\\
0, &\quad\text{otherwise}.
\end{cases}
$$
This is an exact full embedding. It has a left-adjoint which is
denoted by $\Phi^H:\M(G) \to M(N)$.

If the subgroup $H \subset G$ is not normal, we can consider its
normalizer $N_H \subset G$.  Assume that the normalizer $N_H \subset
G$ is cofinite in $G$. Then we can define the functor $\Phi^H:\M(G)
\to \M(N_H/H)$ by first restricting to $N_G(H)$:
$$
\Phi^H(M) = \Phi^H(\Restr^G_{N_H}(M)).
$$
Analogously, $M^H$ has a natural structure of a $(N_H/H)$-Mackey
functor.

\subsection{Products.}\label{mack.prod.subs}

Both in Example~\ref{r.exa} and Example~\ref{h.exa}, the Mackey
functors have an additional structure --- an associative
product. This is axiomatized as follows (this definition is taken
from \cite[Subsection 6.2]{tD}).

\begin{defn}
A {\em Green functor} is a Mackey functor $M \in \M(G)$ equipped
with an associative product in each $M^H$ such that
\begin{enumerate}
\item for any $f$, the map $f^*$ preserves the product,
\item for any $f:[G/H_1] \to [G/H_2]$, $x \in M^{H_2}$, $y \in M^{H_1}$,
  we have
$$
x \cdot f_*(y) = f_*(f^*(x) \cdot y), \qquad f_*(y) \cdot x = f_*(y
\cdot f^*(x)).
$$
\end{enumerate}
\end{defn}

We note that the Cartesian product of finite sets defines a functor
$m:Q\Gamma_G \times Q\Gamma_G \to Q\Gamma_G$; this induces a
symmetric tensor product on the category $\M(G)$ by
$$
M \otimes N = \Add(m_!(M \boxtimes N)).
$$
A Green functor is then the same as a Mackey functor equipped with
an algebra structure in the symmetric tensor category $\M(G)$. For
example, to see the condition \thetag{i}, one can argue as
follows. Consider the natural embedding $i:\Gamma_G^{opp} \to
Q\Gamma_G$. Then for any $M,N,K \in \M(G)$ we have
\begin{align*}
\Hom(M \otimes N,K) &= \Hom(\Add(m_!(M \boxtimes N)),K)\\
&= \Hom(m_!(M \boxtimes N),K) = \Hom(M \boxtimes N,m^*K),
\end{align*}
and since $m^*$ commutes with $i^*$, every map $M \otimes N \to K$
induces a map $i^*M \otimes i^*N \to i^*K$, where the tensor product
on $\Fun(\Gamma_G^{opp},\Ab)$ is again given by the direct image
$m_!$ with respect to the product functor $m:\Gamma_G^{opp} \times
\Gamma_G^{opp} \to \Gamma_G^{opp}$. However, on $\Gamma_G^{opp}$,
this product functor is left-adjoint to the diagonal embedding
$\delta:\Gamma^{opp}_G \to \Gamma^{opp}_G \times \Gamma_G^{opp}$;
therefore $m_! \cong \delta^*$, and the tensor structure on
$\Fun(\Gamma_G^{opp},\Ab)$ is in fact given by the pointwise
product, so that for any Green functor $M \in \M(G)$, the
restriction $i^*M$ is simply a functor from $\Gamma_G^{opp}$ to the
category of rings.

Since the point orbit $[G/G] \in Q\Gamma_G$ is the unit object for the
product functor $m$, the Mackey functor $\A$, $[G/H] \mapsto
\B^G([G/G],[G/H])$ it represents is the unit object for the tensor
product of Mackey functors and in particular, a Green functor. This
Green functor $\A$ is called the {\em Burnside ring Green functor}
--- indeed, by \eqref{b.g.hom}, the component $\A^H$ is exactly the
Burnside ring of the finite group $H$. Every Mackey functor $M \in
\M(G)$ is then a module over the Burnside ring Green functor $\A$.

\section{The derived version.}\label{mack.der}

Since the category $\M(G)$ of $G$-Mackey functors is abelian, one
can consider its derived category $\D(\M(G))$. However, as we have
explained in the Introduction, its formal properties are somewhat
deficient. The goal of this paper is to suggest a cure for this by
defining a certain triangulated category which contains $\M(G)$ but
differs from $\D(\M(G))$, and has all the properties one would like
to have.

The idea behind the construction is very simple. In the definition
of the Burnside ring $\A^G$, and more generally, in the definition
of the additive category $\B^G$, we take abelian groups spanned by
isomorphism classes of certain objects -- in other words, we take
the $0$-th homology group $H_0$ of a certain groupoid. The correct
thing to take at the level of triangulated categories is the full
homology, not just its degree-$0$ part.

There are several ways to make this precise. In this section, we
give a construction which uses $A_\infty$ methods and
bar-resolutions.

\subsection{The quotient construction.}\label{quo.subs}

Assume given a small category $\C$ which has fibered products. Then
we can obviously define a category $Q\C$ as follows: objects are
objects of $\C$, maps from $c_1 \in Q\C$ to $c_2 \in Q\C$ are
given by isomorphism classes of diagrams $c_1 \gets c \to
c_2$, and compositions are given by the fibered products, as in
Subsection~\ref{m.defn.subs}.

However, we can refine the construction. Let $\Q\C$ be the {\em
$2$-category} whose objects are again the objects of $\C$, and such
that for any $c_1,c_2 \in \C$, the category $\Q\C(c_1,c_2)$ of
maps from $c_1$ to $c_2$ in $\Q\C$ is the category of diagrams
$c_1 \gets c \to c_2$ and their isomorphisms. If $\C$ also has
a terminal object, thus all products, we can equivalently set
$\Q(c_1,c_2) = \Q(c_1 \times c_2)$, where $\Q(c)$ for an
object $c \in \C$ is the category of objects in $\C$ equipped with
a map to $c$ and their isomorphisms. The composition is again
given by fibered products.

Since fibered products are associative up to a canonical
isomorphism, $\Q\C$ is a well-defined $2$-category. We denote by
$\B_\idot^\C$ the $A_\infty$-category in the sense of
Subsection~\ref{dgcat.subs} associated to $\Q\C$ by the procedure
described in Subsection~\ref{2cat.subs}. Thus the objects in
$\B^\C_\idot$ are again the same as in $\C$, maps from $c_1$ to
$c_2$ are given by
$$
\B^\C_\idot(c_1,c_2) = C_\idot(\Q\C(c_1,c_2),\Z),
$$
where $\Z$ in the right-hand side is the constant functor with value
$\Z$, and compositions are induced by the $2$-category structure on
$\Q\C$. We note that by definition, this $A_\infty$-category is
concentrated in non-positive cohomological degrees.

\begin{defn}
The derived category $\D(\B^\C_\idot,\Ab)$ of $A_\infty$-functors
from $\B^\C_\idot$ to $\Ab$ is denoted by $\DQ(\C,\Ab)$.
\end{defn}

We note that every map $f:c_1 \to c_2$ in the category $\C$
canonically defines a $1$-map from $c_1$ to $c_2$ in the
$2$-category $\Q\C$, and we have a $2$-functor $\C \to \Q\C$, where
$\C$ is understood as a discrete $2$-category. The bar complex of a
discrete category with the set of objects $S$ is canonically
quasiisomorphic to the free abelian group $\Z[S]$ generated by $S$;
therefore by restriction, we obtain a canonical functor $\DQ(\C,\Ab)
\to \D(\C,\Ab)$. Analogously, we have a canonical functor
$\DQ(\C,\Ab) \to \D(\C^{opp},\Ab)$.

Assume that the category $\C$, in addition to fibered products, has
finite coproducts.

\begin{defn}\label{add.dq}
An $A_\infty$-functor $M_\idot \in \DQ(\C,\Ab)$ is {\em additive} if
its restriction $\ol{M}_\idot \in \D(\C^{opp},\Ab)$ is additive in
the sense of Subsection~\ref{m.defn.subs}: for any $c,c' \in \C$,
the natural map
$$
\ol{M}_\idot(c \copr c') \to \ol{M}_\idot(c)
\oplus \ol{M}_\idot(c')
$$
induced by the embeddings $c \to c \copr c'$, $c' \to c
\copr c'$ is a quasiisomorphism. The full subcategory in
$\DQ(\C,\Ab)$ spanned by additive $A_\infty$-functors is denoted by
$\DQ_{add}(\C,\Ab) \subset \DQ(\C,\Ab)$.
\end{defn}

In particular, let $G$ be a group, and let $\Gamma_G$ be the
category of finite $G$-sets.

\begin{defn}\label{d.m.defn}
A {\em derived $G$-Mackey functor} is an additive object in the
category $\DQ(\Gamma_G,\Ab)$. The category $\DQ_{add}(\Gamma_G,\Ab)$
of derived $G$-Mackey functors is denoted by $\DM(G,\Ab)$. If $\Ab =
\Z\amod$ is the category of abelian groups $\DM(G,\Z\amod)$ is
denoted simply by $\DM(G)$.
\end{defn}

\subsection{Example: the trivial group.}\label{gamma.subs}

To illustrate the general notion of a derived Mackey functor,
consider the case of the trivial group $G = \{e\}$, so that
$\Gamma_G = \Gamma$. Of course, $\M(\{e\},\Ab) = \Ab$; we would
expect the same to hold of the derived level. Let us check that this
is indeed so.

To do this, consider the subcategory $\Gamma_+ \subset \Q\Gamma$
with the same objects, and those $1$-morphisms $S_1 \gets S \to S_2$
for which the map $S \to S_1$ is injective. We note that such
diagrams have no non-trivial automorphisms; therefore the
$2$-category structure on $\Gamma_+$ is trivial and we can treat it
as a usual category. Here are two equivalent descriptions of
$\Gamma_+$.
\begin{enumerate}
\item The category whose objects are finite sets $S$, and whose
  morphisms from $S_1 \to S_2$ are ``partial maps'' $f:S_1 \ratto
  S_2$ -- that is, maps to $S_2$ defined on a subset $S \subset
  S_1$.
\item The category of finite sets $S$ with a fixed point $1 \in S$.
\end{enumerate}
Here \thetag{i} is just a restatement of the definition, and the
passage from \thetag{i} to \thetag{ii} is by formally adding the
fixed point (on morphisms, all elements in the set $S_1$ where a
partial map $f:S_1 \ratto S_2$ is undefined go into the added fixed
point). Denote by
$$
\lambda^*:\DQ(\Gamma,\Ab) \to \D(\Gamma_+,\Ab)
$$
the functor given by restriction with respect to the embedding
$\lambda:\Gamma_+ \to \Q\Gamma$.

For any finite set $S$, denote by $T(S) = \Z[S]$ the free abelian
group it generates. Then the correspondence $S \mapsto T(S)$
obviously defines an object $T \in \DQ(\Gamma,\Z\amod)$: for any map
$f:S_1 \to S_2$, the map $f_*:T(S_1) \to T(S_2)$ is induced by $f$,
and the map $f^*:T(S_2) \to T(S_1)$ is the adjoint map, 
$$
f^*([s]) = \sum_{s' \in f^{-1}(s)} [s']
$$
for any element $s \in S_2$. The object $T \in \DQ(\Gamma,\Z\amod)$
is additive in the sense of Definition~\ref{add.dq}. Restricting it
to $\Gamma_+$ gives an object $\lambda^*(T) \in
\Fun(\Gamma_+,\Z\amod)$ which we will denote by the same letter $T$
by abuse of notation.

\begin{lemma}\label{DQ.gamma}
\begin{enumerate}
\item For any $M_\idot \in \D(\Ab)$ and any $M'_\idot \in
  \DQ(\Gamma,\Ab)$, the natural map
\begin{equation}\label{QG.Gplus}
\RHom^\hdot_{\DQ(\Gamma,\Ab)}(M'_\idot,T \otimes M_\idot) \to
\RHom^\hdot_{\D(\Gamma_+,\Ab)}(\lambda^*(M'_\idot),T \otimes M_\idot)
\end{equation}
is an isomorphism.
\item The functor $D(\Ab) \to \DQ(\Gamma,\Ab)$ given by
  $M_\idot \mapsto T \otimes M_\idot$ is the full embedding onto
  $\DQ_{add}(\Gamma,\Ab)$.
\end{enumerate}
\end{lemma}

\proof{} We will need some semi-obvious facts on the structure of
the category $\Fun(\Gamma_+,\Z\amod)$ (see e.g. \cite[Section
3.2]{K}). A standard set of projective generators of this category
is given by representable functors $T_n$, $n \geq 0$, explicitly
described by
$$
T_n(S) = \Z[(S \copr \{1\})^n].
$$
We have $T_n = T_1^{\otimes n}$. In particular, $T_0 = \Z$, the
constant functor with value $\Z$. Moreover, we have a direct sum
decomposition $T_1 = T \oplus T_0$. Therefore the tensor powers
$T^{\otimes n}$ are also projective, and give another set of
generators for the category $\Fun(\Gamma_+,\Z\amod)$. These
generators are semi-orthogonal: we have $\Hom(T^{\otimes
n},T^{\otimes m}) = 0$ when $n > m$. In addition,
$\Hom(T_0,T^{\otimes n})=0$ for any $n \geq 1$. Explicitly,
\begin{equation}\label{T.S}
T^{\otimes n}(S) = \Z[S^n]
\end{equation}
for any $S \in \Gamma_+$. We also note that we have $\Hom(T,T) =
\Z$, which immediately implies that $M \mapsto T \otimes M$ gives
fully faithful embeddings $\Ab \subset \Fun(\Gamma_+,\Ab)$, $\D(\Ab)
\subset \D(\Gamma_+,\Ab)$.

The category $\DQ(\Gamma,\Ab)$ is generated by representable
$A_\infty$-functors $M^S_\idot$ of the form
$$
M^S_\idot(S') = C_\idot(\Q(S,S'),\Z) \otimes M
$$
for all $M \in \Ab$, $S \in \Gamma$. Therefore it is sufficient to
prove \thetag{i} for objects of this form. Fix a finite set $S' \in
\Gamma$ and an object $M' \in \Ab$, and let $M'_\idot =
M^{'S'}_\idot$. Explicitly, we have
$$
M'_\idot(S) = \bigoplus_{\wt{S} \in
  \Gamma}C_\idot(\Sigma_{\wt{S}},M' \otimes \Z[\Gamma(\wt{S},S
  \times S')]),
$$
where $\Sigma_{\wt{S}}$ is the group of automorphisms of the finite
set $\wt{S}$. This direct sum decomposition is not functorial with
respect to $S$. However, if we restrict to $\Gamma_+$, then the
increasing filtration $F_\idot$ on $\lambda^*(M'_\idot)$ given by
$$
F_n\lambda^*(M'_\idot)(S) = \bigoplus_{|\wt{S}| \leq
n}C_\idot(\Sigma_{\wt{S}},M' \otimes \Z[\Gamma(\wt{S},S \times
S')])
$$
is functorial ($|\wt{S}|$ denotes the cardinality of the set
$\wt{S}$). The associated graded quotient is given by
\begin{equation}\label{gr.Tn}
\gr^F_n\lambda^*(M'_\idot) = C_\idot(\Sigma_{\wt{S}},M' \otimes
\Z[\Gamma(\wt{S},S')] \otimes T^{\otimes n}),
\end{equation}
where we have used \eqref{T.S}, and $\wt{S}$ is the set of
cardinality $n$. By semi-or\-tho\-go\-na\-lity of the generators
$T^{\otimes n}$, this implies that
$$
\RHom^\hdot(\gr^F_n\lambda^*(M'_\idot),T \otimes M_\idot) = 0
$$
for $n \neq 1$, so that
\begin{align*}
\RHom^\hdot(\lambda^*(M'_\idot),T \otimes M_\idot) &=
\RHom^\hdot(\gr^F_1\lambda^*(M'_\idot),T \otimes M_\idot)\\
&= \RHom^\hdot_{\D(\Gamma_+,\Ab)}(T
\otimes M',T \otimes M_\idot) \otimes \Z[S']\\
&= \RHom^\hdot_{\Ab}(M',M_\idot) \otimes \Z[S'].
\end{align*}
Since $M'_\idot = M^{'S'}_\idot$ is representable, this is exactly
the left-hand side of \eqref{QG.Gplus}, so that we have proved
\thetag{i}.

As for \thetag{ii}, \thetag{i} immediately implies that the functor
$\D(\Ab) \to \DQ(\Gamma,\Ab)$ is a full embedding, and since $T \in
\DQ(\Gamma,\Z\amod)$ is additive, we in fact have a full embedding
$\D(\Ab) \subset \DQ_{add}(\Gamma,\Ab)$. To prove that it is
essentially surjective, it suffices to prove that its contains all
objects $M_\idot \in \DQ_{add}(\Gamma,\Ab)$ which are concentrated
in a single cohomological degree. But such objects are Mackey
functors in the usual non-derived sense, so they are all of the form
$T \otimes M$, $M \in \Ab$.
\endproof

\subsection{Wreath products.}

Definition~\ref{d.m.defn} is a DG version of the first definition of
a Mackey functor given in Subsection~\ref{m.defn.subs}. To get a
more explicit description of the category $\DM(G)$, we need to
somehow use the additivity condition and replace the
$A_\infty$-category $\B^{\Gamma_G}_\idot$ with an
$A_\infty$-category whose objects are $G$-orbits, not all finite
$G$-sets. We do it by using the structure of a ``wreath product'' of
the category $\Gamma_G$ of finite $G$-sets.

For any small category $\C$, by the {\em wreath product} $\C \wr
\Gamma$ of $\C$ with the category $\Gamma$ of finite sets we will
understand the category of pairs $\langle S,\{c_s\} \rangle$ of a
finite set $S$ and a collection of objects $c_s \in \C$, one for
each element $s \in S$, with maps from $\langle S,\{c_s\} \rangle$
to $\langle S',\{c'_s\} \rangle$ being a pair $\langle f,f_s
\rangle$ of a map $f:S \to S'$ and a collection of maps $f_s:c_s \to
c'_{f(s)}$, one for each $s \in S$.

By definition, we have a forgetful functor $\rho:\C \wr \Gamma \to
\Gamma$, $\langle S,\{c_s\}\rangle \mapsto S$. The functor $\rho$ is
a fibration; its fiber $\rho_S$ over a finite set $S \in \Gamma$ is
canonically identified with $\C^S$, the product of copies of $\C$
indexed by elements of the set $S$. In particular, $\C$ itself is
naturally embedded into $\C\wr\Gamma$ as the fiber over the
one-element set $\ppt \in \Gamma$. We will denote this embedding by
$j^\C:\C \to \C \wr \Gamma$.

Irrespective of the properties of the category $\C$, the category
$\C \wr \Gamma$ has finite coproducts. In a sense, it is obtained by
formally adjoining finite coproducts to $\C$ -- this can be
formulated precisely as a certain universal property of wreath
products, but we will not need this. Another way to characterize
$\C\wr\Gamma$ by a universal property is the following: for any
category $\C'$ fibered over $\Gamma$, any functor $f:\C'_\ppt \to
\C$ from the fiber $\C'_\ppt$ over the one-element set $\ppt \in
\Gamma$ extends uniquely to a Cartesian functor $\C' \to
\C\wr\Gamma$. We will need the following easy corollary of this
fact.

\begin{lemma}\label{wr.fib.lem}
Assume given a small category $\C$ and an object $S \in
\C\wr\Gamma$. Then the category $(\C\wr\Gamma)/S$ of objects $S' \in
\C\wr\Gamma$ equipped with a map $S' \to S$ is naturally equivalent
to the wreath product $(\C/S)\wr\Gamma$, where $\C/S$ is the
category of objects $c \in \C$ equipped with a map $j^\C(c) \to
S$.
\end{lemma}

\proof{} The projection $(\C\wr\Gamma)/S \to \Gamma$ which sends $S'
\to S$ to $\rho(S')$ is obviously a fibration; the universal
property of wreath products then gives a Cartesian comparison
functor
$$
(\C\wr\Gamma)/S \to (\C/S)\wr\Gamma,
$$
which is obviously an equivalence on all the fibers, thus an
equivalence.
\endproof

In the assumptions of Lemma~\ref{wr.fib.lem}, denote by $\Q^\wr(S)
\subset (\C\wr\Gamma)/S$ the subcategory with the same objects as
$(\C\wr\Gamma)/S$ and those maps which are Cartesian with respect to
the fibration $(\C\wr\Gamma)/S \to \Gamma$. Equivalently, we have
$\Q^\wr(S) \cong \ol{\C/S}\wr\Gamma$, where $\ol{\C/S}
\subset \C/S$ is the subcategory whose objects are all objects in
$\C/S$, and whose maps are isomorphisms in $\C/S$.

More generally, given two objects $S_1,S_2 \in \C\wr\Gamma$, denote
by $(\C\wr\Gamma)/(S_1,S_2)$ the category of objects $S \in \C\wr\Gamma$
equipped with maps $S \to S_1$, $S \to S_2$, let $\C/(S_1,S_2)
\subset (\C\wr\Gamma)/(S_1,S_2)$ be the fiber of the
projection $(\C\wr\Gamma)/(S_1,S_2) \to \Gamma$ given by $\langle
S_1 \gets S \to S_2 \rangle \mapsto \rho(S)$, so that
$(\C\wr\Gamma)/(S_1,S_2) \cong (\C/(S_1,S_2))\wr\Gamma$, and let
$\ol{\C/(S_1,S_2)} \subset \C/(S_1,S_2)$ be the groupoid of
isomorphisms of the category $\C/(S_1,S_2)$. Denote
$$
\Q^\wr(S_1,S_2) = \ol{\C/(S_1,S_2)}\wr\Gamma \subset
(\C\wr\Gamma)/(S_1,S_2),
$$
and assume that the category $\C\wr\Gamma$ has fibered
products. Then these fibered products define associative composition
functors
\begin{equation}\label{m.123}
m:(\C\wr\Gamma)/(S_1,S_2) \times (\C\wr\Gamma)/(S_2,S_3) \to
(\C\wr\Gamma)/(S_1,S_3)
\end{equation}
which induce composition functors on the categories
$\Q^\wr(-,-)$. This allows to define a $2$-category which we denote
by $\Q^\wr\C$: it objects are the objects of $\C\wr\Gamma$, and its
categories of morphisms are $\Q^\wr(-,-)$.

Note that for any $S_1,S_2 \in \C\wr\Gamma$, we have a natural
embedding $F_\C(S_1,S_2):\Q(S_1,S_2) \to \Q^\wr(S_1,S_2)$, and these
embeddings are compatible with fibered products, thus glue together
to a $2$-functor
$$
\F_\C:\Q(\C\wr\Gamma) \to \Q^\wr(\C).
$$
Both $2$-categories here have the same objects and the same
$1$-morphisms; the difference is that the right-hand side has more
$2$-morphisms -- the embeddings $F_\C(-,-)$ are identical on objects
and faithful, but not full.

\bigskip

We prove right away one technical results on the categories
$\Q^\wr(-,-)$ which we will need later on. Any object $c \in \C
\subset \C\wr\Gamma$ defines a functor $j^c:\Gamma \to \C\wr\Gamma$
which sends a finite set $S$ to the union of $S$ copies of $c$. This
functor preserves fibered products, thus gives a $2$-functor
$j^c:\Q\Gamma \to \Q(\C\wr\Gamma)$, a restriction functor
$j^{c*}:\DQ(\C\wr\Gamma,\Ab) \to \DQ(\Gamma,\Ab)$, and a
right-adjoint functor $j^c_*:\DQ(\Gamma,\Ab) \to
\DQ(\C\wr\Gamma,\Ab)$.

\begin{lemma}\label{adj.DQ}
For any $E_\idot \in \DQ(\Gamma,\Ab)$ and any $S \in \C\wr\Gamma$,
we have a natural quasiisomorphism
$$
j^c_*E_\idot(S) \cong C_\idot(\Q^\wr(c,S)^{opp},\rho^{opp *}E_\idot),
$$
where $E_\idot$ in the right-hand side is restricted to
$\Gamma^{opp} \subset \Q\Gamma$ and then pulled back to
$\Q^\wr(c,S)^{opp}$ by the opposite $\rho^{opp}$ to the projection
$$
\rho:\Q^\wr(c,S)=\overline{\C/(c,S)}\wr\Gamma \to \Gamma.
$$
\end{lemma}

\proof{} Denote the embedding $\Gamma^{opp} \to \Q\Gamma$ by
$\lambda$, and let $\lambda_!:\D(\Gamma^{opp},\Ab) \to
\DQ(\Gamma,\Ab)$ be the left-adjoint functor to the restriction
functor $\lambda^*$. Then for any $M_\idot \in \D(\Ab)$ we have
$$
\RHom^\hdot(M_\idot,C_\idot(\Q^\wr(c,S)^{opp},\rho^{opp *}E_\idot))
\cong \RHom^\hdot(\lambda_!\rho^{opp}_!\tau^*M_\idot,E_\idot),
$$
where $\tau:\Q^\wr(c,S)^{opp} \to \ppt$ is the tautological
projection. Thus by adjunction, it suffices to prove that
$$
\lambda_!\rho^{opp}_!\tau^*M_\idot \cong j^{c*}M_\idot^S.
$$
To construct a map $f:\lambda_!\rho^{opp}_!\tau^*M_\idot \to
j^{c*}M_\idot^S$, it suffices by adjunction to construct a map
$\tau^*M_\idot \to \rho^{opp*}\lambda^*j^{c*}M_\idot^S$, that is, a
compatible system of maps
$$
M_\idot \to \rho^{opp*}\lambda^*j^{c*}M_\idot^S(S') \cong
H_\idot(\Q(S,j^c(\rho(S'))),M_\idot)
$$
for any $S' \in \Q^\wr(c,S)$; these maps are induced by obvious
tautological maps $\Z \to H_\idot(\Q(S,j^c(\rho(S'))),\Z)$.

To prove that the map $f$ is an isomorphism, we first need a way to
control the functor $\lambda_!$. For any finite set $\ol{S} \in
\Gamma$, let $\ol{\Gamma}/\ol{S}$ be the category of all finite sets
$\ol{S'}$ equipped with a map $\ol{S'} \to \ol{S}$ and isomorphisms
between them, and let $\kappa^{\ol{S}}:(\ol{\Gamma}/\ol{S})^{opp}
\to \Gamma^{opp}$ be the natural projection which sends $[\ol{S'}
\to \ol{S}]$ to $\ol{S'}$ Then we obviously have
$$
\lambda^*\Z_{\ol{S}}^\hdot \cong \kappa^{\ol{S}}_*\Z \in
\D(\Gamma^{opp},\Z\amod),
$$
and by adjunction, this yields a canonical isomorphism
$$
\lambda_!N_\idot(\ol{S}) \cong
H_\idot((\ol{\Gamma}/\ol{S})^{opp},\kappa^{\ol{S}*}N_\idot),
$$
for any $N_\idot \in \D(\Gamma^{opp},\Ab)$.

Now apply this to $N_\idot = \rho^{opp}_!\tau^*M_\idot$. Since
$\rho^{opp}$ is a cofibration, we may compute
$\kappa^{\ol{S}*}\rho^{opp}_!$ by base change; this gives a
canonical isomorphism
$$
H_\idot((\ol{\Gamma}/\ol{S})^{opp},\kappa^{\ol{S}*}\rho^{opp}_!\tau^*M_\idot)
\cong H_\idot(\Q'(c,S,\ol{S})^{opp},M_\idot),
$$
where $\Q'(c,S,\ol{S})$ is the category obtained as the Cartesian
product
$$
\begin{CD}
\Q'(c,S,\ol{S}) @>>> \ol{\Gamma}/\ol{S}\\
@VVV @VVV\\
\Q^\wr(c,S) @>>> \Gamma.
\end{CD}
$$
It remains to notice that the category $\Q'(c,S,\ol{S})$ is
canonically identified with $\Q(S,j^c(\ol{S}))$, so that
$$
H_\idot(\Q'(c,S,\ol{S})^{opp},M_\idot) \cong
H_\idot(\Q(S,j^c(\ol{S}))^{opp},M_\idot) \cong j^{c*}M_\idot^S(\ol{S}).
$$
Thus the map $f$ becomes an isomorphism after evaluating at every
object $\ol{S} \in \Q\Gamma$, as required.
\endproof

\subsection{Additivization of the quotient construction.}

Now let us consider again the functor $T \in \DQ(\Gamma,\Z\amod)$ of
Subsection~\ref{gamma.subs}, and let us restrict it to an object $T
\in \Fun(\Gamma,\Z\amod)$ by the embedding $\Gamma \to
\Q\Gamma$. This $T \in \Fun(\Gamma,\Z\amod)$ is isomorphic to the
functor $\Z_{[1]}$ represented by the set $[1] \in \Gamma$ with a
single element. For any $S \in \Gamma$, let $\tau_S:\Z \to \Z[S] =
T(S)$ be the diagonal embedding. The maps $\tau_S$ are {\em not}
functorial with respect to arbitrary maps of finite sets $S$;
however, they are functorial with respect to isomorphisms. Thus if
we denote by $\ol{\Gamma} \subset \Gamma$ the category of finite
sets and their isomorphisms, then we have a map of functors
\begin{equation}\label{tau.t}
\tau:\Z \to t^*T,
\end{equation}
where $\Z \in \Fun(\ol{\Gamma},\Z\amod)$ is the constant
functor with value $\Z$, and $t:\ol{\Gamma} \to \Gamma$ is the
embedding.

For any small category $\C$, we will denote the pullback $\rho^*T
\in \Fun(\C \wr \Gamma,\Z)$ with respect to the forgetful functor
$\rho:\C \wr \Gamma \to \Gamma$ by the same letter $T$. By base
change, we have
$$
T \cong \rho^*\Z_{[1]} \cong j^\C_!\Z_\C,
$$
where $\Z_\C \in \Fun(\C,\Z)$ is the constant functor with value
$\Z$, and $L^ij^\C_!\Z_\C = 0$ for $i \geq 1$. Therefore the natural
map
\begin{equation}\label{T.qis}
H_\idot(\C,\Z_\C) \to H_\idot(\C \wr \Gamma,T)
\end{equation}
is an isomorphism.

Assume that the wreath product category $\C\wr\Gamma$ has fibered
products, and consider the $2$-category $\Q^\wr(\C)$. Then for any
$S_1,S_2 \in \C\wr\Gamma$, the category $\Q^\wr(S_1,S_2) \cong
\ol{\C/(S_1,S_2)}\wr\Gamma$ carries a natural $\Z\amod$-valued
functor
$$
T \cong \rho^*T \cong j^{S_1,S_2}_!\Z \in
\Fun(\Q^\wr(S_1,S_2),\Z\amod),
$$
where $j^{S_1,S_2}:\ol{\C/(S_1,S_2)} \to \Q^\wr(S_1,S_2)$
denotes the natural embedding. Moreover, for any $S_1,S_2,S_3 \in
\C\wr\Gamma$, the composition functors \eqref{m.123} induce functors
$$
\ol{m} = (j^{S_1,S_2} \times j^{S_2,S_3}) \circ m:
\ol{\C/(S_1,S_2)} \times \ol{\C/(S_2,S_3)} \to
\Q^\wr{S_1,S_2},
$$
and by definition, all the maps in the category
$\ol{\C/(S_1,S_2)}$ are invertible, so that the composition
$\rho \circ \ol{m}$ actually goes into $\ol{\Gamma}
\subset \Gamma$. Therefore the canonical map $\tau$ of \eqref{tau.t}
induces a map $\tau_{S_1,S_2,S_3}:\Z \to \ol{m}^*T$. These
maps are associative on triple products in the obvious sense. By
adjunction, they induce maps
\begin{equation}\label{mu.boxtimes}
\mu_{S_1,S_2,S_3}:T \boxtimes T \cong (j^{S_1,S_2} \times
j^{S_2,S_3})_1\Z \to m^*T,
\end{equation}
and these maps are also associative on triple products.

As in the Subsection~\ref{quo.subs}, the procedure of
Subsection~\ref{2cat.subs} gives an $A_\infty$-category
$\B^{\wr\C}_\idot$ with the same objects as $\C \wr \Gamma$ by
setting
$$
\B^{\wr\C}_\idot(S_1,S_2) = C_\idot(\Q^\wr\C(S_1,S_2),T),
$$
with compositions induced by the $2$-category structure on
$\Q^\wr\C$ and the canonical maps $\mu$ of \eqref{mu.boxtimes}.

\begin{defn}\label{dq.wr.defn}
The derived category $\D(\B^{\wr\C}_\idot,\Ab)$ of $A_\infty$-functors
from $\B^{\wr\C}_\idot$ to $\Ab$ is denoted by $\DQ^\wr(\C,\Ab)$.
\end{defn}

We note that for any $S_1,S_2 \in \C\wr\Gamma$, the natural functor
$F:\Q(S_1,S_2) \to \Q^\wr(S_1,S_2)$ again goes into $\ol{\Gamma}
\subset \Gamma$ when composed with the projection
$\rho:\Q^\wr(S_1,S_2) \to \Gamma$; therefore the map $\tau$ of
$\eqref{tau.t}$ induces maps $\Z \to F^*T$, and the $2$-functor
$\F_\C:\Q(\C\wr\Gamma) \to \Q^\wr\C$ extends to an
$A_\infty$-functor $\F_\C:\B^{\C\wr\Gamma}_\idot \to
\B^{\wr\C}_\idot$. The main comparison result that we want to prove
is the following.

\begin{prop}\label{wr.nowr}
The functor $\F^*_\C:\DQ^\wr(\C,\Ab) \to \DQ(\C\wr\Gamma,\Ab)$
induced by $\F_\C$ gives an equivalence
$$
\DQ^\wr(\C,\Ab) \to \DQ_{add}(\C\wr\Gamma,\Ab).
$$
\end{prop}

\begin{remark}\label{loop}
The appearance of the functor $T$ in the definition of the
$A_\infty$-category $\B^{\wr\C}_\idot(-,-)$ looks like a trick. One
motivation for this comes from topology. The categories $\Q\C(-,-)$
are symmetric monoidal categories with respect to the disjoint
union. Applying group completion to their classifying spaces
$|\Q\C(-,-)|$, one obtains infinite loop spaces. Then the complex
$\B^{\wr\C}_\idot(-,-)$ simply computes the homology of the
corresponding $\Omega$-spectrum (as opposed to $\B^\C(-,-)$, which
computes the homology of the classifying space $|\Q\C(-,-)|$). We do
not prove this since we will not need it, but the proof is not
difficult (for example, it can be done along the lines of
\cite[Section 3.2]{K}).
\end{remark}

\subsection{The proofs of the comparison results.}

Before we prove Proposition~\ref{wr.nowr}, let us explain why it is
useful. Assume given a small category $\C$ such that $\C\wr\Gamma$
has fibered products, and denote by
\begin{equation}\label{wr.defn}
\wt{\B}^\C_\idot \subset \B^{\wr\C}_\idot
\end{equation}
the full subcategory spanned by $\C \subset \C\wr\Gamma$. Let
$\wt{\DQ}(\C,\Ab)$ be the derived category of $A_\infty$-functors
from $\wt{\B}^\C_\idot$ to $\Ab$.

\begin{lemma}\label{add.wr}
Restriction with respect to the natural embedding $\wt{\B}^\C_\idot
\to \B^{\wr\C}_\idot$ induces an equivalence
$$
R:\DQ^\wr(\C,\Ab) \overset{\sim}{\longrightarrow} \wt{\DQ}(\C,\Ab).
$$
\end{lemma}

\proof{} Let $E:\wt{\DQ}(\C,\Ab) \to \DQ^\wr(\C,\Ab)$ be the
left-adjoint functor to $R$. It suffices to prove that $E$ is
essentially surjective, and that the adjunction map $\Id \to R \circ
E$ is an isomorphism. The second fact is obvious: by adjunction, $E$
sends representable functors into representable functors, and since
$\wt{\B}^\C_\idot \subset \B^{\wr\C}_\idot$ is a full subcategory,
$E$ does not change the spaces of maps between them.  It remains to
prove that $E$ is essentially surjective. Since $\DQ^\wr(\C,\Ab)$ is
generated by representable functors $M^S_\idot$, $S \in
\C\wr\Gamma$, $M \in \Ab$, it suffices to prove that all these
functors lie in the essential image of $E$. By induction on the
cardinality $|\rho(S)|$, it suffices to prove that
$$
M^{S_1 \copr S_2}_\idot \cong M^{S_1}_\idot \oplus M^{S_1}_\idot
$$
for any $M \in \Ab$, $S_1,S_2 \in \C\wr\Gamma$.

Indeed, assume given such an $M$ and $S_1$, $S_2$. Explicitly, we
have
$$
M^{S_1 \copr S_2}_\idot(S') = C_\idot(\Q^\wr(S,S'),T) \otimes M
$$
for any $S' \in \C\wr\Gamma$. We have
\begin{equation}\label{pr.copr}
\Q^\wr(S_1 \copr S_2,S') \cong \Q^\wr(S_1,S') \times
\Q^\wr(S_2,S'),
\end{equation}
and
\begin{equation}\label{T.copr}
T \cong (\pi_1^*T \boxtimes \pi_2^*\Z) \oplus (\pi_1^*\Z
\boxtimes \pi_2^*T),
\end{equation}
where $\pi_1$ and $\pi_2$ are the projections onto the factors of
the decomposition \eqref{pr.copr}, and $\Z$ means the constant
functor with value $\Z$. The direct sum decomposition \eqref{T.copr}
is functorial with respect to $S'$, thus induces a certain direct sum
decomposition $\M^{S_1 \copr S_2}_\idot \cong M^1_\idot \oplus
M^2_\idot$. Here $M^1_\idot$ is given by
$$
M^1_\idot(S') = C_\idot(\Q^\wr(S_1,S') \times
\Q^\wr(S_2,S'),\pi_1^*T \boxtimes \pi_2^*\Z),
$$
and by the K\"unneth formula, this is canonically quasiisomorphic to
$$
M^{S_1}_\idot(S') \otimes C_\idot(\Q^\wr(S_2,S'),\Z).
$$
Since the category $\Q^\wr(S_2,S')$ is a wreath product, it has an
initial object, thus no homology with constant coefficients,
$H_\idot(\Q^\wr(S_2,S'),\Z) = \Z$, and we conclude that $M^1_\idot
\cong M^{S_1}_\idot$. Analogously $M^2_\idot \cong M^{S_2}_\idot$.
\endproof

Lemma~\ref{add.wr} shows that Proposition~\ref{wr.nowr} allows one
to get rid of the additivity assumption in the definition of the
category $\DQ_{add}(\C\wr\Gamma,\Ab)$ and reduces everything to an
$A_\infty$-category whose objects are those of $\C$, not of
$\C\wr\Gamma$. This has the following immediate corollary.

\begin{corr}\label{C.wr}
Assume that the small category $\C$ itself has fibered
products. Then there exists a natural equivalence of triangulated
categories
$$
\DQ(\C,\Ab) \cong \DQ^\wr(\C,\Ab).
$$
\end{corr}

\proof{} By Lemma~\ref{add.wr}, it suffices to construct an
equivalence $\DQ(\C,\Ab) \cong \wt{\DQ}(\C,\Ab)$, or equivalently, a
quasiisomorphism of $A_\infty$-categories $\B^\C_\idot \cong
\wt{\B}^\C_\idot$. Both $A_\infty$-categories have the same objects,
the objects of the category $\C$. For any $c_1,c_2 \in \C$, the
natural embedding $\Q(c_1,c_2) \to \Q^\wr(c_1,c_2) =
\Q(c_1,c_2)\wr\Gamma$ induces a map
$$
\B^\C_\idot(c_1,c_2) = C_\idot(\Q(c_1,c_2),\Z) \to
\wt{\B}^\C(c_1,c_2) = C_\idot(\Q^\wr(c_1,c_2),T),
$$
these maps are obviously compatible with compositions, and they are
all quasiisomorphisms by \eqref{T.qis}.
\endproof

We will now prove Proposition~\ref{wr.nowr}. To do this, recall that
for any object $c \in \C$, we have the embedding $j^c:\Q\Gamma \to
\Q(\C\wr\Gamma)$ and the corresponding restriction functor
$j^{c*}:\DQ(\C\wr\Gamma,\Ab) \to \DQ(\Gamma,\Ab)$. Moreover, $j^c$
obviously extends to an embedding $j^c_\wr:\Q^\wr(\ppt) \to
\Q^\wr(\C)$, and we have a restriction functor
$j_\wr^{c*}:\DQ^\wr(\C,\Ab) \to \DQ^\wr(\ppt,\Ab)$. These are
compatible with the functors $\F_\C^*$ of Proposition~\ref{wr.nowr}:
we have a commutative diagram
\begin{equation}\label{dia.wr.nowr}
\begin{CD}
\DQ(\Gamma,\Ab) @<{j^{c*}}<< \DQ(\C\wr\Gamma,\Ab)\\
@A{\F_{\ppt}^*}AA @AA{\F_\C^*}A\\
\DQ^\wr(\ppt,\Ab) @<{j_\wr^{c*}}<< \DQ^\wr(\C,\Ab).
\end{CD}
\end{equation}
Of course, by Corollary~\ref{C.wr} we have $\DQ^\wr(\ppt,\Ab) \cong
\DQ(\ppt,\Ab) \cong \D(\Ab)$, and $j_\wr^{c*}:\DQ^\wr(\C,\Ab) \to
\D(\Ab)$ is simply the evaluation at $c \in \Q^\wr(\C)$.

\begin{lemma}\label{2.bc}
Denote by $\F_{\C!}$, $\F_{\ppt!}$ the functors left-adjoint to
$\F_{\C}^*$ and $\F_{\ppt}^*$. Then the base change map
$$
\F_{\ppt!} \circ j^{c*} \to j_\wr^{c*} \circ \F_{\C!}
$$
obtained by adjunction from \eqref{dia.wr.nowr} is an isomorphism.
\end{lemma}

\proof{} Since $\DQ(\C\wr\Gamma,\Ab)$ is generated by representable
objects $M^S_\idot$, it suffices to prove that
$$
\RHom^\hdot(\F_{\ppt!}j^{c*}M_\idot^S,M'_\idot) \cong 
\RHom^\hdot(j_\wr^{c*}\F_{\C!}M_\idot^S,M'_\idot)
$$
for any $M_\idot,M'_\idot \in \D(\Ab)$, $S \in \C$. By adjunction,
$\F_{\C!}M^\idot_S = M_\idot^{\F_\C(S)}$, so that the right-hand
side is isomorphic to
$$
\RHom^\hdot_{\Q^\wr(S,c)}(T,\Z)) \otimes
\RHom^\hdot(M_\idot,M'_\idot).
$$
The left-hand side by adjunction is isomorphic to
$$
\RHom^\hdot(j^{c*}M^S_\idot,\F_{\ppt}^*M'_\idot) \cong
\RHom^\hdot(j^{c*}M^S_\idot,M'_\idot \otimes T),
$$
which is isomorphic to
$$
\RHom^\hdot_{\Q^\wr(S,c)^{opp}}(\Z,T)) \otimes
\RHom^\hdot(M_\idot,M'_\idot)
$$
by Lemma~\ref{adj.DQ}. It remains to notice that
$$
\RHom^\hdot_{\Q^\wr(S,c)^{opp}}(\Z,T)) =
\RHom^\hdot_{\Q^\wr(S,c)}(T,\Z))
$$
identically, if both sides are understood in the sense of
\eqref{rhom.bar.eq}.
\endproof

\proof[Proof of Proposition~\ref{wr.nowr}.] It suffices to check that
\begin{enumerate}
\item the adjunction map $\F_{\C!} \circ \F_{\C}^* \to \Id$ is an
isomorphism, so that $\F_{\C}^*$ is fully faithful, and 
\item an additive object $M_\idot \in \DQ_{add}(\C\wr\Gamma,\Ab)$
with trivial $\F_{\C!}(M_\idot)$ is itself trivial, so that
$F_{\C}^*$ is essentially surjective.
\end{enumerate}
By Lemma~\ref{add.wr}, it suffices to check \thetag{i} after
evaluating on all $c \in \C \subset \C\wr\Gamma$, and by
Lemma~\ref{2.bc}, this amount to checking \thetag{i} with $\C$
replaced by $\ppt$. Analogously, the restriction functor $j^{c*}$
obviously sends $\DQ_{add}(\C\wr\Gamma,\Ab)$ into
$\DQ_{add}(\Gamma,\Ab)$, and an object $M_\idot \in
\DQ(\C\wr\Gamma,\Ab)$ with trivial $j^{c*}(M_\idot)$ for all $c \in
\C$ is itself trivial; thus by Lemma~\ref{2.bc}, \thetag{ii} can
also be only checked for $\C = \ppt$. Conclusion: it suffices to
prove the Proposition for $\C=\ppt$. This has been done already --
combine Lemma~\ref{add.wr} and Corollary~\ref{C.wr}, on one hand,
and Lemma~\ref{DQ.gamma}, on the other hand.  \endproof

\subsection{Derived Burnside rings.}

Now again fix a finite group $G$ and take $\C = \Gamma_G$, the
category of finite $G$-sets. Define a functor $\rho:\Gamma_S \to
\Gamma$ by setting $S \mapsto S/G$, the set of $G$-orbits on
$S$. This functor is a fibration, and moreover, we actually have an
identification $\Gamma_G \cong O_G\wr\Gamma$, where $O_G$ is the
category of finite $G$-orbits, and $\rho$ is the tautological
projection $O_G \wr\Gamma \to \Gamma$. Therefore we can also
consider the $2$-category $\Q^\wr(O_G)$ and the associated
$A_\infty$-category $\wt{\B}^{O_G}_\idot$ of \eqref{wr.defn} whose
objects are finite $G$-orbits. To simplify notation, denote
$$
\B^G_\idot = \wt{\B}^{O_G}_\idot.
$$
then the following is a reformulation of Proposition~\ref{wr.nowr}
and Lemma~\ref{add.wr} (with $\C = O_G$).

\begin{prop}
The triangulated category $\DM(G)$ of derived $G$-Mackey functors is
equivalent to the derived category of $A_\infty$-functors
$\B^G_\idot \to \Ab$.\endproof
\end{prop}

This Proposition allows us to do some computations in the category
$\DM(G)$; in particular, spelling out the definitions, we can now
define a derived version of the Burnside ring $\A^G$. Let $T =
\rho^*T \in \Fun(\Gamma_G,\Z\amod)$, so that $T(S) = \Z[S/G] =
\Z[S]^G$. For every $S_1,S_2 \in \Gamma_G$, let
$$
\mu_{S_1,S_2}:\Z[S_1]^G \otimes \Z[S_2]^G =
\Z[S_1 \times S_2]^{G \times G} \to
\Z[S_1 \times S_2]^G
$$
be the natural embedding. Taken together, these maps give a map
\begin{equation}\label{mu}
\mu:T \boxtimes T \to m^*T,
\end{equation}
where, as in Section~\ref{mack.rec}, $m:\Gamma_G \times \Gamma_G \to
\Gamma_G$ is the product functor (this map $\mu$ is of course the
special case of \eqref{mu.boxtimes} for $S_1 = S_2 = S_3 =
\{\ppt\}$).

\begin{defn}
The {\em derived Burnside ring} $\A_\idot^G$ of the group $G$ is the
complex $C_\idot(\bGamma_G,T)$, with the $A_\infty$-structure
induced by the Cartesian product of $G$-sets and the canonical map
$$
m^{T}_* \circ \mu:C_\idot(\bGamma_G \times \bGamma_G,T \boxtimes T)
\to C_\idot(\bGamma_G \times \bGamma_G,m^*T) \to
C_\idot(\bGamma_G,T),
$$
where $\mu$ is as in \eqref{mu}.
\end{defn}

By definition, $\A_\idot^G$ is a $A_\infty$-algebra over $\Z$ (and
its homology algebra $H_\idot(\A_\idot^G)$ is commutative).  It is
isomorphic to $\B_G([G/G],[G/G])$, where $[G/G]$ is the trivial
$G$-orbit (the point set $\ppt$ with the trivial $G$-action).

\begin{lemma}
Assume given a group $G$.
\begin{enumerate}
\item The $0$-th homology $H_0(\A^G_\idot)$ of the derived Burnside
  ring $\A^G_\idot$ is isomorphic to the usual Burnside ring $\A^G$,
  and the $0$-th homology $H_0(\B^G_\idot)$ of the
  $A_\infty$-category $\B^G_\idot$ is isomorphic to the additive
  category $\B^G$ of Subsection~\ref{m.defn.subs}.
\item For any two $G$-orbits $S_1$, $S_2$, we have a natural
  quasiisomorphism
\begin{equation}\label{b.g.idot.hom}
\B^G_\idot(S_1,S_2) = \bigoplus_{p \in (S_1 \times S_2)/G} \A^{H_p}_\idot,
\end{equation}
where $S_p = [G/H_p]$ are the components in the decomposition
\eqref{unio} of the product $S = S_1 \times S_2$, and this
quasiisomorphism induces the isomorphism \eqref{b.g.hom} on $0$-th
homology.
\end{enumerate}
\end{lemma}

\proof{} The quasiisomorphism \eqref{T.qis} in our new notation
reads as
$$
\A^G_\idot \cong C_\idot(\ol{O}_G,\Z),
$$
so that the $0$-th homology of $A^G_\idot$ is the $0$-th homology of
the groupoid $\ol{O}_G$ of $G$-orbits; this is precisely the
Burnside ring $\A^G$. The decomposition \eqref{b.g.idot.hom} follows
from the decomposition
$$
\Q^\wr(S_1,S_2) \cong \prod_{p \in (S_1 \times S_2)/G} \Q^\wr(S_p)
$$
by the same argument as in the proof of Lemma~\ref{add.wr}.
Combining together \eqref{b.g.idot.hom}, \eqref{b.g.hom} and the
isomorphism $H_0(\A^G_\idot) \cong \A^G$ gives the isomorphism
$H_0(\B^G_\idot) \cong \B^G$. It remains to prove that this
isomorphism is compatible with the compositions; this is a
straightforward check which we leave to the reader.
\endproof

\section{Waldhausen-type description.}\label{wald.sec}

The construction of the triangulated category $\DM(G)$ of derived
Mackey functors given in the last Section is very explicit but
somewhat deficient, since it relies on explicit resolutions and
$A_\infty$ methods. We will now give a more invariant
construction. To do this, we modify the quotient construction
$\Gamma_G \mapsto Q\Gamma_G$ of Subsection~\ref{quo.subs} in a way
which is similar to the passage from Quillen's $Q$-construction to
Waldhausen's $S$-construction in algebraic $K$-theory.

\subsection{Heuristic explanation.}\label{heur.subs}

Let us first explain informally what we are going to do (this
Subsection is purely heuristic and may be skipped, formally, nothing
in the rest of the paper depends on it). Recall that a {\em
simplicial set} $X$ is by definition a contravariant functor from
the category $\Delta$ of finite non-empty totally ordered sets to
the category of all sets. For any non-negative integer $n \geq 1$,
we will denote by $[n] \in \Delta$ the set with $n$ elements, or, to
be specific, the set of all integers $i$, $1 \leq i \leq n$; we will
also denote $X_i = X([i])$ for any $i \geq 1$. An {\em $\Ab$-valued
sheaf} $M$ on $X$ is a collection of
\begin{enumerate}
\item a functor $M_n:X_n \to \Ab$ for every $n \geq 1$, and
\item a map $M(f):M_{n'} \to X(f)^*M_n$ for every map $f:[n] \to
[n']$,
\end{enumerate}
subject to standard compatibility conditions. Here in \thetag{i}, we
treat the set $X_n$ as a discrete category, so that $M_n$ is
effectively just an $X_n$-graded object in $\Ab$; in \thetag{ii},
$X(f):X_{n'} \to X_n$ is value of the functor $X:\Delta^{opp} \to
\Sets$ on the map $f$.

There is the following convenient way to pack together these data
\thetag{i}, \thetag{ii}, and also the compatibility conditions. Let
us not only treat the sets $X_n$ as discrete categories, but also
treat $X$ as a functor $\Delta^{opp} \to \Cat$. Then we can apply
the Grothendieck construction of Subsection~\ref{fib.subs}. The
result is a category $S(X)$ fibered over $\Delta$; explicitly,
objects in $S(X)$ are pairs $\langle n,x \in X_n \rangle$, and a map
from $\langle n,x \in X_n \rangle$ to $\langle n',x' \in X_n
\rangle$ is given by a map $f:[n] \to [n']$ such that $X(f)(x') =
x$. One immediately checks that in this notation, a sheaf $M$ on $X$
is exactly the same thing as a functor $M:S(X) \to \Ab$.

Now assume given a small category $\C$. Recall that to $\C$, one
canonically associates a simplicial set $N(\C)$ called the {\em
nerve} of $\C$, and to any functor $E \in \Fun(\C,\Ab)$, one
associates an $\Ab$-valued sheaf $\wt{E}$ on the nerve -- in other
words, we have a natural embedding
\begin{equation}\label{nerve.emb}
\Fun(\C,\Ab) \to \Fun(S(N(\C)),\Ab).
\end{equation}
Explicitly, $N(\C)_n$ is the set of diagrams $c_1 \to \dots \to c_n$
in the category $\C$; the functor $\wt{E}:S(N(\C)) \to \Ab$ sends a
diagram $c_1 \to \dots \to c_n$ to $E(c_n)$. For any map $f:[n'] \to
[n]$, the map $N(\C)(f)$ sends a diagram $c_1 \to \dots \to c_n$ to
the diagram $c_{f(1)} \to \dots \to c_{f(n')}$, and the
corresponding map $\wt{E}(f):E(c_{f(n')}) \to E(c_n)$ is induces by
the natural map $c_{f(n')} \to c_n$.

The embedding \eqref{nerve.emb} is fully faithful but not
essentially surjective -- not all sheaves on $N(\C)$ come from
functors $E \in \Fun(\C,\Ab)$. Indeed, for any sheaf of the form
$\wt{E}$, the map $\wt{E}(f)$ is an isomorphism whenever the map
$f:[n'] \to [n]$ sends $n'$ to $n$ -- that is, preserves the last
elements of the totally ordered sets. However, this is the only
condition: the essential image of embedding \eqref{nerve.emb}
consists of such $M \in \Fun(S(N(\C)),\Ab)$ that $M(f)$ is a
quasiisomorphism whenever $f$ preserves the last elements. Indeed,
it is easy to see that such a sheaf $M$ is completely defined by the
following part of \thetag{i}, \thetag{ii}:
\begin{enumerate}
\item the functor $M_1:N(\C)_1 \to \Ab$, and
\item the map 
\begin{equation}\label{act.map}
M(s):s^*M_1 \to M_2 \cong t^*M_1,
\end{equation}
where the maps $s,t:[1] \to [2]$ send $1 \in [1]$ to $1 \in [2]$,
resp. to $2 \in [2]$,
\end{enumerate}
subject to some conditions. This is exactly the same as a functor
from $\C$ to $\Ab$: $M_0$ gives its values on objects of the
category $\C$, and the map $s^*M_1 \to t^*M_1$ encodes the action of
its morphisms.

The same comparison result is true on the level of derived
categories: the functor
\begin{equation}\label{nerve.emb.d}
\D(\C,\Ab) \to \D(S(N(\C)),\Ab)
\end{equation}
is a fully faithful embedding, and its essential image consists of
such $M_\idot \in \B(S(N(\C)),\Ab)$ that $M_\idot(f)$ is a
quasiisomorphism whenever $f$ preserves the last elements.

\bigskip

Now, the observation that we would like make is that small
$2$-categories also have nerves. Of course, the nerve $N(\C)$ of a
$2$-category $\C$ is a simplicial category rather than a simplicial
set; however, the Grothendieck construction still applies, and the
fibered category $S(N(\C))/\Delta$ is perfectly well defined. The
only new thing is that the fibers of this fibration are no longer
discrete. As it happens, this changes nothing: if we define the
triangulated category $\D(\C,\Ab)$ by the bar-construction, as in
Subsection~\ref{2cat.subs}, then we still have a fully faithful
embedding \eqref{nerve.emb.d}, with the same characterization of its
essential image. However, since $S(N(\C))$ is a usual category, not
a $2$-category, the right-hand side $\D(S(N(\C)),\Ab)$ can be
defined in the usual way, with no recourse to $2$-categorical
machinery and $A_\infty$ methods.

\bigskip

Return now to the situation of Subsection~\ref{quo.subs} -- assume
given a small category $\C$ which has fibered products. In
principle, we could simple apply the above discussion to the
$2$-category $\Q(\C)$. However, we will actually do something
slightly different. Namely, we can define nerves in an even greater
generality: instead of a $2$-category, we can consider a ``category
in categories'', where not only {\em morphisms} form a category, not
just a set, but also {\em objects} do the same thing.

It is probably not worth the effort to axiomatize the situation;
instead, let us give the specific example we will use. Let
$\C^{[2]}$ be the category of arrows $c_1 \to c_2$ in $\C$, with
morphisms given by Cartesian squares
$$
\begin{CD}
c_1 @>>> c_2\\
@VVV @VVV\\
c_1' @>>> c'_2.
\end{CD}
$$
Consider the opposite category $\C^{[2]opp}$. We have two
projections $s$, resp. $t$ from $\C^{[2]opp}$ to $\C^{opp}$ which
send an arrow $c_1 \to c_2$ to its source $c_1$, resp. its target
$c_2$. Moreover, composition of the arrows defines a functor
$$
m:\C^{[2]opp} \times_{\C^{opp}} \C^{[2]opp} \to \C^{[2]opp},
$$
where in the left-hand side, the projection to $\C^{opp}$ is by $t$
in the left factor, and by $s$ in the right factor. This is our
``category in categories'': $\C^{opp}$ is its category of objects,
$\C^{[2]opp}$ is its category of morphisms, and the functor $m$
defines the composition.

A functor from this ``category in categories'' to $\Ab$ is, then,
given by the following data:
\begin{enumerate}
\item a functor $M_1:\C^{opp} \to \Ab$, and
\item a map $s^*M_1 \to t^*M_1$ (the analog of the action map
  \eqref{act.map}).
\end{enumerate}
We now notice that these data define exactly a functor $M$ from the
quotient category $Q\C$ of Section~\ref{mack.der} to $\Ab$. The
functor $M_1$ defines the values of $M$ at objects of the category
$\C$ and the action of the maps $f^*$, $f$ a morphism of $\C$, while
the action map $s^*M_1 \to t^*M_1$ adds the maps $f_*$ to the
picture.

This description breaks the symmetry between $f_*$ and $f^*$
manifestly present in the definition of the category $Q\C$, but this
is not necessarily a bad thing: constructions such as the tensor
product of Subsection~\ref{mack.prod.subs} also break this symmetry,
and one may hope that some constructions actually look better in the
new description. As we shall see, this is indeed the case.

We now resume rigorous exposition. We start by constructing the
nerve of our ``category in categories''; we will denote this nerve
simply by $S\C$.

\subsection{The $S$-construction.}\label{S.subs}

Assume given a small category $\C$ which has fibered products. For
every integer $n \geq 2$, let $\C^{[n]}$ be the category of diagrams
of the form $c_1 \to \dots \to c_n$ in the category $\C$ --- in
other words, $\C^{[n]}$ is the category of $\C$-valued functors from
the totally ordered set $[n]$ with $n$ elements considered as a
small category in the usual way. Sending a diagram $c_1 \to \dots
\to c_n$ to $c_n \in \C$ defines a projection $\C^{[n]} \to \C$, and
since $\C$ has fibered products, this projection is a fibration in
the sense of \cite{SGA}. Denote by $\bC^{[n]} \subset \C^{[n]}$ the
subcategory with the same objects as $\C^{[n]}$ and the maps which
are Cartesian with respect to the fibration $\C^{[n]} \to \C$. Thus
for example, $\bC^{[2]}$ is the category whose objects are arrows in
$\C$, and whose morphisms are Cartesian squares.

As usual, denote by $\Delta$ the category of non-empty finite
totally ordered sets. Then for any order-preserving map $f:[n] \to
[n']$, we have a natural functor $f^*:\bC^{[n']} \to
\bC^{[n]}$. Moreover, we can also treat it as a functor
$f^*:\bC^{[n']opp} \to \bC^{[n]opp}$ between the opposite
categories. Then the collection $\bC^{[n]opp}$ with the transition
functors $f^*$ forms a simplicial category, that is, a category
fibered over $\Delta$. We denote this category by $S\C$. 

Explicitly, objects of $S\C$ are pairs of a finite non-empty totally
ordered set $[n]$ and a diagram $c_1 \to \dots \to c_n$ in $\C$; a
map from $\langle [n],c_1 \to \dots \to c_n\rangle$ to $\langle
[n'],c'_1 \to \dots \to c'_{n'}\rangle$ is given by an
order-preserving map $f:[n] \to [n']$ and a collection of maps
$f_i:c'_{f(i)} \to c_i$ in $\C$, one for each $i \in [n]$, such that
the square
$$
\begin{CD}
c'_{f(i)} @>{f_i}>> c_i\\
@VVV @VVV\\
c'_{f(j)} @>{f_j}>> c_j
\end{CD}
$$
is commutative and Cartesian for any $i,j \in [n]$, $i \leq j$. In
particular, we have a natural embedding $\C^{opp} \to S\C$ which
sends $c \in \C^{opp}$ to $\langle [1],c \rangle$.

We can define a natural projection functor $\SQ:S\C \to Q\C$ as
follows: an object $\langle [n],c_1 \to \dots \to c_n\rangle \in
S\C$ goes to $c_n \in Q\C$, and a map $\langle f,
\{f_i\}\rangle:\langle [n],c_1 \to \dots \to c_n\rangle \to \langle
[n'],c'_1 \to \dots \to c'_{n'}\rangle$ goes to a map represented by
the diagram
$$
\begin{CD}
c_n @<{f_n}<< c'_{f(n)} @>>> c'_{n'}.
\end{CD}
$$
We will say that a map $f:[n'] \to [n]$ in the category $\Delta$ is
{\em special} if it is an isomorphism between $[n']$ and a final
segment of the ordinal $[n]$ (in other words, we have $f(l) =
n+l-n'$ for any $l \in [n']$). We will say that a map $\langle f,
\{f_i\} \rangle$ in the category $S\C$ is {\em special} if it is
Cartesian with respect to the fibration $S\C \to \Delta^{opp}$ ---
that is, all the maps $f_i$ are invertible --- and the component
$f:[n'] \to [n]$ is special.

One checks easily that if a map $f$ is special, then $\SQ(f)$ is
invertible in $Q\C$. Moreover, we will say that a functor $F:S\C \to
\Ab$ is special if $F(f)$ is invertible for every special $f$. Then
setting $F \mapsto F \circ \SQ$ gives an equivalence
$$
\Fun(Q\C,\Ab) \cong \Fun_{sp}(S\C,\Ab),
$$
where $\Fun_{sp}(S\C,\Ab) \subset \Fun(S\C,\Ab)$ is the full
subcategory spanned by special functors. In this sense, the category
$Q\C$ is obtained from the category $S\C$ by inverting all special
maps. Explicitly, a functor $M \in \Fun(S(\C),\Ab)$ is given by the
following data:
\begin{enumerate}
\item A functor $M_n \in \Fun(\C^{[n]opp},\Ab)$ for every $[n] \in
  \Delta$.
\item A transition map
\begin{equation}\label{M.trans}
(f^*)^*M_n \to M_{n'}
\end{equation}
for any map $f:[n] \to [n']$, where $f^*:\C^{[n']opp} \to
\C^{[n]opp}$ is the transition functor corresponding to the map $f$.
\end{enumerate}
The functor $M$ is special if the transition map $M(f)$ is an
isomorphism for every special map $f$ (it is clearly sufficient to
check this for the maps $f:[1] \to [n]$, $f(1) = n$).

\bigskip

\noindent
Consider now the derived category $\D(S\C,\Ab)$.

\begin{defn}\label{ds.defn}
An object $M \in \D(S\C,\Ab)$ is called {\em special} if for any
special map $f$ in the category $S\C$, the corresponding map $M(f)$
is a quasiisomorphism. The full subcategory in $\D(S\C,\Ab)$ spanned
by special complexes is denoted by $\DS(\C,\Ab) \subset
\D(S\C,\Ab)$.
\end{defn}

The triangulated category $\DS(\C,\Ab)$ obviously contains the
derived category $\D(\Fun_{sp}(S\C,\Ab)) \cong \D(Q\C,\Ab)$;
however, there is no reason why they should be same, or indeed, even
why the functor $\D(Q\C,\Ab) \subset \DS(\C,\Ab)$ should be full and
faithful. It turns out that in general, it is not. This is exactly
the difference between the naive derived category $\D(Q(\C),\Ab)$
and the triangulated category $\DQ(\C,\Ab)$ of
Subsection~\ref{quo.subs}.

\begin{theorem}\label{wald}
For any small category $\C$ which has fibered products, we have a
natural equivalence
$$
\DS(\C,\Ab) \cong \DQ(\C,\Ab).
$$
\end{theorem}

\subsection{Digression: complementary pairs.}\label{comple.subs}

To prove Theorem~\ref{wald}, we need to develop some combinatorial
machinery for inverting special maps in the category
$S\C$. Unfortunately, the class of special maps does not admit a
calculus of fractions in the usual sense. However, there is the
following substitute.

\begin{defn}\label{frac}
Assume given a category $\Phi$ and two classes of maps $P$, $I$ in
$\Phi$. Then $\langle P,I \rangle$ is a {\em complementary pair} if
the following conditions are satisfied.
\begin{enumerate}
\item The classes $P$ and $I$ are closed under the composition and
  contain all isomorphisms.
\item For every object $b \in \Phi$, the category $\Phi^I_b$ of
  diagrams $i:b' \to b$, $i \in I$ has an initial object
  $i_b:\iota(b) \to b$.
\item Every map $f$ in $\Phi$ factorizes as $f = p(f) \circ i(f)$,
  $p(f) \in P$, $i(f) \in I$, and such a factorization is unique up
  to a unique isomorphism.
\item Every diagram $b_1 \overset{p}{\gets} b \overset{i}{\to} b_2$
  in $\Phi$ with $p \in P$, $i \in I$ fits into a cocartesian square
$$
\begin{CD}
b @>{p}>> b_1\\
@V{i}VV @VV{i''}V\\
b_2 @>{p'}>> b_{12}
\end{CD}
$$
in $\Phi$ with special $p' \in P$, $i' \in I$.
\end{enumerate}
\end{defn}

\begin{remark}
If $\Phi$ has an initial object $0$, then \thetag{ii} follows from
\thetag{iii} (the map $i_c:\iota(c) \to c$ is a part of the
decomposition $0 \to \iota(c) \to c$ of the map $0 \to c$).
\end{remark}

\begin{remark}
A complementary pair in the sense of Definition~\ref{frac} is an
example of a {\em factorization system} in the sense of \cite{bou}.
\end{remark}

Assume given a small category $\Phi$ and a complementary pair
$\langle P,I \rangle$ of classes of maps in $\Phi$. Let
$\D_I(\Phi,\Ab)$ be the full subcategory in the derived category
$\D(\Phi,\Ab)$ spanned by those $M_\idot \in \D(\Phi,\Ab)$ for which
$M_\idot(i)$ is a quasiisomorphism for any $i \in I$. Let $\R\Phi$
be the category of diagrams $c_1 \gets c \to c_2$ in $\Phi$ with
maps $i_1:c \to c_1$, $i_2:c \to c_2$ in the class $I$.  Let
$\pi_i,\pi_2:\R\Phi \to \Phi$, be the projections which send $c_1
\gets c \to c_2$ to $c_1$, resp. to $c_2$, and denote by
$\Sp:\D(\Phi,\Ab) \to \D(\Phi,\Ab)$ the functor given by $\Sp =
L^\hdot\pi_{1!}\pi_2^*$.

\begin{lemma}\label{sp}
There exists a map of functors $\psi:\Id \to \Sp$ such that $\psi
\circ \Sp =\id$, and the map $\psi:M_\idot \to \Sp(M_\idot)$ for an
object $M_\idot \in \D(\Phi,\Ab)$ is a qua\-si\-isomorphism if and
only if $M_\idot$ lies in $\D_I(\Phi,\Ab) \subset \D(\Phi,\Ab)$.
\end{lemma}

In other words, $\Sp:\D(\Phi,\Ab) \to \D(\Phi,\Ab)$ is the canonical
projection onto the left-admissible full subcategory
$\D_I(\Phi,\Ab)$ --- that is, the composition of the embedding
$\D_I(\Phi,\Ab) \subset \D(\Phi,\Ab)$ and the left-adjoint functor
$\D(\Phi,\Ab) \to \D_I(\Phi,\Ab)$.

\proof{} Both projections $\pi_1$, $\pi_2$ have a common section
$\delta:\Phi \to \R\Phi$ which sends $c \in \Phi$ to the diagram $c
\gets c \to c$ with identity maps. The isomorphism $\Id \cong
\delta^*\pi_2^*$ induces by adjunction a map $L^\hdot\delta_! \to
\pi_2^*$; the functorial map $\psi$ is obtained by applying
$L^\hdot\pi_{1!}$ to this map.

To prove the required properties of the map $\psi$, we start by
noting that due to the conditions \thetag{iii} and \thetag{iv} of
Definition~\ref{frac}, the projection $\pi_1:\R\Phi \to \Phi$ is a
cofibration. Therefore for any $M_\idot \in \D(\Phi,\Ab)$ and any
object $c \in \Phi$, we have a canonical base change isomorphism
$$
\Sp(M_\idot)(c) \cong H_\idot(\R\Phi_c,\pi_2^*M_\idot),
$$
where $\R\Phi_c \subset \R\Phi$ is the fiber of the cofibration
$\pi_1$ --- that is, the category of diagrams $c \gets c' \to c''$
in $\Phi$ with maps $i':c' \to c$, $i'':c' \to c''$ lying in
$I$. This fiber $\R\Phi_c$ projects to the category $\Phi^I_c$ of
Definition~\ref{frac}~\thetag{iii} by forgetting $c''$ and $i''$;
denote by $F\Phi_c$ the fiber of this projection over the initial
object $\iota(c) \in \Phi^I_c$. Explicitly, $F\Phi_c$ is the
category of objects $c'' \in \Phi$ equipped with a map $i:\iota(c)
\to c''$, $i \in I$.

By definition, for any special map $i':c' \to c$, the map
$i_c:\iota(c) \to c$ canonically factorizes as $i_c = i' \circ
i_o:\iota(c) \to c' \to c$. Moreover, it is easy to see that
if we have two maps $i_1, i_2 \in \Phi$ such that $i_1 \in I$ and
$i_1 \circ i_2 \in I$, then $i_2 \in I$ (take its factorization $i_2
= i_2' \circ p$ of \thetag{iii}, compose with $i_1$ to obtain a
factorization of $i_1 \circ i_2$, use the uniqueness to deduce that
$p$ is invertible). This implies that $\iota(c) \cong \iota(c')$
canonically, and $i_o = i_{c'}$. Then sending a diagram $c \gets
c' \to c''$ to $i'' \circ i_{c'}:\iota(c) \to c''$ defines a
projection $\ol{\pi}:\R\Phi_c \to F\Phi_c$, and this
projection is right-adjoint to the embedding $\iota:F\Phi_c \to
\R\Phi_c$.

The projection $\pi_2$ factors as $\pi_2 = \pi \circ
\ol{\pi}:\R\Phi_c \to F\Phi_c \to \Phi$, where $\pi:F\Phi_c
\to \Phi$ sends a diagram $\iota(c) \to c''$ to
$c''$. Therefore we have
\begin{align*}
H_\idot(\R\Phi_c,\pi_2^*M_\idot) &\cong
H_\idot(\R\Phi_c,\ol{\pi}^*\pi^*M_\idot) \cong
H_\idot(\R\Phi_c,L^\hdot\iota_!\pi^*M_\idot)\\
&\cong H_\idot(F\Phi_c,\pi^*M_\idot),
\end{align*}
so that $\Sp(M_\idot)(c) \cong H_\idot(F\Phi_c,\pi^*M_\idot)$. The
canonical map $\psi:M_\idot(c) \to \Sp(M_\idot)(c)$ is then
induced by the inclusion $\ppt \to F\Phi_c$ which sends the point to
the diagram $\langle\iota(c) \to c\rangle \in F\Phi_c$.

Now, for every map $i:c \to c'$, $i \in I$, the categories $F\Phi_c$
and $F\Phi_{c'}$ are canonically equivalent; therefore the map
$\Sp(M_\idot)(c) \to \Sp(M_\idot)(c')$ is a quasiisomorphism, and we
have $\Sp(M_\idot) \in \D_I(\Phi,\Ab)$. And if we know in advance
that $M_\idot \in \D_I(\Phi,\Ab)$, then the pullback $\pi^*M_\idot
\in \D(F\Phi_C,\Ab)$ is constant, so that $\Sp(M_\idot)(c) \cong
M_\idot(c) \otimes H_\idot(F\Phi_c,\Z)$. Since $F\Phi_c$ has an
initial object, we have $H_\idot(F\Phi_c,\Z) \cong \Z$.
\endproof

For any two objects $c_1,c_2 \in \Phi$, let $\Q_I(c_1,c_2)$
be the category of diagrams $c_1 \to c \gets c_2$ in $\Phi$ such
that the map $c_1 \to c$ lies in the class $I$, and the map $c_2
\to c$ lies in the class $P$. Then for any three objects $c_1$,
$c_2$, $c_2$, we have natural composition functors
\begin{equation}\label{qq}
\Q_I(c_1,c_2) \times \Q_I(c_2,c_3) \to \Q_I(c_1,c_3)
\end{equation}
given by the cocartesian squares which exist by
Definition~\ref{frac}~\thetag{iv}. These functors are
associative. Say that an object $c \in \Phi$ is {\em simple} if the
canonical map $i_c:\iota(c) \to c$ is an isomorphism, and let
$\Q_I(\Phi)$ be the $2$-category whose objects are simple objects in
$\Phi$, and whose morphism categories are $\Q_I(c_1,c_2)$, with
composition functors \eqref{qq}. Applying the procedure of
Subsection~\ref{2cat.subs}, construct an $A_\infty$-category
$\B^I_\idot(\Phi)$ with the same objects as $\Q_I(\Phi)$, and with
morphisms given by
$$
\B^I_\idot(c_1,c_2) = C_\idot(\Q_I(c_1,c_2),\Z).
$$
Let $\DQ_I(\Phi,\Ab)$ be the triangulated category of $A_\infty$-functors
from the $A_\infty$-category $\B^I_\idot(\Phi)$ to $\Ab$.

\begin{prop}\label{comp.prop}
There exists a natural equivalence of triangulated categories
$$
\DQ_I(\Phi,\Ab) \cong \D_I(\Phi,\Ab).
$$
\end{prop}

\proof{} To define a comparison functor $\SQ:\D_I(\Phi,\Ab) \to
\DQ_I(\Phi,\Ab)$, let $\wt{\Phi}$ be $\Phi$ with a formally added
initial object $\emptyset$, and declare that the map $\emptyset \to
c$ is in the class $P$ for every $c \in \Phi$. Then all the
conditions of Definition~\ref{frac} are satisfied for $\wt{\Phi}$,
so that we can form the $2$-category $\Q_I(\wt{\Phi})$. For
any simple $c \in \Phi$, the category $\Q_I(\emptyset,c)$ is the
category $F\Phi_c$ of the proof of Lemma~\ref{sp}, and for any
simple $c_1,c_2 \in \Phi$, we have the composition functors
\begin{align*}
m_{c_1,c_2}: F\Phi_{c_1} &\times \Q_I(c_1,c_2) =
\Q_I(\emptyset,c_1) \times \Q_I(c_1,c_2) \\
&\longrightarrow F\Phi_{c_2} = \Q_I(\emptyset,c_2).
\end{align*}
A functor $M \in \Fun(\Phi,\Ab)$ gives by restriction a functor $M_c
= \pi_2^*M \in \Fun(F\Q_c,\Ab)$ for every simple $c \in \Phi$. For
every simple $c_1,c_2 \in \Phi$, we have a natural map
$$
\mu_{c_1,c_2}:M_{c_2} \boxtimes \Z \to m_{c_1,c_2}^*M_{c_1},
$$
and these maps are associative in the obvious sense. We define
$\SQ(M)_\idot \in \DQ_I(\Phi,\Ab)$ by
$$
\SQ(M)_\idot(c) = C_\idot(F\Phi_c,M_c),
$$
with the $A_\infty$-functor structure induced by the functors $m$ and
the maps $\mu$. This extends to a functor $\SQ:\D(\Phi,\Ab) \to
\DQ_I(\Phi,\Ab)$. Comparing the definitions of the functors $\Sp$
and $\SQ$, we see that the canonical map $\psi:M_\idot \to
\Sp(M_\idot)$ induces a quasiisomorphism
$$
\SQ(M_\idot) \cong \SQ(\Sp(M_\idot))
$$
for any $M_\idot \in \D(\Phi,\Ab)$, so that the functor $\SQ$
factors through the projection $\Sp:\D(\Phi,\Ab) \to
\D_I(\Phi,\Ab)$.

To show that the functor $\SQ:\D_I(\Phi,\Ab) \to \DQ_I(\Phi,\Ab)$ is
an equivalence, we have to prove that it is full, faithful, and
essentially surjective. The triangulated category $\DQ_I(\Phi,\Ab)$
is generated by the representable $A_\infty$-functors $M^Q_c$, for
all simple $c \in \Phi$ and complexes $M \in \Ho(\Ab)$, given by
$M^Q_c(c') = \B^I(c,c') \otimes M$. The triangulated category
$\DQ_I(\Phi,\Ab)$ is generated by the objects $\Sp(M_c)$, for
all $c \in \Phi$ and $M \in \Ho(\Ab)$, where $M_c \in
\Fun(\Phi,\Ab)$ is the representable functor given by $M_c(c') =
\Z[\Phi(c,c')] \otimes M$. Therefore it suffices to prove that
\begin{enumerate}
\item for any $c \in \Phi$ and $M \in \Ho(\Ab)$, we have $\Sp(M_c) \cong
  \Sp(M_{\iota(c)})$,
\item for any simple $c \in \Phi$ and $M \in \Ho(\Ab)$, we have
  $\SQ(M_c) \cong M^Q_c$, and
\item for any simple $c,c' \in \Phi$ and $M \in \Ho(\Ab)$, the
  natural map
$$
\RHom^\hdot(\Sp(M_c),\Sp(M_{c'})) \to \RHom^\hdot(M^Q_c,M^Q_{c'})
$$
induced by the functor $\SQ$ is a quasiisomorphism.
\end{enumerate}
To prove \thetag{i}, note that by adjunction, we have
$$
\RHom^\hdot(\Sp(M_c),M') \cong \RHom^\hdot(M_c,M') \cong
\RHom^\hdot(M,M'(c))
$$
for any $c \in \Phi$, $M \in \Ho(\Ab)$, and any $M' \in
\D_I(\Phi,\Ab) \subset \D(\Phi,\Ab)$. In particular,
$\RHom^\hdot(\Sp(M_c),M') \cong
\RHom^\hdot(\Sp(M_{\iota(c)},M')$. Since this is true for any $M'
\in \D_I(\Phi,\Ab)$, this implies \thetag{i}. To prove \thetag{ii},
note that for any simple $c' \in \Phi$, we have an embedding
$j:\Q_I(c,c') \to F\Phi_{c'}$, and by
Definition~\ref{frac}~\thetag{iii} there is a natural isomorphism
$$
M_c|_{F\Phi_{c'}} \cong j_!M_{const},
$$
where $M_{const} \in \Fun(\Q_I(c,c'),\Ab)$ is the constant functor
with value $M$. Therefore
$$
\Sp(M_c)(c') \cong C_\idot(F\Phi_{c'},j_!M_{const}) \cong
C_\idot(\Q_I(c,c'),M_{const}) = \B^I(c,c') \otimes M,
$$
as required. Finally, for \thetag{iii}, note that by adjunction
$$
\RHom^\hdot(\Sp(M_c),\Sp(M_{c'})) \cong \RHom^\hdot(M_c,\Sp(M_{c'}))
\cong \Sp(M_{c'})(c),
$$
so that \thetag{iii} follows from \thetag{ii}.
\endproof

\subsection{The comparison theorem.}\label{cosp.subs}

Now again, assume given a small category $\C$ which has fibered
products, and take $\Phi = S\C$. Say that a map $\langle f,
\{f_i\}\rangle$ in $S\C$ is {\em co-special} if the underlying map
$f:[n] \to [n']$ sends the first element in $[n]$ to the first
element in $[n']$, $f(1) = 1$. Let $I$, $P$ be the classes of
special, resp. co-special maps in the category $S\C$.

\begin{lemma}\label{sc.phi}
The pair $\langle P,I \rangle$ is a complementary pair in the sense
of Definition~\ref{frac}.
\end{lemma}

\proof{} Definition~\ref{frac}~\thetag{i} is obvious. For
\thetag{ii}, let $c = \langle [n], c_1 \to \dots \to c_n \rangle$;
then $\iota(c) = \langle [1], c_n \rangle$, with the obvious map
$i_c:\iota(c) \to c$. For \thetag{iii}, take a map $f = \langle f,
\{f_l\} \rangle:\langle [n], \{c_l\}\rangle \to \langle
[n'],\{c'_l\}\rangle$, and let $n_0 = f(1) \in [n']$; then the
decomposition is given by
$$
\begin{CD}
\langle [n], \{c_l\}\rangle @>>> \langle [n'-n_0], \{c''_l\}\rangle
@>>> \langle [n'],\{c'_l\}\rangle,
\end{CD}
$$
where $c''_l = c'_{l+n_0}$, $1 \leq l \leq n'-n_0$, with the obvious
maps. Finally, for \thetag{iv}, let $b = \langle
[n],\{b_l\}\rangle$, $b_1 = \langle [n_1],\{b^1_l\}\rangle$, $b_2 =
\langle [n_2],\{b^2_l\}\rangle$. Then $b_{12} = \langle [n_{12}],
\{b^{12}_l\}\rangle$ is given by $n_{12} = n_1+n_2-n$,
$b^{12}_l = b^1_{l+n-n_2}$ for $l = n-n_2+1,\dots,n_{12}$, and for
$l=1,\dots,n-n_2$, $b^{12}_l$ is obtained as the fibered product
$$
\begin{CD}
b^{12}_l @>>> b^2_l\\
@VVV @VVV\\
b^1_1 @>>> b_1 \cong b^2_{n_2+1-n}
\end{CD}
$$
in the category $\C$.
\endproof

\proof[Proof of Theorem~\ref{wald}.] By Lemma~\ref{sc.phi},
Proposition~\ref{comp.prop} can be applied to the category $S\C$. An
object $c = \langle [n], c_i \rangle$ is simple if and only if
$n=1$, thus simple objects in $S\C$ are the same as objects in $\C
\subset S\C$. Then comparing the definitions of $A_\infty$-categories
$\B^\C_\idot$ of Subsection~\ref{quo.subs} and $\B^I_\idot(S\C)$ of
Subsection~\ref{comple.subs}, we see that it remains to prove the
following: for any $c_1,c_2 \in \C \subset S\C$, there exists a
natural quasiisomorphism
$$
\eta^{c_1,c_2}:C_\idot(\Q_I(c_1,c_2),\Z) \cong
C_\idot(\Q(c_1,c_2),\Z),
$$
and these quasiisomorphisms extend to an $A_\infty$-functor. The
functor $\SQ:S\C \to Q\C$ of Subsection~\ref{quo.subs} obviously
extends to a $2$-functor $\Q_I(S\C) \to \Q\C$; this gives maps
$\eta^{c_1,c_2}$ which form an $A_\infty$-functor. To prove that
$\eta^{c_1,c_2}$ is a quasiisomorphism, it remains to notice that
the functor $\SQ:\Q_I(c_1,c_2) \to \Q\C(c_1,c_2)$ has a left-adjoint
which sends a diagram $c_1 \gets c \to c_2$ in $\C$ to the diagram
$c_1 \to \langle [2], c_2 \to c_1 \rangle \gets c_2$ in $S\C$.
\endproof

\subsection{Additivity.}

Assume now given a small category $\C$ such that the wreath product
$\C \wr \Gamma$ has fibered products. Then we can form the category
$S(\C\wr\Gamma)$ and the derived category $\DS(\C\wr\Gamma,\Ab)$. By
definition, $(\C\wr\Gamma)^{opp}$ is embedded into $S(\C\wr\Gamma)$
as the fiber over $[1] \in \Delta$, so that we have the restriction
functor $\DS(\C\wr\Gamma,\Ab) \to \D((\C\wr\Gamma)^{opp},\Ab)$.

\begin{defn}\label{add.defn.ds}
An object $M_\idot \in \DS(\C\wr\Gamma)$ is called {\em additive} if
its restriction to $\D((\C\wr\Gamma)^{opp},\Ab)$ is additive in the
sense of Definition~\ref{add.dq}.
\end{defn}

\begin{prop}\label{ds.cons}
Assume that the small category itself $\C$ has fibered
products. Then the full subcategory $\DS_{add}(\C\wr\Gamma,\Ab)
\subset \DS(\C\wr\Gamma,\Ab)$ spanned by additive objects is
canonically equivalent to the category $\DS(\C,\Ab)$.
\end{prop}

\proof{} The statement immediately follows from Theorem~\ref{wald},
Proposition~\ref{wr.nowr} and Corollary~\ref{C.wr}.
\endproof

\begin{defn}\label{ds.defn.bis}
Assume given a small category $\C$ such that $\C\wr\Gamma$ has
fibered products. Then the full subcategory $\DS_{add}(C \wr
\Gamma,\Ab) \subset \DS(C \wr \Gamma,\Ab)$ spanned by additive
objects is denoted by $\DS(\C,\Ab)$.
\end{defn}

We note that by Proposition~\ref{ds.cons}, this is consistent with
our earlier Definition~\ref{ds.defn}.

\begin{remark}
For the sake of methodological purity, it would be nice to have a
direct proof of Proposition~\ref{ds.cons} which does not use the
material of Section~\ref{mack.der}. Unfortunately, I was not able to
find such a proof. All I could come up with essentially repeats the
proof of Proposition~\ref{wr.nowr}, with additional complications
(which arise because we do not have a Waldhausen-type interpretation
of the category $\DQ^{\wr}(\C,\Ab)$ of Definition~\ref{dq.wr.defn}).
\end{remark}

In the particular case $\C = O_G$, the category of finite $G$-orbits
for a group $G$, we have $\C \wr \Gamma \cong \Gamma_G$. Combining
Theorem~\ref{wald} and Proposition~\ref{ds.cons}, we get the
following.

\begin{corr}
For any finite group $G$, let $O_G$ be the category of
$G$-orbits. Then the category $\DM(G)$ of derived $G$-Mackey
functors is naturally equivalent to the category $\DS(O_G,\Ab)$ of
Definition~\ref{ds.defn.bis}.\endproof
\end{corr}

\section{Functoriality and products.}\label{pro.sec}

We will now describe some basic properties of derived Mackey functors,
mostly analogous to the material in Subsection~\ref{mack.func.subs}
and Subsection~\ref{mack.prod.subs}.

\subsection{Functoriality.}\label{der.func}

Assume given two small categories $\C$, $\C'$ such that
$\C\wr\Gamma$ and $\C'\wr\Gamma$ have fibered products. Then any
functor $\gamma:\C' \to \C\wr\Gamma$ uniquely extends to a
coproduct-preserving functor $\gamma:\C'\wr\Gamma \to
\C\wr\Gamma$. In either of the two constructions of the category
$\DS(-,\Ab) \cong \DQ_{add}(-\wr\Gamma,\Ab)$ we have the following
obvious functoriality property.
\begin{itemize}
\item If the extended functor $\gamma:\C'\wr\Gamma \to \C\wr\Gamma$
preserves fibered products, then we have a natural restriction
functor
$$
\gamma^*:\DS(\C,\Ab) \cong \DQ_{add}(\C\wr\Gamma,\Ab) \to
\DS(\C',\Ab) \cong \DQ_{add}(\C'\wr\Gamma,\Ab),
$$
and the left-adjoint induction functor $\gamma_!:\DS(\C',\Ab) \to
\DS(\C,\Ab)$.
\end{itemize}
It turns out that this yields the derived versions of both the
functor $\Psi$ and the functor $\Phi$ of
Subsection~\ref{mack.func.subs}.

In fact, the functor $\Psi$ has already appeared in
Section~\ref{mack.der} under a different name. For any object $c \in
\C$, we have its embedding functor $j^c:\ppt \to \C\wr\Gamma$,
and the extended functor $j^c:\Gamma \to \C\wr\Gamma$ preserves
fibered products for semi-trivial reasons.

\begin{defn}
The {\em naive fixed point functor} $\Psi^c$ is the functor
$$
\Psi^c = j^{c*}:\DS(\C,\Ab) \to \DS(\ppt,\Ab) \cong \D(\Ab).
$$
\end{defn}

In the Mackey functor case $\C = O_G$, $c = [G/H]$, the functor
$\Psi^c$ is the derived version of the functor $\Psi^H$ of
Subsection~\ref{mack.func.subs}.

We can also slightly refine the construction. Let $\br{c} \subset
\C$ be the groupoid of objects in $\C$ isomorphic to $c$, and
invertible maps between them. Then $j^c$ extends to an embedding
$\wt{j}^c:\br{c} \to \C\wr\Gamma$ whose natural extension
$\wt{j}^c:\br{c}\wr\Gamma \to \C\wr\Gamma$ still preserves fibered
products. We denote the corresponding restriction functor by
\begin{equation}\label{wpsi.eq}
\wPsi^c = \wt{j}^{c*}:\DS(\C,\Ab) \to \DS(\br{c},\Ab) \cong
\D(\Aut(c),\Ab).
\end{equation}
It takes values in the category $\DS(\br{c},\Ab)$, which is
obviously equivalent to the derived category $\D(\Aut(c),\Ab)$ of
$\Ab$-valued representations of the group $\Aut(c)$ of automorphisms
of the object $c$.

\medskip

To proceed further, we need to impose a restriction on the category
$\C$.

\begin{defn}
A category $\C$ is called {\em $\Hom$-finite} if the set of maps
$\C(c,c')$ is finite for any two objects $c,c' \in \C$.
\end{defn}

Assume that the category $\C$ is $\Hom$-finite. Then so is the
wreath product $\C\wr\Gamma$, and an object $c$ also represents a
functor $\tau_c:\C \to \Gamma$,
$$
\tau_c(c') = \C(c,c').
$$
Its natural extension $\tau_c:\C\wr\Gamma \to \Gamma$ is represented
by the same object $c$, thus preserves fibered products.

\begin{defn}\label{phi.defn}
The {\em geometric fixed points functor} $\Phi^c$ is given by
$$
\Phi^c = \tau_{c!}:\DS(\C,\Ab) \to \DS(\ppt,\Ab) \cong \D(\Ab).
$$
\end{defn}

In the case $\C = O_G$, $c = [G/H]$, this is the derived version of
the functor $\Phi^H$ of Subsection~\ref{mack.func.subs}. The adjoint
functor $\tau_c^*:\D(\Ab) \to \DS(\C,\Ab)$ can be explicitly
described as as follow. Denote
$$
T^c = \tau_c^*(T) \in \Fun_{sp}(S(\C\wr\Gamma),\Z\amod) \subset
\DS(\C,\Z\amod),
$$
where $T \in \Fun_{sp}(S(\Gamma),\Z\amod)$ is the natural generator
of the category $\DQ_{add}(\Gamma,\Z\amod)$ introduced in
Lemma~\ref{DQ.gamma}. Then for any $M_\idot \in \D(\Ab)$, we have
$$
\tau_c^*(M_\idot) \cong M_\idot \otimes T^c.
$$
This can also be refined to incorporate $\Aut(c)$. Indeed, the group
$\Aut(c)$ obviously acts on the object $T^c$. We define the {\em
inflation functor} 
$$
\Infl^c:\Fun(\Aut(c),\Ab) \to \Fun(S(C\wr\Gamma),\Ab)
$$
by $\Infl^c(M) = \left(M \otimes T^c\right)^{\Aut(c)}$, and we note
that the functor $R^\hdot\Infl^c$ sends $\D(\Aut(c),\Ab)$ into the
category $\DS(\C,\Ab) \subset \D(S(\C\wr\Gamma),\Ab)$. By abuse of
notation, we will drop $R^\hdot$ and denote the derived functor
simply by
$$
\Infl^c:\D(\Aut(c),\A) \cong \D(\br{c},\Ab) \to \DS(\C,\Ab),
$$
and we denote by
$$
\wPhi^c:\DS(\C,\Ab) \to \D(\Aut(c),\Ab) \cong \D(\br{c},\Ab)
$$
its left-adjoint.

\begin{remark}
The names {\em naive} and {\em geometric} attached to fixed points
functors come from equivariant stable homotopy theory, whose part in
the story we will explain in Section~\ref{spectra}. It seems that
the geometric fixed point functor is the more important of the two;
thus from now, ``fixed point functor'' without an adjective will
mean the functor $\Phi^c$.
\end{remark}

We can further refine this construction by the following
observation: the group $\Aut(c)$ acts not only on the functor $T^c$,
but on the functor $\tau_c$, too -- for any $c' \in \C\wr\Gamma$,
$\Aut(c)$ naturally act on the finite set $\tau_c(c) =
\C(c,c')$. Thus we actually have a functor $\wh{\tau}_c:\C \to
\Gamma_{\Aut(c)} = O_{\Aut(c)}\wr\Gamma$ which induces functors
\begin{align*}
\hPhi^c &= \wh{\tau}_{c!}:\DS(\C,\Ab) \to \DS(O_\Aut(c),\Ab) =
\DM(\Aut(c),\Ab),\\
\wh{\Infl}^c &= \wh{\tau}_c^*:\DM(\Aut(c),\Ab) \to \DS(\C,\Ab).
\end{align*}
In the case $\C = O_G$, $c = [G/H]$, we have $\Aut(c) = N_H/H$,
where $N_H \subset G$ is the normalizer of the subgroup $G$. Then
the extended inflation functor takes the form
$$
\wh{\Infl}^N_G:\DM(N,\Ab) \to \DM(G,\Ab),
$$
where $N = N_H/H$. If $H \subset G$ is normal, so that $N_H = G$ and
$N = G/H$, this is the derived version of the fully faithful functor
$\Infl^N_G$ of Subsection~\ref{mack.func.subs}. As we will prove in
Section~\ref{tate.sec}, Lemma~\ref{infl.ff}, in this case the
functor $\wh{\Infl}^N_G$ is fully faithful.

\subsection{Products.}\label{der.prod}

To introduce a tensor product on the category $\DS(\C,\Ab)$, we
impose the same additional assumption as in
Subsection~\ref{mack.prod.subs} -- we require that the small
category $\C$ has a terminal object. Then so does $\C\wr\Gamma$; of
course, we still assume that $\C\wr\Gamma$ has fibered products, so
that altogether, it has all finite limits (in particular, products).

It is convenient to identify $\DS(\C,\Ab) \cong
\DQ_{add}(\C\wr\Gamma,\Ab)$ and use the $A_\infty$ methods of
Section~\ref{mack.der}. By the definition of the $2$-category
$\Q(\C\wr\Gamma)$, we have an isomorphism
$$
\Q((\C\wr\Gamma) \times (\C\wr\Gamma)) \cong \Q(\C\wr\Gamma) \times
\Q(\C\wr\Gamma).
$$
Since the bar construction is compatible with products, this means
that for any two objects $M_\idot,M'_\idot \in
\DQ(\C\wr\Gamma,\Ab)$, we have a well-defined box product
$$
M_\idot \boxtimes M'_\idot \in \DQ((\C\wr\Gamma)\times(\C\wr\Gamma),\Ab)
$$
We also have the product functor
$$
m:(C\wr\Gamma) \times (C\wr\Gamma) \to C\wr\Gamma,
$$
and since it preserves fibered products, it induces the restriction
functor
\begin{align*}
m^*:\DS(\C,\Ab) &\cong \DQ_{add}(\C\wr\Gamma,\Ab) \subset
\DQ(\C\wr\Gamma,\Ab) \\
&\longrightarrow \DQ((\C\wr\Gamma)\times(\C\wr\Gamma),\Ab).
\end{align*}
Consider the left-adjoint functor
$$
m_!^{add}:\DQ((\C\wr\Gamma)\times(\C\wr\Gamma),\Ab) \to \DS(\C,\Ab).
$$

\begin{defn}\label{prod.defn}
The {\em tensor product} $M_\idot \otimes M'_\idot$ of two objects
$M_\idot,M'_\idot \in \DS(\C,\Ab)$ is given by
$$
M_\idot \otimes M'_\idot = m_!^{add}(M_\idot \boxtimes M'_\idot).
$$
\end{defn}

Under our assumptions on $\Ab$, this obviously gives a well-defined
symmetric tensor product structure on the triangulated category
$\DS(\C,\Ab)$.

\medskip

We have the following result on compatibility between tensor
products and fixed points. Assume that the small category $\C$ is
$\Hom$-finite, so that for any object $c \in \C$, we have a
well-defined geometric fixed points functor $\Phi^c:\DS(\C,\Ab) \to
\D(\Ab)$.

\begin{prop}\label{phi.tensor}
For any $c \in \C$, the geometric fixed points functor $\Phi^c$ is a
tensor functor -- that is, for any two objects $M_\idot,M'_\idot \in
\DS(\C,\Ab)$ we have an isomorphism
$$
\Phi^c(M_\idot \otimes M'_\idot) \cong \Phi^c(M_\idot) \otimes
\Phi^c(M'_\idot),
$$
and this isomorphism is functorial in $M_\idot$ and $M'_\idot$.
\end{prop}

\proof{} The functor $\tau_c:\C\wr\Gamma \to \Gamma$, being
representable, commutes with products, so that we have a commutative
diagram
$$
\begin{CD}
\Q(\C\wr\Gamma) \times \Q(\C\wr\Gamma) @>{m}>> \Q(\C\wr\Gamma)\\
@V{\tau_c \times \tau_c}VV @VV{\tau_c}V\\
\Q(\Gamma) \times \Q(\Gamma) @>{m}>> \Q(\Gamma),
\end{CD}
$$
and $m^* \circ \tau_c^* \cong (\tau_c \times \tau_c)^* \circ
m^*$. By adjunction, this gives a functorial isomorphism
$$
\Phi^c(M_\idot \otimes M'_\idot) \cong \Phi^c(M_\idot) \otimes
\Phi^c(M'_\idot),
$$
where the product in the right-hand side is taken in the category
$\DS(\ppt,\Ab)$ in the sense of Definition~\ref{prod.defn}. It
remains to notice that the canonical equivalence of
Lemma~\ref{DQ.gamma} identifies this product with the tensor
product in the derived category $\D(\Ab)$.
\endproof

\begin{remark}
The situation with the tensor product in $\DQ(\C,\Ab)$ is somewhat
reminiscent of the tensor product of $\D$-modules on a smooth
algebraic variety $X$: while a $\D$-modules is simply a sheaf of
modules over the ring $\D_X$ of differential operators on $X$, this
ring is not commutative, and one has to take special care to define
a symmetric tensor product on the category $\D_X\amod$. The product
becomes much more natural if one passes to a Koszul-dual
interpretation and replaces $\D$-modules with sheaves of DG modules
over the de Rham complex $\Omega^\hdot_X$. We will obtain an
analogous Koszul-dual interpretation of Mackey functors in
Section~\ref{galois.sec}.
\end{remark}

\subsection{Induction.}\label{der.indu}

We finish this section with one more functoriality result. It seems
that it does not appear in the standard theory of Mackey functors;
however, it will be very useful in Section~\ref{spectra}.

We again fix a small category $\C$ such that $\C\wr\Gamma$ has
fibered products. By the definition of the $2$-category
$\Q(\C\wr\Gamma)$, we have a natural functor $\C\wr\Gamma \to
\Q(\C\wr\Gamma)$ which is identical on objects, and sends a map
$f:S' \to S$ to the diagram $S' \overset{\id}{\gets} S'
\overset{f}{\to} S$. Composing this with the natural embedding
$j^\C:\C \to \C\wr\Gamma$, we obtain a functor
$$
q:\C \to \Q(\C\wr\Gamma).
$$
Since the $2$-category $\Q(\C\wr\Gamma)$ is manifestly self-dual, we
also have the opposite functor $q^{opp}:\C^{opp} \to \Q\wr\Gamma$,
and the corresponding restriction functor
$$
q^{opp *}:\DS(\C,\Ab) \cong \DQ_{add}(\C\wr\Gamma,\Ab) \subset
\DQ(\C\wr\Gamma,\Ab) \to \D(\C^{opp},\Ab).
$$
We will call the adjoint functor
$$
q^{opp}_!:\D(\C^{opp},\Ab) \to \DS(\C,\Ab)
$$
the {\em induction functor}, and we will say that an object $M_\idot
\in \DS(\C,\Ab)$ is {\em induced} if $M_\idot \cong
q^{opp}_!(E_\idot)$ for some $E_\idot \in \D(\C^{opp},\Ab)$.

If the target category $\Ab$ is a symmetric tensor category, then
the functor category $\Fun(\C^{opp},\Ab)$ has a natural
``pointwise'' tensor product that induces a tensor product on the
derived category $\D(\C^{opp},\Ab)$.

\begin{prop}\label{indu.prod}
\begin{enumerate}
\item Assume that the category $\C$ is $\Hom$-finite. Then for any
  object $c \in \C$ and any $E_\idot \in \D(\C^{opp},\Ab)$, we have
  an isomorphism
$$
\Phi^c(q^{opp}_!(E_\idot)) \cong E_\idot(c),
$$
and this isomorphism is functorial in $E_\idot$.
\item Assume that the category $\C$ has a terminal object and that
  the target category $\Ab$ is a symmetric tensor category. Then the
  induction functor $q^{opp}_!$ is a tensor functor -- we have an
  isomorphism
$$
q^{opp}_!(E_\idot \otimes E_\idot') \cong q^{opp}_!E_\idot \otimes
q^{opp}_!E_\idot'
$$
for any $E_\idot,E'_\idot \in \D(\C^{opp},\Ab)$, and this
isomorphism is functorial in $E_\idot$ and $E'_\idot$.
\end{enumerate}
\end{prop}

\proof{} For \thetag{i}, note that by adjunction, it suffices to
prove that $q^{opp *}T^c \in \Fun(\C^{opp},\Z\amod)$ is the object
co-represented by $c \in \C^{opp}$, that is,
$$
T^c(c') \cong \Z[\C^{opp}(c',c)] = \Z[\C(c,c')].
$$
This is the definition of the object $T^c$.

For \thetag{ii}, we repeat the argument of
Subsection~\ref{mack.prod.subs}. Firstly, as noted several times
already, we have a natural equivalence between $\D(\C^{opp},\Ab)$
and the full subcategory $\D_{add}((\C\wr\Gamma)^{opp},\Ab) \subset
\D((\C\wr\Gamma)^{opp},\Ab)$ spanned by additive objects. Secondly,
we have a commutative diagram
$$
\begin{CD}
(\C\wr\Gamma)^{opp} \times (\C\wr\Gamma)^{opp} @>{m}>>
  (\C\wr\Gamma)^{opp}\\
@V{q^{opp}}VV @VV{q^{opp}}V\\
\Q(\C\wr\Gamma) \times \Q(\C\wr\Gamma) @>{m}>> \Q(\C\wr\Gamma).
\end{CD}
$$
Therefore it suffices to prove that for any $E_\idot,E'_\idot \in
\D(\C^{opp},\Ab) \cong \D_{add}((\C\wr\Gamma)^{opp},\Ab)$ we have a
functorial isomorphism
$$
m_!(E_\idot \boxtimes E'_\idot) \cong E_\idot \otimes E'_\idot.
$$
This is obvious: the product functor $m:(\C\wr\Gamma)^{opp} \times
(\C\wr\Gamma)^{opp} \to (\C\wr\Gamma)^{opp}$ is left-adjoint to the
diagonal embedding $\delta:(\C\wr\Gamma)^{opp} \to
(\C\wr\Gamma)^{opp} \times (\C\wr\Gamma)^{opp}$, so that $m_! \cong
\delta^*$, and by definition, $E_\idot \otimes E'_\idot \cong
\delta^*(E_\idot \boxtimes E'_\idot)$.
\endproof

\section{Categories of Galois type.}\label{galois.sec}

We will now give yet another description of the category
$\DS(\C,\Ab)$ of Definition~\ref{ds.defn.bis}, to complement those
given in Section~\ref{mack.der} and Section~\ref{wald.sec}. It is
based on an explicit DG model, as in Section~\ref{mack.der};
however, this new DG model is more economical and more convenient
for applications. It will require some additional assumptions on
$\C$ (which are satisfied in the Mackey functor case $\C = O_G$). In
addition, from now we will assume that the target category $\Ab$ is
the category of modules over a ring, so that Lemma~\ref{weib}
applies.

\subsection{Galois-type categories and fixed points.}

We begin by imposing our conditions on the small category $\C$.

\begin{defn}
A category $\C$ is {\em lattice-like} if it is $\Hom$-finite, and
all its morphisms are surjective.
\end{defn}

\begin{lemma}\label{lattice}
In a lattice-like category $\C$, every right-inverse $f':c \to c'$
to a morphism $f:c' \to c$ is also a left-inverse. Moreover, every
endomorphism $f:c \to c$ of an object $c \in \C$ is invertible.
\end{lemma}

\proof{} For the first claim, note that $(f' \circ f) \circ f' = f'
\circ (f \circ f') = f' \circ \id = \id \circ f'$; since $f'$ is
surjective, this implies $f' \circ f = \id$. For the second claim,
note that since $f$ is surjective, the natural map
$$
\begin{CD}
\C(c,c) @>{- \circ f}>> \C(c,c)
\end{CD}
$$
is injective. Since $\C(c,c)$ is a finite set, this map must also be
surjective, so that there exists $f' \in \C(c,c)$ such that $f'
\circ f = \id$.
\endproof

\begin{defn}
A category $\C$ is {\em of Galois type} if it is lattice-like, and
the wreath product $\C\wr\Gamma$ has fibered products.
\end{defn}

\begin{exa} Here are some examples of categories of Galois type.
\begin{enumerate}
\item $\C=O_G$, the category of finite $G$-orbits for a group $G$.
\item $\C$ is the category of finite separable extensions of a field
  $k$.
\item $\C$ is the category $\Gamma_-$ of finite sets and surjective
  maps between them.
\end{enumerate}
Of course, Galois theory shows that \thetag{ii} is a particular case
of \thetag{i}, thus the name ``Galois-type''; however, the very
interesting example \thetag{iii} is not of this form.
\end{exa}

Assume given a small category $\C$ of Galois type, and consider the
category $\DS(\C,\Ab)$ of Definition~\ref{ds.defn.bis}. Since $\C$
is $\Hom$-finite, we can apply the constructions of
Subsection~\ref{der.func}; thus for any $c \in \C$, we have the
object $T^c \in \DS(\C,\Z\amod)$, the inflation functor
$$
\Infl^c:\D(\br{c},\Ab) \to \DS(\C,\Ab)
$$
and the left-adjoint fixed point functor
\begin{equation}\label{wphi.eq}
\wPhi^c:\DS(\C,\Ab) \to \D(\br{c},\Ab).
\end{equation}
In general, computing the functors $\wPhi^c$ explicitly seems to be
rather difficult; however, in the case of Galois-type categories,
there is a drastic simplification.  Namely, let $q:\C \to
\C\wr\Gamma \to \Q(\C)$ be the natural embedding of
Subsection~\ref{der.indu}, and let $q^*:\DS(\C,\Ab) \to \D(\C,\Ab)$
be the corresponding restriction functor.  Moreover, define an
embedding
$$
\nu_c:\D(\br{c},\Ab) \to \D(\C,\Ab)
$$
by
$$
\nu_c(M_\idot)([c']) =
\begin{cases}
M_\idot(c'), &\quad c' \in \br{c},\\
0, &\quad \text{otherwise},
\end{cases}
$$
and let $\wphi^c:\D(\C,\Ab) \to \D(\br{c},\Ab)$ be the
left-adjoint functor. Then we have the following.

\begin{prop}\label{nothing.to.do}
We have a functorial isomorphism
$$
\wPhi^c \cong \wphi^c \circ q^*.
$$
\end{prop}

In order to prove this proposition, we need some preliminary results
and constructions.

\medskip

Consider the simplicial category $S(\C\wr\Gamma)/\Delta$. Note that
by the definition of the wreath product, we have a natural fibration
$S(\C\wr\Gamma) \to \Gamma$, so that altogether, $S(\C\wr\Delta)$ is
fibered over $\Delta \times \Gamma$. Let $\ol{S(\C\wr\Gamma)}
\subset S(\C\wr\Gamma)$ be the subcategory with the same objects,
and those morphisms which are Cartesian with respect to the
fibration $S(\C\wr\Gamma) \to \Delta \times \Gamma$. By definition,
$\ol{S(\C\wr\Gamma)}$ is also fibered over $\Delta$, and we have a
natural Cartesian embedding
\begin{equation}\label{eta.emb}
\eta:\ol{S(\C\wr\Gamma)} \to S(\C\wr\Gamma).
\end{equation}
While the category $\ol{S(\C\wr\Gamma)}$ is not of the form $S(\C')$
for some small category $\C'$, most of the material of
Section~\ref{wald.sec} still applies to it, with appropriate
changes. In particular, say that a map $f$ in $\ol{S(\C\wr\Gamma)}$
is {\em special}, resp. {\em co-special} if so is $\eta(f)$. Then it
is easy to check that the classes of special and co-special maps
form a complementary pair in the sense of Definition~\ref{frac}. We
can thus consider the full subcategory
$$
\ol{\DS}^\wr(\C,\Ab) \subset \D(\ol{S(\C\wr\Gamma)},\Ab)
$$
spanned by special functors. Since $\eta$ is a Cartesian functor
which respects special maps, we have a restriction functor
$\eta^*:\DS(\C\wr\Gamma,\Ab) \to \ol{\DS}^\wr(\C,\Ab)$, and the
adjoint functor 
$$
R^\hdot\eta_*:\D(\ol{S(\C\wr\Gamma)},\Ab) \to \D(S(\C\wr\Gamma),\Ab)
$$
sends special functors into special functors.

Moreover, we also have a $2$-category description of
$\ol{DS}^\wr(\C,\Ab)$. Namely, the fiber $\ol{S(\C\wr\Gamma)}_1$
over $[1] \in \Delta$ is obviously equivalent to
$(\bC\wr\Gamma)^{opp}$ (where $\bC \subset \C$ is the subcategory
with the same objects and invertible maps between them). Consider
the $2$-category $\Q(\C\wr\Gamma)$ of Section~\ref{mack.der}, and
let $\ol{\Q}(\C\wr\Gamma) \subset \Q(\C\wr\Gamma)$ be the
$2$-subcategory with the same objects, those $1$-morphisms $c_1
\gets c \to c_2$ for which the map $c \to c_1$ actually lies in
$\bC\wr\Gamma \subset \C\wr\Gamma$, and all $2$-morphisms between
them. One checks easily that this condition is compatible with the
pullbacks, so that we indeed have a well-defined
$2$-category. Applying the machinery of Subsection~\ref{2cat.subs},
we produce the $A_\infty$-category $\ol{\B}_\idot$ and the derived
category $\ol{\DQ}(\C\wr\Gamma,\Ab)$ of $A_\infty$-functors from
$\ol{\B}_\idot$ to $\Ab$. Note that we have a natural functor
$(\bC\wr\Gamma)^{opp} \to \ol{\Q}(\C\wr\Gamma)$, thus a restriction
functor $\ol{\DQ}(\C\wr\Gamma,\Ab) \to
\D((\bC\wr\Gamma)^{opp},\Ab)$. By definition, we also have a natural
embedding $\eta:\ol{\Q}(\C\wr\Gamma) \to \Q(\C\wr\Gamma)$.

\begin{lemma}\label{wald.wrr}
There exists a natural equivalence
$$
\ol{\DS}^\wr(\C,\Ab) \cong \ol{\DQ}(\C\wr\Gamma,\Ab)
$$
which is compatible with $\eta^*$ and with the natural restrictions
to $(\bC\wr\Gamma)^{opp}$.
\end{lemma}

\proof{} Same as Theorem~\ref{wald}.
\endproof

Now fix an object $c \in \C$, and define a $2$-functor
$\eps_c:\ol{\Q}(\C\wr\Gamma) \to \Q(\br{c}\wr\Gamma)$ as follows. On
objects, $\eps_c$ sends $\langle S,\{c_s\} \rangle \in \C\wr\Gamma$
to the formal union of those components $c_s$ which are isomorphic
to $c$. On morphisms, $\eps_c$ sends a diagram $c_1 \gets c' \to c_2$
to $\eps_c(c_1) \gets \eps_c(c') \to \eps_c(c_2)$ (if $c_1,c_2 \in
\br{c}\wr\Gamma$, this is well-defined). Moreover, let
$$
\overline{T}^c = \eps_c^*(T^c_{\br{c}}) \in
\ol{\DQ}(\C\wr\Gamma,\Z\amod)
$$
be the pullback of the standard object $T^c_{\br{c}} \in
\ol{\DQ}(\br{c}\wr\Gamma,\Z\amod)$, $T^c_{\br{c}}(c') =
\Z[(\br{c}\wr\Gamma)(c,c')]$, and define an embedding
$\nu^\wr_c:\D(\br{c},\Ab) \to \ol{\DQ}(\C\wr\Gamma,\Ab)$ by
$$
\nu^\wr_c(M_\idot) = (\ol{T}^c \otimes M_\idot)^{\Aut(c)}
$$
(since $\ol{T}^c(c')$ is a free $\Aut(c)$-module for any $c' \in
\C$, the group $\Aut(c)$ has no higher cohomology with coefficients
in $\ol{T}^c \otimes M_\idot$, so that taking $\Aut(c)$-invariants is
well-defined on the level of derived categories).

\begin{lemma}\label{infla}
We have a natural isomorphism
$$
T^c \cong R^\hdot\eta_*\ol{T}^c \in \DS(\C\wr\Gamma,\Ab) \cong
\DQ(\C\wr\Gamma,\Ab).
$$
\end{lemma}

\proof{} To construct a map $T^c \to R^\hdot\eta_*\overline{T}^c$,
it suffices by adjunction to construct a map $\eta^*T^c \to
\overline{T}^c$. Such a map is obvious -- it is identical on $c' \in
\br{c}\wr\Gamma \subset \C\wr\Gamma$, and $0$ otherwise. To prove
that the induced map $T^c \to R^\hdot\eta_*\ol{T}^c$ is an
isomorphism, use Theorem~\ref{wald} and Lemma~\ref{wald.wrr} to pass
to the Waldhausen-type interpretation. Then since both $T^c$ and
$R^\hdot\eta_*\ol{T}^c$ are special, it suffices to check that $T^c
\to R^\hdot\eta_*\ol{T}^c$ is an isomorphism on the fiber
$(\C\wr\Gamma)^{opp}$ over $[1] \in \Delta$. This is again obvious:
on this fiber, $T^c$ is the functor represented by $c \in
(\C\wr\Gamma)^{opp}$, and $\overline{T}^c$ is the functor
represented by the same $c$ considered as an object in
$\ol{S(\C\wr\Gamma)}_1 = (\bC\wr\Gamma)^{opp} \subset
(\C\wr\Gamma)^{opp}$.
\endproof

It remains to analyse the $2$-category $\ol{\Q}(\C\wr\Gamma)$. As in
Subsection~\ref{gamma.subs}, consider the full subcategory in
$\ol{\Q}(\C\wr\Gamma)$ whose $1$-morphisms are diagrams $S \gets S_1
\to S'$ with injective map $S_1 \to S$. This is actually a
$1$-category which we denote by $\C\wr\Gamma_+$ by abuse of
notation. The embedding $q:\C \to \Q(\C\wr\Gamma)$ factors through
the embedding $\eta:\ol{\Q}(\C\wr\Gamma)$ by means of the embeddings
$$
\begin{CD}
\C @>{j}>> \C\wr\Gamma_+ @>{\lambda}>> \ol{\Q}(\C\wr\Gamma),
\end{CD}
$$
and we have and the corresponding restriction functors $\lambda^*$,
$j^*$. Say that an object $M_\idot \in \ol{\DS}^\wr(\C,\Ab) \cong
\ol{\DQ}(\C\wr\Gamma,\Ab)$ is {\em additive} if its restriction to
$\bC\wr\Gamma^{opp}$ is additive in the sense of
Definition~\ref{add.dq}.

\begin{lemma}\label{addi.wr}
Assume given two objects $M_\idot,M'_\idot \in
\ol{\DQ}(\C\wr\Gamma,\Ab)$. If $M'_\idot$ is additive, then the
natural map
\begin{equation}\label{addi.wr.eq}
\Hom(M_\idot,M'_\idot) \to \Hom(\lambda^*M_\idot,\lambda^*M_\idot')
\end{equation}
is an isomorphism. If $M_\idot$ is also additive, then the natural
map
$$
\Hom(\lambda^*M_\idot,\lambda^*M_\idot') \to
\Hom(j^*\lambda^*M_\idot,j^*\lambda^*M_\idot')
$$
is also an isomorphism.
\end{lemma}

\proof{} As in the proof of Lemma~\ref{DQ.gamma}, we first note that
the category $\Fun(\C\wr\Gamma_+,\Z\amod)$ is generated by
representable functors $\Z^S$, $S \in \C\wr\Gamma_+$. For any such $S
= \coprod c_s$, $c_s \in \C$, we have
$$
\Z^S \cong \bigotimes \Z^{c_s},
$$
and for any $c \in \C$, we have a canonical direct sum decomposition
$\Z^c = T^c \oplus \Z$. This induces direct sum decompositions
$$
\Z^S = \bigoplus_{S' \subset S}T^{S'},
$$
where the objects
$$
T^S = \bigotimes T^{c_s}, \qquad S = \coprod c_s \in \C\wr\Gamma_+
$$
give a smaller set of projective generators of the category
$\Fun(\C\wr\Gamma_+,\Z\amod)$. For any $S_1,S_2 \in \C\wr\Gamma_+$,
we have $\Z^{S_1 \coprod S_2} \cong \Z^{S_1} \otimes \Z^{S_2}$;
therefore the additivity condition on an object $E_\idot \in
\D(\C\wr\Gamma_+,\Ab)$ means that the natural map $\Z^{S_1} \oplus
\Z^{S_2} \to \Z^{S_1 \coprod S_2} = \Z^{S_1} \otimes \Z^{S_2}$
induces an isomorphism
$$
\RHom^\hdot(\Z^{S_1} \otimes \Z^{S_2} \otimes M_\idot,E_\idot) \to
\RHom^\hdot(\Z^{S_1} \otimes M_\idot,E_\idot)
\oplus
\RHom^\hdot(\Z^{S_2} \otimes M_\idot,E_\idot)
$$
for any $M_\idot \in \D(\Ab)$. In terms of the generators $T^S$,
this is equivalent to saying that $\RHom^\hdot(T^S \otimes
M_\idot,E_\idot)=0$ unless the decomposition $S = \coprod c_s$ has
exactly one term, $S = c$ for some $c \in \C \subset \C\wr\Gamma_+$:
\begin{equation}\label{addi.wr.eq2}
\RHom^\hdot(T^S \otimes M_\idot,E_\idot) =
\begin{cases}
\RHom^\hdot(M,E_\idot(c)), &\quad S = c \in \C \subset
\C\wr\Gamma_+,\\
0, &\quad \text{otherwise}.
\end{cases}
\end{equation}
(If $\C = \{\ppt\}$, as in Lemma~\ref{DQ.gamma}, these are exactly the
orthogonality conditions on the generators $T^{\otimes n}$.) Then as
in the proof of Lemma~\ref{DQ.gamma}, it suffices to check the first
claim for a representable $M_\idot = M^S_\idot \in \DS(\C,\Ab)$; we
have a filtration on $\Hom(M^S_\idot,M'_\idot)$ with associated
graded quotient of the form \eqref{gr.Tn}, and checking
that the map \eqref{addi.wr.eq} is an isomorphism amounts to
applying \eqref{addi.wr.eq2} to $E_\idot = \lambda^*M'_\idot$. To prove the
second claim, it remains to notice that for any
additive object $E_\idot \in \D(\C\wr\Gamma_+,\Ab)$, the adjunction
map $L^\hdot j_!j^*E_\idot \to E_\idot$ is an isomorphism.
\endproof

\proof[Proof of Proposition~\ref{nothing.to.do}.] By definition, for
any $M_\idot \in \DS(\C,\Ab)$, $M'_\idot \in \D(\br{c},\Ab)$ we have
$$
\Hom(\wPhi^cM_\idot,M'_\idot) \cong
\Hom(M_\idot,\Infl^c(M'_\idot)),
$$
and by Lemma~\ref{infla}, we have
$$
\Infl^c(M'_\idot) \cong R^\hdot\eta_*(\ol{T}^c \otimes
M'_\idot)^{\Aut(c)} = R^\hdot\eta_*\nu^\wr_c(M'_\idot),
$$
so that
$$
\Hom(\wPhi^cM_\idot,M'_\idot) \cong
\Hom(\eta^*M_\idot,\nu^\wr_c(M'_\idot)).
$$
It remains to apply Lemma~\ref{wald.wrr} and Lemma~\ref{addi.wr},
and notice that by definition, $q^*\nu^\wr_c(M'_\idot) \cong
\nu_c(M'_\idot)$.
\endproof

\subsection{Filtration by support.}

We will now give some corollaries of
Proposition~\ref{nothing.to.do}.

First, the following simple observation. Suppose we are not
interested in the functor $\wPhi^c$ of \eqref{wphi.eq}, but only in
the fixed point functor $\Phi^c$ of Definition~\ref{phi.defn} ---
that is, we are prepared to forget the $\Aut(c)$-action on
$\wPhi^c$. Then there is a more convenient way to compute
$\Phi^c:\DS(\C,\Ab) \to \D(\Ab)$. Namely, denote by $\C_c$ the
category of objects $c' \in \C$ equipped with a map $c' \to c$. We
have a projection $s_c:\C_c \to \C$ which forgets the map. Composing
$s_c$ with the embedding $q$ gives a functor 
$$
q_c:\C_c \to \Q(\C\wr\Gamma)
$$ 
and a corresponding restriction functor $q_c^*:\DS(\C,\Ab) \to
\D(\C_c,\Ab)$. Let $\ol{T}^c \in \Fun(\C_c,\Z\amod)$ be given by
\begin{equation}\label{oT}
\ol{T}^c(c') =
\begin{cases}
\Z, &\quad c' \cong c,\\
0, &\quad\text{otherwise},
\end{cases}
\end{equation}
and let $\phi^c:\D(\C_c,\Ab) \to \D(\Ab)$ be the functor
left-adjoint to the embedding $\D(\Ab) \to \D(\C_c,\Ab)$, $M_\idot
\mapsto M_\idot \otimes \ol{T}$.

\begin{corr}\label{0.corr}
We have
$$
\Phi^c \cong \phi^c \circ q_c^*.
$$
In particular, $\Phi^cM_\idot=0$ if $q_c^*M_\idot = 0$.
\end{corr}

\proof{} The right-adjoint to the forgetful functor $\D(\br{c},\Ab)
\to \D(\Ab)$ sends $M_\idot \in \D(\Ab)$ to $M_\idot \otimes
\Z[\Aut(c)]$, where $\Z[\Aut(c)]$ is the regular representation of
the group $\Aut(c)$. Then by Proposition~\ref{nothing.to.do}, it
suffices to prove that $R^\hdot s_c(\ol{T}) \cong
\nu_c(\Z[\Aut(c)])$. This is clear: $s^c$ is obviously a discrete
fibration, so that $R^\hdot s^c_*$ can be computed fiberwise,
$\ol{T}$ is non-trivial on the fiber over $c' \in \C$ if and only
if $c'$ lies in $\br{c}$, and for every such $c'$, the fiber is a
torsor over $\Aut(c)$.
\endproof

Now, by Lemma~\ref{lattice}, the set $[\C]$ of isomorphism classes
of objects of the lattice-like category $\C$ has a natural partial
order, given by $[c] \geq [c']$ if and only if there exists a map $c
\to c'$. 

\begin{defn}
A subset $U \subset [\C]$ is {\em closed} if it is closed with
respect to the standard $T_0$-topology associated to the partial
order -- that is, for any $[c],[c'] \in [\C]$ with $[c] \geq [c']$,
$[c] \in U$ implies $[c'] \in U$.
\end{defn}

\begin{defn}\label{supp.defn}
For any subset $U \subset [\C]$, an object $E_\idot \in \DS(\C,\Ab)$
is {\em supported in $U$} if $E_\idot(c)=0$ for any $c \in \C$ whose
isomorphism class is outside $U$. An object $E_\idot \in \DS(\C\Ab)$
{\em has finite support} if it is supported in a finite closed
subset $U \subset [\C]$ The full subcategory in $\DS(\C,\Ab)$
spanned by objects supported in a subset $U \in [\C]$ is denoted by
$$
\DS_U(\C,\Ab) \subset \DS(\C,\Ab),
$$
and the full subcategory in $\DS(\C,\Ab)$ spanned by objects with
finite support is denoted by $\DS_{fs}(\C,\Ab) \subset
\DS(\C,\Ab)$.
\end{defn}

For any two subsets $U' \subset U \subset [\C]$, we obviously have a
full embedding $\DS_{U'}(\C,\Ab) \subset \DS_U(\C,\Ab)$.

\begin{corr}\label{1.corr}
For any closed subset $U \subset [\C]$, the fixed points functor
$$
\wPhi^c:\DS_U(\C,\Ab) \to \D(\br{c},\Ab)
$$
is trivial unless $[c] \in U$. If $[c] \in U$ is maximal, so that
$U' = U \setminus \{[c]\}$ is also closed, then $\wPhi^c \cong
\wPsi^c$ on $\DS_U(\C,\Ab)$, where $\wPsi^c$ is as in
\eqref{wpsi.eq}, the subcategory $\DS_{U'}(\C,\Ab) \subset
\DS_U(\C,\Ab)$ is right-admissible, and $\wPhi^c \cong \wPsi^c$
factors through an equivalence
\begin{equation}\label{quo.eq.bis}
\DS_U(\C,\Ab)/\DS_{U'}(\C,\Ab) \cong \D(\br{c},\Ab).
\end{equation}
\end{corr}

\proof{} Consider an object $E_\idot \in \DS_U(\C,\Ab)$. By
definition, for any $c \in \C$ whose isomorphism class is not in
$U$, we have $E_\idot(c) = 0$. Since $U$ is closed, we also have
$E_\idot(c')=0$ for any $c'$ with $[c'] \geq [c]$, so that $q_c^*E =
0$. By Corollary~\ref{0.corr}, this proves that $\Phi^c(E_\idot)=0$;
since the forgetful functor $\D(\br{c},\Ab) \to \D(\Ab)$ is
obviously conservative, we also have $\wPhi^c(E_\idot)$. This proves
the first claim.

For an object $c \in \C$ whose isomorphism class is maximal in $U$,
we have $q_c^*(E_\idot)(c') = 0$ unless $c' \in \C_c$ is isomorphic
to $c$, so that $q_c^*(E_\idot) \cong \Psi^c(E_\idot) \otimes
\ol{T}^c$. Since the correspondence $E_\idot \mapsto E_\idot \otimes
\ol{T}^c$ is obviously a full embedding from $\D(\Ab)$ to
$\D(\C_c,\Ab)$, we conclude that
$$
\phi^c(q_c^*E_\idot) \cong \Psi^c(E_\idot),
$$
so that by Corollary~\ref{0.corr}, the natural map $\Psi^c(E_\idot)
\to \Phi^c(E_\idot)$ is an isomorphism. Since the forgetful functor
$\D(\br{c},\Ab) \to \D(\Ab)$ is conservative, this proves the second
claim.

To prove the rest of the claims, note that the inflation functor
$\Infl^c:\D(\br{c},\Ab) \to \DS(\C,\Ab)$ obviously sends
$\D(\br{c},\Ab)$ into the full subcategory $\DS_U(\C,\Ab)$. By
definition, $\Infl^c$ is right-adjoint to $\Phi^c \cong
\Psi^c$. Moreover, we obviously have
$$
\Psi^c \circ \Infl^c \cong \Id.
$$
Therefore $\Infl^c:\D(\br{c},\Ab) \to \DS_U(\C,\Ab)$ is a fully
faithful embedding with left-admissible image. But by definition,
for any $E \in \DS_U(\C,\Ab)$, $\Psi^c(E) = 0$ if and only if $E$
lies in $DS_{U'}(\C,\Ab)$. Thus the orthogonal
${}^\perp\Infl^c(\D(\br{c},\Ab))$ is exactly $\DS_{U'}(\C,\Ab)$, and
the category $\DS_U(\C,\Ab)$ has a semiorthogonal decomposition
$\langle \Infl^c(\D(\br{c},\Ab)),\DS_{U'}(\C,\Ab) \rangle$, which
finishes the proof.
\endproof

\begin{remark}
In the particular case $\C = O_G$, $c = [G/H]$ for a cofinite
subgroup $H \subset G$, and $U = \{[c'] \in [\C]| [c'] \leq [c]\}$,
\eqref{quo.eq.bis} is exactly \eqref{quo.eq} of the Introduction.
\end{remark}

\begin{corr}\label{2.corr}
For any two finite closed subsets $U' \subset U \subset [\C]$, the
category $\DS_{U'}(\C,\Ab)$ is an admissible full subcategory in
$\DS_U(\C,\Ab)$, and the left orthogonal ${}^\perp\DS_{U'}(\C,\Ab)$
consists of those objects $E \in \DS_U(\C,\Ab)$ which are supported
in $U \setminus U'$.
\end{corr}

\proof{} The fact that $\DS_{U'}(\C,\Ab) \subset \DS_U(\C,\Ab)$ is
right-admissible immediately follows by induction from
Corollary~\ref{1.corr}. Moreover, the semiorthogonality statement of
Corollary~\ref{1.corr} implies that $\DS_{U'}(\C,\Ab)$ is generated
by objects of the form $\Infl^{c'}(M)$, $c' \in U'$, $M \in
\D(\br{c'},\Ab)$. Thus $E \in \DS_U(\C,\Ab)$ lies in the orthogonal
${}^\perp\DS_{U'}(\C,\Ab)$ if and only if it is orthogonal to
all such $\Infl^{c'}(M)$. By adjunction this is equivalent to
$\wPhi^{c'}(E) = 0$, $c' \in U'$, which by definition means exactly
that $E$ is supported in $U \setminus U'$.

Finally, to prove that $\DS_{U'}(\C,\Ab) \subset \DS_U(\C,\Ab)$ is
left-admissible, it suffices by induction to consider the case $U' =
U \setminus \{[c]\}$ for a maximal $c \in U$. By Lemma~\ref{orth}, it
suffices to prove that the projection
$$
{}^\perp\DS_{U'}(\C,\Ab) \to \DS_U(\C,\Ab)/\DS_{U'}(\C,\Ab) \cong
\D(\br{c},\Ab)
$$
is essentially surjective. Define a functor $I:\Fun(\Aut(c),\Ab) \to
\Fun(\C\wr\Gamma,\Ab)$ by
$$
I(M) = (T^c \otimes M)_{\Aut(c)},
$$
and let
$$
L^\hdot I:\D(\Aut(c),\Ab) \cong \D(\br{c},\Ab) \to \DS_U(\C,\Ab)
\subset \D(\Fun(\C\wr\Gamma,\Ab))
$$
be its left-derived functor (this is dual to the definition of the
inflation functor $\Infl^c$ in Subsection~\ref{der.func}). We
obviously have $\wPhi^c \circ L^\hdot I \cong \wPsi^c \circ L^\hdot
I \cong \Id$, so that it suffices to prove that $L^\hdot I$ sends
$\D(\Aut(c),\Ab)$ into the orthogonal ${}^\perp\DS_{U'}(\C,\Ab)$. As
we have already proved, this is equivalent to proving that
$$
\Phi^{c'} \circ L^\hdot I = 0
$$
for any $c' \in U'$. Thus by Corollary~\ref{0.corr}, we have to
prove that
\begin{equation}\label{t.p}
\phi^{c'} \circ q^*_{c'} \circ L^\hdot I = 0
\end{equation}
for any $c' \in U'$. By definition, we have
$$
q_{c'}^*(T^c \otimes M) \cong s_{c'}^*j^c_!M
$$
for any $M \in \Fun(\br{c},\Ab)$, where $j^c:\br{c} \to \C$ is the
natural embedding. Passing to the derived functors, we obtain
$$
q^*_{c'} \circ L^\hdot I \cong s_{c'}^* \circ L^\hdot j^c_!,
$$
and by the definition of the functor $\phi^{c'}$, \eqref{t.p} is
equivalent to
$$
j^c_!R^\hdot s_{c'*}\ol{T}^{c'} = 0,
$$
where $\ol{T}^{c'} \in \Fun(\C_{c'},\Z\amod)$ is as in
\eqref{oT}. Since the projection $s_{c'}:\C_{c'} \to \C$ is a
fibration, this immediately follows by base change.
\endproof

\subsection{DG models.} Corollary~\ref{1.corr} and Corollary~\ref{2.corr}
show that the category $\DS_{fs}(\C,\Ab)$ admits a filtration by
admissible full triangulated subcategories $\DS_U(\C,\Ab)$ indexed
by finite closed subsets $U \subset [\C]$, and the associated graded
quotients
$$
\DS_c(\C,\Ab) = \DS_U(\C,\Ab)/\DS_{U \setminus \{[c]\}}(\C,\Ab),
\quad c \in U \quad \text{maximal}
$$
of this filtration are identified with the derived categories
$\D(\Aut(c),\Ab)$. In the remainder of this section, we will obtain
a full description of the category $\DS_{fs}(\C,\Ab)$ in terms of
the functors $\wPhi^c$.  Among other things, this will later allow
us to describe the gluing data between these associated graded
quotients.

\medskip

The general strategy is rather straightforward -- we apply a version
of the standard Tanakian formalism. Namely, we treat the collection
$\wPhi^c$ as a ``fiber functor'' for the category
$\DS_{fs}(\C,\Ab)$, and
\begin{enumerate}
\item lift the functors $\wPhi^c$ to DG functors $\wPhi^c_\idot$, and
\item present a natural $\C$-graded $A_\infty$-coalgebra $T^\C_\idot(-)$
  which acts on $\wPhi^c_\idot(-)$ (in particular, this will include
  the $\Aut(c)$-action).
\end{enumerate}
The functor $\Phi^c:\DS_{fs}(\C,\Ab) \to \D(\Ab)$ is by definition
left-adjoint to the functor $\D(\Ab) \to \DS_{fs}(\C,\Ab)$ given by
$M \mapsto M \otimes T^c$; the ``natural $A_\infty$-coalgebra'' in
\thetag{ii} is then an $A_\infty$ model for
$$
\Phi^c(T^{c'}), \qquad c,c' \in \C,
$$
with the comultiplication given by adjunction. The same adjunction
gives a comparison functor from $\DS_{fs}(\C,\Ab)$ to the derived
category of $A_\infty$-comodules over $\T^\C_\idot$, and we check that this
comparison functor is an equivalence (this is
Theorem~\ref{galois.comp}).

We note that here it is crucial to work with $A_\infty$-coalgebras
rather than $A_\infty$-algebras. In fact, the complex
$T^\C_\idot(f)$ is acyclic for any non-invertible map $f:c \to
c'$. However, as we have noted in Subsection~\ref{dgcat.subs},
quasiisomorphic $A_\infty$-coalgebras may have different categories
of $A_\infty$-comodules, and in particular, acyclic complexes
$\T^\C_\idot(f)$ cannot be replaced with $0$.

\subsubsection{The complexes.}
Fix an object $c \in \C$. For any $n \geq 1$, let $\C'_n(c)$ be the
groupoid of diagrams $c_1 \to \dots \to c_n \to c$ in $\C$ and
invertible maps between them, and let $\C_n(c) \subset \C_n'(c)$ be
the full subcategory spanned by such diagrams that
\begin{itemize}
\item the map $c_n \to c$ is not an isomorphism.
\end{itemize}
As in Subsection~\ref{bar.gene}, let $\sigma_n:\C'_n \to \C$ be the
functor which sends a diagram to $c_1 \in \C$. For any $E \in
\Fun(S(\C\wr\Gamma),\Ab)$, let
\begin{equation}\label{phi.eq}
\Phi^c_{n,\idot}(E) = C_\idot(\C_n(c),\sigma_n^*q^*E).
\end{equation}
For $n = 0$, set $\Phi^c_{0,0}(E) = E(c)$ and $\Phi^c_{0,i}(E) = 0$
for $i \neq 1$. As in Subsection~\ref{bar.gene}, define a second
differential $d:\Phi^c_{n+1,\idot}(E) \to \Phi^c_{n,\idot}(E)$ as
the alternating sum of maps $d^i$, $1 \leq i \leq n$, where $d^i$
removes the object $c_i$ from the diagram, and let $\Phi^c_\idot(E)$
be the total complex of the resulting bicomplex.

Any map $f:c \to c$ acts on the diagrams $c_1 \to \dots \to c_n \to
c$ by composition on the right-hand side; this turns
$\Phi^c_\idot(E)$ into a complex $\wPhi^c_\idot(E)$ of
representations of the group $\Aut(c)$.

\begin{lemma}\label{phi.idot}
The complex $\wPhi^c_\idot(E)$ computes $\wPhi^c(E)$.
\end{lemma}

\proof{} By Proposition~\ref{nothing.to.do}, it suffices to
construct a functorial quasiisomorphism
\begin{equation}\label{t.p.1}
\RHom^\hdot_{\br{c}}(\wPhi^c_\idot(E),M) \cong
\RHom^\hdot(q^*E,\nu_c(M))
\end{equation}
for any $M \in \Fun(\br{c},\Ab)$. Computing the left-hand side of
\eqref{t.p.1} by \eqref{rhom.bar.eq.1}, we get
\begin{multline}\label{lh.I}
\bigoplus_{0 \leq i} \RHom^{\hdot-i}_{\br{c}}(\wPhi^c_i(E),M)\cong\\
\begin{aligned}
&\cong\bigoplus_{0 \leq n}\bigoplus_{0 \leq i < n}C^\hdot(\br{c}_{n-i},
\RHom^{\hdot+1-n}(C_\idot(\C_i(c), \sigma_i^*q^*E),M))\\
&\cong\bigoplus_{0 \leq n}\bigoplus_{0 \leq i < n}C^\hdot(\C_i(c) \times
\br{c}_i,\RHom^{\hdot+1-m}(\sigma_n^*q^*E,M)),
\end{aligned}
\end{multline}
where $\br{c}_i$ is as in Subsection~\ref{bar.gene}. For every
diagram $c_1 \to \dots \to c_n \to c$ in $\C'_n(c)$, there exists a unique
$i$, $0 \leq i < n$ such that $c_i \to c_{i+1}$ is not an
isomorphism, but  $c_j \cong c$ for $j \geq i+1$. Therefore we have
a functorial decomposition
$$
\C'_n(c) = \coprod_{0 \leq i < n}\C_i(c) \times \br{c}_{n-i},
$$
and \eqref{lh.I} can be rewritten as
$$
\bigoplus_{n \geq
  0}C^\hdot(\C'_n(c),\RHom^{\hdot+1-n}(\sigma_n^*q_n^*E,M)).
$$
To compute the right-hand side of \eqref{t.p.1}, we can also use
\eqref{rhom.bar.eq.1}; this gives the same formula, except that
$\C'_n(c)$ is replaced with the groupoid $\C_n$ of
Subsection~\ref{bar.gene}, and $M$ is replaced by
$\tau_n^*\nu_c(M)$. It remains to notice that by the definition of
the functor $\nu_c$, the terms corresponding to $\C_n \setminus
\C'_n(c)$ are all equal to $0$.
\endproof

Now consider two objects $c,c' \in \C$, and apply Lemma~\ref{phi.idot}
to the object $T^{c'} \in \Fun(S(\C\wr\Gamma),\Z\amod)$. Recall that
by definition, for any diagram $\alpha=[c_1 \to \dots \to c_n \to c]
\in \C_n(c)$ the group $T^{c'}(\alpha) = T^{c'}(c_1)$ is the free
abelian group generated by the space of maps from $c'$ to
$c_1$. These maps can be incorporated into the diagram. For any map
$f:c' \to c$ and any $n \geq 1$, denote by $\C_n(f)$ the groupoid of all
diagrams $c' \to c_1 \to \dots \to c_n \to c$ such that the
composition of all the maps is equal to $f$, and the last map $c_n
\to c$ is not an isomorphism. Let
\begin{equation}\label{t.n.eq}
\T^\C_{n,\idot}(f) = C_\idot(\C_n(f),\Z),
\end{equation}
and extend it to $n=0$ by setting $\T^\C_{0,i}(f) = \Z$ for $i = 0$
and $0$ otherwise. For any $1 \leq i \leq n$, forgetting the object
$c_i$ in the diagram gives functors $\C_n(f) \to \C_{n-1}(f)$ and
the corresponding maps $d_i:\T^\C_n(f) \to \T^\C_{n-1}(f)$; define a
second differential $d:\C_n(f) \to \C_{n-1}(f)$ as their alternating
sum, and let $\T^\C_\idot(f)$ be the total complex of the
corresponding bicomplex. By definition, we have an isomorphism of
complexes
$$
\Phi^c_\idot(T^{c'}) \cong \bigoplus_{f \in \C(c',c)}\T^\C_\idot(f).
$$

\subsubsection{Coaction.}

Consider a diagram $\alpha = [c_1 \to \dots \to c_n \to c] \in
\C_n(c)$. Since $\C\wr\Gamma$ has fibered products, we can form the
fibered product diagram $c_1 \times_c c' \to \dots \to c_n \times_c
c' \to c'$. By definition, for any $1 \leq i \leq n$, $c_i \times_c
c'$ is the formal union
$$
c_i \times_c c' = \coprod_{s \in S_i} c_{i,s}
$$
of some objects $c_{i,s} \in \C$, numbered by elements in a finite
set $S_i$. Denote $S(\alpha) = S_1$. Picking up an element $s \in
S(\alpha)$ uniquely determines an element $s \in S_i$ for any $i$,
and we obtain a diagram $c_{1,s} \to \dots \to c_{n,s} \to
c'$. Moreover, we have a commutative diagram
$$
\begin{CD}
c_1 @>>> \dots @>>> c_n @>>> c\\
@AA{f_{1,s}}A @. @AA{f_{n,s}}A @AA{f}A\\
c_{1,s} @>>> \dots @>>> c_{n,s} @>>> c'.
\end{CD}
$$
By assumption, the map $c_i \to c$ is not an isomorphism for any
$i$, but the map $c_{i,c} \to c'$ might well be. So, let $i$ be the
smallest integer such that all the maps $c_{j,s} \to c$, $j > i$ are
isomorphisms. Then $c_{1,s} \to \dots \to c_{i,s} \to c'$ is a
well-defined object of the groupoid $\C_i(c)$, which we denote by
$\alpha^{(1)}_s$. Moreover, composing the map $f_{i,s}:c_{i,s} \to
c_i$ with the isomorphism $c' \to c_{i,s}$ inverse to the map
$c_{i,s} \to c'$, we obtain a diagram
$$
\begin{CD}
c' \cong c_{i,s} @>{f_i,s}>> c_i @>>> \dots @>>> c_n @>>> c
\end{CD}
$$
which gives an object in the groupoid $\C_{n-i}(f)$; we denote it by
$\alpha^{(2)}_s$. 

\medskip

These constructions are functorial, so that sending a diagram
$\alpha$ to the formal union of the products $\alpha_s^{(1)} \times
\alpha_s^{(2)}$, $s \in S_1$ gives a functor
$$
\beta_n:\C_n(c) \to \left(\coprod_{0 \leq i \leq n}\C_i(c') \times
\C_{n-i}(f)\right) \wr \Gamma.
$$
For any $E \in \Fun(S(\C\wr\Gamma),\Ab)$, we have a natural map
$$
\sigma_n^*q^*E \to \beta_n^*\sigma_n^{'*}E,
$$
where $\sigma_n^{'*}$ is obtained by extending the natural
functor
$$
\begin{CD}
\C_i(c') \times \C_{n-1}(f) @>>> \C_n(c') @>{\sigma_n}>> \C
\end{CD}
$$
to wreath products. If $E$ is additive in the sense of
Definition~\ref{add.defn.ds}, then $\sigma_n^{'*}E$ is also
additive, so that we have a canonical decomposition
$$
\sigma_n^{'*}E\left(\copr_s \alpha^{(1)}_s \times \alpha^{(2)}_s\right)
\cong \bigoplus_s E(\alpha^{(1)}_s).
$$
We can then compose the map $C_\idot(\beta_\idot)$ induced
by the functors $\beta_\idot$ with the projections onto the terms of
this decomposition and obtain a canonical map of complexes
\begin{equation}\label{m.f}
b_f:\Phi^c_\idot(E) \to \Phi^{c'}_\idot(E) \otimes \T^\C_\idot(f).
\end{equation}

\subsubsection{Comultiplication and higher operations.}

Now assume given a third object $c'' \in \C$ and a map $g:c'' \to
c'$, and apply \eqref{m.f} to $E = T^{c''}$. This gives a
comultiplication map
$$
b_{f,g}:\T^\C_\idot({f \circ g}) \to T_\idot(f) \otimes \T^\C_\idot(g)
$$
by the following procedure. Consider again a diagram $\alpha = [c''
\to c_1 \to \dots \to c_n \to c] \in \C_n(f \circ g)$, and form
the pullback diagram $c_1 \times_c c' \to \dots \to c_n \times_c c'
\to c'$. In this case, we also have a natural map $c'' \to c_1
\times_c c'$; this distinguishes a component $c_1'$ in the
decomposition $c_1 \times_c c' = \coprod c_{1,s}$, hence also a
component $c_i'$ in the decomposition $c_i \times_c c' = \coprod
c_{i,s}$ for any $1 \leq i \leq n$, and we have a commutative
diagram
$$
\begin{CD}
c'' @>>> c_1 @>>> \dots @>>> c_n @>>> c\\
@| @AA{f_1}A @. @AA{f_n}A @AA{f}A\\
c'' @>>> c_1' @>>> \dots @>>> c_n' @>>> c'.
\end{CD}
$$
We again take the smallest integer $i$ such that $c_{i+1}' \to c'$ is an
isomorphism, and obtain a functor
$$
\beta'_n:\C_n(f \circ g) \to \coprod_{0 \leq i \leq n}\C_i(g) \times
\C_{n-i}(f).
$$
This functor induces our map $b_{f,g}$.

Moreover, since our functors $\beta_n$, $\beta'_n$ are defined by
pullbacks, they are associative up to a canonical isomorphism. This
means that for any composable $l$-tuple of maps $f_1,\dots,f_l$ in
the category $\C_l$, we have canonical functors
$$
\beta'_{n,l}:\C_n(f_1 \circ \dots \circ f_l) \times I_l  \to
  \coprod_{n_1 + \dots + n_l = n}\C_{n_1}(f_1) \times \dots \times
  \C_{n_l}(f_l),
$$
where $I_l$ are the terms of the asymmetric operad of categories
considered in Subsection~\ref{2cat.subs}, and these functors are
compatible with the operad structure on $I_l$. Taking the bar
complexes and using the operad map $\Ass_\infty \to
C_\idot(I_\idot,\Z)$ of Subsection~\ref{2cat.subs}, we turn
$\T^\C_\idot(f)$ into a $\C$-graded $A_\infty$-coalgebra in the sense
of Definition~\ref{DG.coa}. The analogous functors $\beta_{n,l}$
turn the collection $\Phi^c_\idot(E)$, $c \in \C$, into an
$\Ab$-valued $A_\infty$-comodule over $\T^\C_\idot(-)$.

\subsection{The comparison theorem.} We can now prove the
comparison theorem expressing $\DS_{fs}(\C,\Ab)$ in terms of the
coalgebra $\T^\C_\idot(-)$.

\begin{defn}\label{T.defn}
The category $\DT(\C,\Ab) = \D(\T^\C_\idot,\Ab)$ is the derived
category of $\Ab$-valued $A_\infty$-comodules over the $\C$-graded
$A_\infty$-coalgebra $\T^\C_\idot(-)$. An object $E_\idot \in
\DT(\C,\Ab)$ is {\em finitely supported} if $E_\idot(c) = 0$ except
for $c$ with isomorphism class in a finite closed subset $U \subset
[\C]$; the full subcategory spanned by finitely supported objects is
denoted by $\DT_{fs}(\C,\Ab) \subset \DT(\C,\Ab)$.
\end{defn}

In the last Subsection, we have constructed a functor
$\Phi_\idot:\Fun(S(\C\wr\Gamma),\Ab) \to \DT(\C,\Ab)$. Taking the
total complex of a double complex, we extend it to a functor
$$
\Phi_\idot:\DS(\C,\Ab) \to \DT(\C,\Ab).
$$
(There is the usual ambiguity in taking the total complex of a
possibly infinite bicomplex -- we can either take the sum of the
terms on a diagonal, or the product of these terms. Here we take the
sum.)

\begin{theorem}\label{galois.comp}
The functor $\Phi_\idot$ induces an equivalence
$$
\Phi_\idot:\DS_{fs}(\C,\Ab) \to \DT_{fs}(\C,\Ab).
$$
\end{theorem}

\proof{} For any subset $U \subset [\C]$, let $\C_U \subset \C$ be
the full subcategory spanned by objects with isomorphism classes in
$\C$, with the embedding functor $\iota_U:\C_U \to \C$, and let
$\T^{\C_U}_\idot(-) = \iota_U^*\T^\C(-)$ be the $\C_U$-graded
$A_\infty$-coalgebra made up of the complexes $\T^\C_\idot(f)$ with
$f$ in $\C_U$. Let $\DT_U(\C,\Ab)$ be the triangulated category of
$A_\infty$-comodules over $\T^{\C_U}_\idot$. We have the restriction
functor $\iota_U^*:\DT(\C,\Ab) \to \DT_U(\C,\Ab)$, and for any two
subsets $U' \subset U$, we have a restriction functor
$\iota_{U,U'}^*:\DT_U(\C,\Ab) \to \DT_{U'}(\C,\Ab)$. If $U$ is
closed, and either $U' \subset U$ or the complement $U \setminus U'
\subset$ is also closed, then we have an obvious inverse functor
$\DT_{U'}(\C,\Ab) \to \DT_U(\C,\Ab)$ given by extension by $0$ ---
explicitly, $M \in \DT_{U'}(\C,\Ab)$ is sent to $M' \in
\DT_U(\C,\Ab)$ such that
$$
M'(c) =
\begin{cases}
M(c), &\quad [c] \in U',\\
0, &\quad \text{otherwise},
\end{cases}
$$
with the obvious structure maps.  If it was $U'$ which was closed,
then this extension functor is right-adjoint to the restriction
functor $i_{U,U'}^*$ (this is obvious for the homotopy categories
$\Ho(\T^{\C_U}_\idot,\Ab)$, $\Ho(\T^{\C_U'}_\idot,\Ab)$, and since
both restriction and extension preserve acyclic comodules, the
adjunction descends to the derived categories).

\begin{lemma}\label{orth.coa}
For any two finite closed subsets $U' \subset U \subset [\C]$, the
category $\DT_U(\C,\Ab)$ admits a semiorthogonal decomposition
$$
\langle \DT_{U \setminus U'}(\C,\Ab),\DT_{U'}(\C,\Ab),
$$
where the embeddings $\DT_{U \setminus U'}(\C,\Ab) \subset
\DT_U(\C,\Ab)$ and $\DT_{U'}(\C,\Ab) \subset \DT_U(\C,\Ab)$ are
given by extension by $0$.
\end{lemma}

\proof{} By adjunction, $\DT_{U'}(\C,\Ab) \subset \DT_U(\C,\Ab)$ is
left-admissible, with the projection given by restriction. On the
other hand, $\DT_{U \setminus U'}(\C,\Ab)$ by definition is exactly
the kernel of the restriction functor.
\endproof

Lemma~\ref{orth.coa} implies by induction that $\DT_{fs}(\C,\Ab) =
\cup_U \DT_U(\C,\Ab)$, the union over all finite closed $U \subset
[\C]$, where $\DT_U(\C,\Ab)$ is identified with the full subcategory
in $\DT(\C,\Ab)$ spanned by objects $E_\idot$ such that $E_\idot(c)
= 0$ unless $[c] \in U$. The comparison functor $\Phi_\idot$
obviously sends $\DS_U(\C,\Ab)$ into $\DT_U(\C,\Ab)$, and it
suffices to prove that it is an equivalence for any finite closed $U
\subset [\C]$. Fix such a subset $U$, take an object $c \in \C$
whose isomorphism class $[c]$ is a maximal element in $U$, and let
$U' = U \setminus \{[c]\}$. By induction on the cardinality of $U$,
we may assume that $\Phi_\idot:\DS_{U'}(\C,\Ab) \to \DT_U(\C,\Ab)$
is an equivalence. Moreover, for any invertible map $f$, the groupoids
$\C_n(f)$, $n \geq 1$ are obviously empty, so that $\T^\C_\idot(f)$
is $\Z$ placed in degree $0$. Therefore the category
$\DT_{\{[c]\}}(\C,\Ab)$ is equivalent to $\D(\br{c},\Ab)$, and the
functor
$$
\Phi_\idot:{}^\perp\DS_{U'}(\C,\Ab) \cong \D(\br{c},\Ab) \to
\DT_{\{[c]\}}(\C,\Ab)
$$
is also an equivalence. Then by Corollary~\ref{2.corr},
$\DS_{U'}(\C,\Ab) \subset \DS_U(\C,\Ab)$ is admissible, so that by
Lemma~\ref{orth.equi}, it suffices to prove that $\Phi_\idot$ sends
$$
\DS_{U'}(\C,\Ab)^\perp = \Infl^c(\D(\br{c},\Ab)) \subset
\DS_U(\C,\Ab)
$$
into $\DT_{U'}(\C,\Ab)^\perp \subset \DT_U(\C,\Ab)$.

Indeed, denote $\Ho_{\{c\}} \cong \Ho(\T^{\C_{\{c\}}}_\idot,\Ab)$,
$\Ho_U = \Ho(\T^{\C_U}_\idot,\Ab)$. Then the restriction functor
$\Ho_U \to \Ho_{\{c\}}$ has an obvious right-adjoint
  $\I:\Ho_{\{[c]\}} \to \Ho_U$ which sends a complex $M_\idot$ of
  $\Aut(c)$-modules into an $A_\infty$-comodule $\I(M_\idot)$ such
  that
\begin{equation}\label{I.gal}
\I(M_\idot)(c') = \left(\bigoplus_{f \in \C(c,c')}T_\idot(f) \otimes
M_\idot\right)^{\Aut(c)}.
\end{equation}
If the complex $M_\idot$ is $h$-injective, then $\I(M_\idot)$ is
also $h$-injective by adjunction, so that its image in the quotient
category $\DT_U(\C,\Ab)$ lies in $\DT_{U'}(\C,\Ab)^\perp \subset
\DT_U(\C,\Ab)$. However, since $\Ho_{\{c\}}$ is equivalent to the
category of unbounded complexes of functors in $\Fun(\br{c},\Ab)$,
every object has an $h$-injective replacement.  Thus we can take any
object $M' = \Infl^c(M) \in \DS_{U'}(\C,\Ab)^\perp$, and choose an
$h$-injective complex $M_\idot$ of $\Aut(c)$-comodules representing
$M = \wPhi^c(M') \in \D(\br{c},\Ab) \cong \D(\Aut(c),\Ab)$. Then by
definition, $M'=\Infl^c(M)$ is represented by the complex
$$
\Infl^c(M_\idot) = (M_\idot \otimes T^c)^{\Aut(c)},
$$
and we see that $\Phi_\idot(M_\idot) \in \DT_U(\C,\Ab)$ is precisely
isomorphic to $\I(M_\idot)$.
\endproof

\begin{remark}
  As the proof of Theorem~\ref{galois.comp} shows, every finitely
  supported object in the category $\Ho(\T_\idot,\Ab)$ does have an
  $h$-injective replacement. As we have mentioned in
  Subsubsection~\ref{coa.subsu}, this property seems rather special.
\end{remark}

\subsection{Induction and products.}\label{DG.indu.subs}

Assume now that the category $\C$ has a terminal object, as in
Subsection~\ref{der.prod}. Then $\DS(\C,\Ab)$ is a tensor category,
and by Proposition~\ref{phi.tensor} all the fixed points functors
$\Phi^c$ are tensor functors, so that the subcategory
$\DS_{fs}(\C,\Ab) \subset \DS(\C,\Ab)$ is closed under tensor
products. To describe the tensor product structure on
$\DS_{fs}(\C,\Ab)$ in terms of the equivalent category
$\DT_{fs}(\C,\Ab)$, one would have to introduce a product on the
$\C$-graded $A_\infty$-coalgebra $\T^\C_\idot$ -- more precisely, an
associative product on the complex $\T^\C_\idot(f)$ for any map $f$
in $\C$. In keeping with the Tanakian formalism, $\T^\C_\idot$ then
becomes a ``$\C$-graded Hopf algebra'', with an induced tensor
product on the category $\DT(\C,\Ab)$.

Unfortunately, it seems that the particular $\C$-graded
$A_\infty$-coalgebra $\T^\C_\idot$ that we have used in
Theorem~\ref{galois.comp} does not have a natural Hopf algebra
structure. Possibly this could be healed by taking some
quasiisomorphic $A_\infty$-coalgebra, or using a Segal-like notion
of a ``special $A_\infty$ $\Gamma_+$-coalgebra'' instead of a
commutative $A_\infty$ Hopf algebra, but I haven't pursued
this. Nevertheless, we have a simpler result which we will need for
applications in Section~\ref{spectra}.

Consider the trivial $\C$-graded coalgebra $\Z^\C$ given by
$\Z^\C(f) = \Z$ for any map $f$ in $\C$. Then since $\C$ is assumed
to be $\Hom$-finite, the derived category $\D(\Z^\C,\Ab)$ of
$\Ab$-valued $A_\infty$-comodules over $\Z^\C$ is equivalent to the
derived category $\D(\C^{opp},\Ab)$. Moreover, since $\T^0_\C(f) =
\Z = \Z^\C(f)$ for any $f$, we have a natural map
$$
\lambda:\Z^\C = F^0\T^\C_\idot \to \T^\C_\idot,
$$
where $F^0$ stands for the $0$-th term of the stupid
filtration. This map induces a corestriction functor
$\lambda^*:\D(\C^{opp},\Ab) \cong \D(\Z^\C,\Ab) \to \DT(\C,\Ab)$.
Explicitly, for any $\Z^\C$-comodule $E_\idot$, we have
$\lambda^*E_\idot(c) \cong E_\idot$, $c \in \C$, and for any map
$f:c' \to c$ in $\C$, the comodule structure map of
Definition~\ref{DG.coa} is the natural map
$$
\begin{CD}
E_\idot(c) @>{E_\idot(f)}>> E_\idot(c') \otimes \Z = E_\idot(c')
\otimes \Z^\C(f) @>{\lambda(f)}>> E_\idot(c') \otimes \T^\C_\idot(f),
\end{CD}
$$
with all the higher maps being equal to $0$. Moreover, for any
complex $E_\idot$ of pointwise-flat functors in $\Fun(\C^{opp},\Ab)$
and any $A_\infty$-comodule $M_\idot$ over $\T^\C_\idot$, we have an
obvious $A_\infty$-comodule $E_\idot \otimes M_\idot$ over
$\T^\C_\idot$, with $(E_\idot \otimes M_\idot)(c) = E_\idot(c)
\otimes M_\idot(c)$, and with the comodule structure maps obtained
as the products of the comodule structure maps for $M_\idot$ and the
canonical maps
$$
E_\idot(c) \to E_\idot(c')
$$
for any composable $n$-tuple $f_1,\dots,f_n$ of maps in $\C$ with
the composition $f_1 \circ \dots \circ f_n$ being a map from $c$ to
$c'$. This preserves quasiisomorphisms, thus descents to a
tensor product functor
\begin{equation}\label{T.C.prod}
\D(\C^{opp},\Ab) \times \DT(\C,\Ab) \to \DT(\C,\Ab).
\end{equation}

\begin{lemma}\label{prod.indu.como}
Under the comparison functor $\Phi_\idot:\DS(\C,\Ab) \to
\DT(\C,\Ab)$ of Theorem~\ref{galois.comp}, the induction functor
$q^{opp}_!:\D(\C^{opp},\Ab) \to \DS(\C,\Ab)$ of
Subsection~\ref{der.indu} goes to the corestriction functor
$\lambda^*$ -- that is, we have a natural functorial
quasiisomorphism
\begin{equation}\label{prod.indu.1}
\Phi_\idot(q^{opp}_!(E_\idot)) \cong \lambda^*E_\idot
\end{equation}
for any $E_\idot \in \D(\C^{opp},\Ab)$. Moreover, for any objects
$E_\idot \in \D(\C^{opp},\Ab)$, $M_\idot \in \DS(\C,\Ab)$, we have a
natural isomorphism
\begin{equation}\label{prod.indu.2}
\Phi_\idot((q^{opp}_!E_\idot) \otimes M_\idot) \cong E_\idot \otimes
\Phi_\idot(M_\idot),
\end{equation}
where the product in the right-hand side is the product
\eqref{T.C.prod}.
\end{lemma}

\proof{} Both claims immediately follow from
Proposition~\ref{indu.prod}.
\endproof

\section{Tate homology description.}\label{tate.sec}

We will now use Theorem~\ref{galois.comp} to obtain an explicit and
useful description of the gluing data in the semiorthogonal
decomposition of Corollary~\ref{1.corr}.

\subsection{Generalized Tate cohomology.}

Assume given a finite group $G$, and a family $\{H_i\}$ of subgroups
$H_i \subset G$. For every $H_i$, we have the induction functor
$$
\Ind_G^{H_i}:\D^b(H_i,\Z) \to \D^b(G,\Z),
$$
where $\D^b(G,\Z)$, resp. $\D^b(H_i,\Z)$ are the bounded derived
categories of $\Z[G]$, resp. $\Z[H_i]$-modules; this functor is
adjoint both on the left and on the right to the natural restriction
functor $\D^b(G,\Z) \to \D^b(H_i,\Z)$. Let
$$
\D^b_{\{H_i\}}(G,\Z) \subset \D^b(G,\Z)
$$
be the full thick triangulated subcategory spanned by direct
summands of objects of the form $\Ind^{H_i}_G(V_\idot)$, $V_\idot
\in \D^b(H_i,\Z)$, $H_i \in \{H_i\}$, and denote by
$$
\D(G,\{H_i\},\Z) = \D^b(G,\Ab)/\D^b_{\{H_i\}}(G,\Z)
$$
the quotient category.

\begin{defn}\label{gen.tate.Z}
The {\em generalized Tate cohomology} $\vH^\hdot(G,\{H_i\},M_\idot)$
of the group $G$ with coefficients in an object $M_\idot \in
\D^b(G\,\Z)$ with respect to the family $\{H_i\}$ is given by
$$
\vH^\hdot(G,\{H_i\},M_\idot) =
\Ext^\hdot_{\D(G,\{H_i\},\Ab)}(\Z,M_\idot),
$$
where $\Z$ is the trivial representation of the group $G$, and
$\Ext^\hdot(-,-)$ is computed in the quotient category
$\D(G,\{H_i\},\Z)$.
\end{defn}

By the definition of the quotient category,
$\vH^\hdot(G,\{H_i\},M_\idot)$ is expressed as follows. Consider the
category $I(G,\{H_i\})$ of objects $V_i \in \D^b(G,\Z)$ equipped
with a map $V_\idot \to \Z$ whose cone lies in $\D^b_{\{H_i\}}(G,\Z)
\subset \D^b(G,\Z)$. Let $I=I(G,\{H_i\})^{opp}$ be the opposite
category, and let
\begin{equation}\label{i.small}
I^{fg} \subset I = I(G,\{H_i\})^{opp}
\end{equation}
be the subcategory spanned by those $V_\idot$ which can be
represented by complexes of finitely generated $\Z[G]$-modules. Then
we have
\begin{equation}\label{tate.eq.1}
\vH^\hdot(G,\{H_i\},M_\idot) = \lim_{\to}\Ext^\hdot(V_\idot,M_\idot),
\end{equation}
where the limit is taken over the category $I^{fg}$, and
$\Ext^\hdot(-,-)$ is computed in the category $\D(G,\Z)$. Sending
$V_\idot$ to the dual complex $V_\idot^*$ identifies $I^{fg}$ with
the category of objects $V_\idot^*$ in $\D^b(G,Z)$ which are
represented by finitely generated $\Z[G]$-modules and equipped with
a map $\Z \to V^*_\idot$ whose cone lies in
$\D^b_{\{H_i\}}(G,\Z)$. Then \eqref{tate.eq.1} can be rewritten as
$$
\vH^\hdot(G,\{H_i\},M_\idot) = \lim_{\to}H^\hdot(G,M_\idot
\otimes V_\idot),
$$
where again the limit is over $I^{fg}$. Since $I^{fg}$ is small and
filtered, the limit is well-defined, and gives a cohomological
functor.

We can use the same expression to define Tate cohomology with
coefficients. Namely, let $\D^b(G,\Ab)$ be the derived category of
representations of the group $G$ in the abelian category $\Ab$.

\begin{defn}\label{gen.tate}
The {\em generalized Tate cohomology} $\vH^\hdot(G,\{H_i\},M_\idot)$
of the group $G$ with coefficients in an object $M_\idot \in
\D^b(G\,\Ab)$ with respect to the family $\{H_i\}$ is given by
$$
\vH^\hdot(G,\{H_i\},M_\idot) = \lim_{\to}H^\hdot(G,M_\idot
\otimes V_\idot),
$$
where the limit is taken over the category $I^{fg} \subset I =
I(G,\{H_i\})$ of \eqref{i.small}
\end{defn}

Again, since $\Ab$ is by assumption a Grothendieck category, and
$I^{fg}$ is small and filtered, generalized Tate cohomology is a
well-defined cohomological functor from $\D(G,\Ab)$ to $\Ab$.

\begin{lemma}\label{no.tate.prod}
For any $M \in \D(G,\Ab)$ and any $W_\idot \in
\D_{\{H_i\}}(G,\Z\amod)$ which can be represented by a bounded
complex of finitely generated $\Z[G]$-modules, we have
$$
\vH^\hdot(G,\{H_i\},M \otimes W_\idot)
= 0.
$$
\end{lemma}

\proof{} Let $W_\idot' = \Z \oplus W_\idot$, and let $w:I^{fg} \to
I^{fg}$ be the functor which sends $V_\idot \in I^{fg}$ to $W'_\idot
\otimes V_\idot = V_\idot \oplus (V_\idot \otimes W_\idot)$. Then
$w$ is obviously cofinal in the sense of \cite{KS}, so that the
natural map $\Id \to w$ induced by $\Z \to W'_\idot$ gives an
isomorphism
$$
\lim_\to H^\hdot(G,M \otimes V_\idot) \cong
\lim_\to H^\hdot(G,M \otimes W'_\idot \otimes V_\idot),
$$
where both limits are over $V_\idot \in I^{fg}$. Since $\vH^\hdot$ is a
cohomological functor, this proves the claim.
\endproof

\subsection{Adapted complexes.}
In practice, the index category $I^{fg}$ is still too big. To be
able to compute generalized Tate cohomology more efficiently, we use
the following gadget.

\begin{lemma}\label{ada.lemma}
Assume given a bounded from above complex $P_\idot$ of finitely
generated $\Z[G]$-modules. For any integer $l$, denote by
$F^lP_\idot$ the $(-l)$-th term of the stupid filtration on
$P_\idot$. Then the following conditions are equivalent.
\begin{enumerate}
\item For any of the subgroups $H_i \in \{H_i\}$ and any integer
  $l$, there exists an integer $l' \geq l$ such that the map
$$
F^lP_\idot \to F^{l'}P_\idot
$$
becomes $0$ after restricting to $\D(H_i,\Z)$.
\item For any $V_\idot \in \D_{\{H_i\}}(G,\Z)$, we have
$$
\lim_{\overset{l}{\to}}H^\hdot(G,V_\idot \otimes F^lP_\idot) = 0.
$$
\item For any $V_\idot \in \D_{\{H_i\}}(G,\Z)$, any Grothendieck
  abelian category $\Ab$, and any $M \in \D(G,\Ab)$, we have
$$
\lim_{\overset{l}{\to}}H^\hdot(G,M \otimes V_\idot \otimes
F^lP_\idot) = 0.
$$
\end{enumerate}
\end{lemma}

\proof{} The condition \thetag{iii} contains \thetag{ii} as a
particular case.  By adjunction, \thetag{i} implies \thetag{iii} for
$V_\idot$ of the form $\Ind^{H_i}_G(V_\idot')$, $V'_\idot \in
\D(H_i,\Z)$. Since $\D_{\{H_i\}}(G,\Z)$ consists of direct summands
of sums of such $V_\idot$, this condition on $V_\idot$ can be
dropped. Finally, \thetag{ii} applied to a bounded complex $V_\idot$
of finitely generated $\Z[G]$-modules can be rewritten as follows:
for any $l$, any map $\kappa:V_\idot^* \to F^lP_\idot$ becomes $0$
after composing with the natural map $F^lP_\idot \to F^{l'}P_\idot$
for sufficiently large $l' \geq l$. Applying this to $V_\idot^* =
\Ind^{H_i}_G(F^lP_\idot)$ with $\kappa$ being the adjunction map
yields \thetag{i}.
\endproof

\begin{defn}\label{adapt}
A complex $P_\idot$ of $\Z[G]$-modules is said to be {\em adapted to
  the family $\{H_i\}$} if
\begin{enumerate}
\item $P_i = 0$ for $i < 0$, $P_0 \cong \Z$, and for any $i \geq 0$,
  $P_i$ is a flat $\Z$-modules of finite rank, and is a direct
  summand of a sum of $\Z[G]$-modules induced from one of the
  subgroups $H_i$, and
\item the complex $P_\idot$ satisfies the equivalent conditions of
  Lemma~\ref{ada.lemma}.
\end{enumerate}
\end{defn}

\begin{prop}
Assume given a finite group $G$, a family $\{H_i\}$ of subgroups
$H_i \subset G$, and complex $P_\idot$ of $\Z[G]$-modules adapted to
the family $\{H_i\}$. Then for any $M_\idot \in
\D^b(G,\Ab)$, we have an isomorphism
\begin{equation}\label{adapt.eq}
\vH^\hdot(G,\{H_i\},M_\idot) \cong
\lim_{\overset{l}{\to}}H^\hdot(G,M_\idot \otimes F^lP_\idot),
\end{equation}
where $F^lP_\idot$ is as in Definition~\ref{adapt}~\thetag{ii}, and
this isomorphism is functorial in $M_\idot$.
\end{prop}

\proof{} Denote for the moment
$$
\vH_o^\hdot(G,\{H_i\},M_\idot) \cong
\lim_{\overset{l}{\to}}H^\hdot(G,M_\idot \otimes F^lP_\idot).
$$
Consider the product $I^{fg} \times \N$ of the index category
$I^{fg}$ of \eqref{i.small} and the partially ordered set of
non-negative integers. Computing the double limit first in one
order, then in another, we obtain an isomorphism
\begin{align*}
\lim_{\overset{l}{\to}}\vH^\hdot(G,\{H_i\},M_\idot \otimes F^lP_\idot)
&\cong \lim_{\to}H^\hdot(G,\{H_i\}M_\idot \otimes F^lP_\idot \otimes
  V_\idot)\\
&\cong \lim_{\to}\vH_o^\hdot(G,\{H_i\}M_\idot \otimes V_\idot),
\end{align*}
where the last limit is over $V_\idot \in I^{fg}$, and the
intermediate limit is over $V_\idot \times l \in I^{fg} \times
\N$. By Lemma~\ref{no.tate.prod} and
Definition~\ref{adapt}~\thetag{i}, $\vH^\hdot(G,\{H_i\}M \otimes
F^lP_\idot)$ does not depend on $l$, so that
$$
\lim_{\overset{l}{\to}}\vH^\hdot(G,\{H_i\}M_\idot \otimes
F^lP_\idot) \cong \vH^\hdot(G,\{H_i\}M_\idot).
$$
By Definition~\ref{adapt}~\thetag{ii} in the form of
Lemma~\ref{ada.lemma}~\thetag{iii}, $\vH_o^\hdot(G,\{H_i\}M_\idot
\otimes V_\idot)$ does not depend on $V_\idot \in I^{fg}$, so that
$$
\lim_{\to}\vH_o^\hdot(G,\{H_i\}M_\idot \otimes V_\idot) \cong
\vH_o^\hdot(G,\{H_i\}M_\idot).
$$
This finishes the proof.
\endproof

\begin{exa}
Take as $\{H_i\}$ the family consisting only of the trivial subgroup
$\{e\} \subset G$. A $\Z[G]$-module is induced from $\{e\}$ if and
only if it is free, so that the category $\D^b_{\{H_i\}}(G,\Z)$
consists of direct summands of finite complexes of free
$\Z[G]$-modules -- equivalently, these are the perfect objects in
$\D^b(G,\Z)$. Thus in this case, $\vH^\hdot(G,\{e\},-)$ is the
standard Tate (co)homology functors. Further,
Definition~\ref{adapt}~\thetag{ii} simply means that the complex
$P_\idot$ is acyclic (thus contractible as a complex of abelian
groups), and Definition~\ref{adapt}~\thetag{i} means that $P_i$ is a
finitely generated projective $\Z[G]$-module for any $i \geq
1$. This gives the standard procedure for computing Tate homology:
take a projective resolution $P_\idot$ of the trivial representation
$\Z$ and let $\wt{P}_\idot$ be the cone of the augmentation map
$P_\idot \to \Z$; then
$$
\vH^\hdot(G,V_\idot) \cong \lim_{\overset{l}{\to}}H^\hdot(G,V_\idot
\otimes F^l\wt{P}_\idot)
$$
for any $V_\idot \in \D^b(G,\Z\amod)$.
\end{exa}

In general, there are different ways to construct an adapted complex
for a finite group $G$ and a family $\{H_i\}$; we will present one
construction later in Subsection~\ref{mack.tate}.

\subsection{Tate cohomology and fixed points functors.}

Now assume given a small category $\C$ of Galois type, and consider
the category $\DS_{fs}(\C,\Ab) \cong \DT_{fs}(\C,\Ab)$, as in
Theorem~\ref{galois.comp}. We want to express the gluing data
between the pieces of the semiorthogonal decomposition of
Corollary~\ref{1.corr} in terms of generalized Tate cohomology.

The main result is as follows. Take a morphism $f:c' \to c$ in $\C$,
and consider the corresponding complex $\T^\C_\idot(f)$ of
\eqref{t.n.eq}. Denote by $\Aut(f) \subset \Aut(c')$ the subgroup of
those automorphisms $\sigma:c' \to c'$ that commute with $f$, $f
\circ \sigma = f$. This group acts on every set $\C_n(f)$ by left
composition, that is, a diagram $c' \overset{f_1}{\longrightarrow}
c_1 \to \dots \to c_n \to c$ goes to $c' \overset{f' \circ
\sigma}{\longrightarrow} c_1 \to \dots \to c_n \to c$. Thus
$\Aut(f)$ also acts on $\T^\C_\idot(f)$.  Moreover, for any diagram
$\alpha = [c' \to c_1 \to c] \in \C_1(f)$, let $\Aut(\alpha) \subset
\Aut(f)$ be its stabilizer in the group $\Aut(f)$.

\begin{prop}\label{T.tate}
The complex $\T^\C_\idot(f)$ equipped with the natural action of the
group $\Aut(f)$ is adapted to the family $\{\Aut(\alpha)| \alpha \in
\C_1(f)\}$ in the sense of Definition~\ref{adapt}.
\end{prop}

\proof{} The condition~\thetag{i} of Definition~\ref{adapt} is
obvious: the stabilizer of any diagram $c' \to c_1 \to \dots \to c_n
\to c$ is contained in the stabilizer of the diagram $c' \to c_1 \to
c$ obtained by forgetting $c_i$, $i \geq 2$, so that it suffices to
use the following trivial observation.

\begin{lemma}\label{no.stab}
Assume given a finite group $G$ with a family of subgroups
$\{H_i\}$, $H_i \subset G$, and let $X$ be a finite $G$-set such
that the stabilizer $G_x$ of any point $x \in X$ is contained in one
of the subgroups $H_i$. Then the set $\Z[X]$ with the natural
$G$-action lies in the subcategory $\D^b_{\{H_i\}}(G,\Z\amod)
\subset \D^b(G,\amod)$.\endproof
\end{lemma}

To check the condition~\thetag{ii} in the form of
Lemma~\ref{ada.lemma}~\thetag{ii}, consider the category $\C_c$ of
Corollary~\ref{0.corr}. Let $\C_f$ be the category of diagrams $c'
\to c_1 \to c$ such that the composition map $c' \to c$ is equal to
$f$, and let $\pi^f:\C_f \to \C_c$ be the projection functor which
sends an object $[c' \to c_1 \to c]$ of $\C_f$ to the object $[c_1
\to c]$ of $\C_c$. This is a discrete cofibration. Let $T^f =
\pi^f_!\Z \in \Fun(\C_c,\Z\amod)$. The group $\Aut(f)$ acts on
$T^f$, and we have
$$
\T^\C_\idot(f) \cong \Phi^c_\idot(T^f),
$$
where $\Phi^c_\idot(-)$ is the complex \eqref{phi.eq}. For any
subgroup $H \subset \Aut(f)$ and any $h$-injective complex $V_\idot$
of $\Z[H]$-modules, we have
$$
\lim_{\overset{l}{\to}}\left(V_\idot \otimes F^l\T^\C_\idot(f)\right)^H
\cong \lim_{\overset{l}{\to}}F^l\Phi_\idot(T^f \otimes V_\idot)^H \cong
\Phi_\idot(T^f \otimes V_\idot)^H,
$$
and by Lemma~\ref{phi.idot}, this computes
\begin{equation}\label{to.kill}
\Phi^c\left(V_\idot \otimes T^f\right)^H.
\end{equation}
It suffices to prove that this is trivial whenever $H =
\Aut(\alpha)$ for some $\alpha = [c' \to c_1 \to c] \in
\C_1(f)$. Indeed, fix such an $\alpha$, and denote its component
maps by $f_1:c' \to c_1$, $f_2:c_1 \to c$. Then we have the category
$\C_{c_1}$ and the functor $T^{f_1} \in \Fun(\C_{c_1},\Z\amod)$, and
the group $\Aut(\alpha) = \Aut(f_1)$ acts on $T^{f_1}$. Moreover,
composition with $f_2$ defines functors $\rho:\C_{c_1} \to \C_c$,
$\rho':\C_{f_1} \to \C_f$, and we have a commutative diagram
$$
\begin{CD}
\C_{f_1} @>{\rho'}>> \C_f\\
@V{\pi^{f_1}}VV @VV{\pi^f}V\\
\C_{c_1} @>{\rho}>> \C_c.
\end{CD}
$$
Since both $\C_{f_1}$ and $\C_f$ have an initial object, and $\rho'$
preserves them, we have $L^\hdot\rho'_!\Z \cong \rho'_!\Z \cong \Z$,
so that we have an isomorphism
$$
L^\hdot\rho_!T^{f_1} \cong \rho_!T^{f_1} \cong \rho_!\pi^{f_1}_!\Z
\cong \pi^f_!\rho'_!\Z \cong \pi^f_!\Z \cong T^f,
$$
and this isomorphism is $\Aut(\alpha)$-equivariant. This yields a
quasiisomorphism
$$
\left(V_\idot \otimes T^f\right)^{\Aut(\alpha)} \cong
L^\hdot\rho_!\left(V_\idot \otimes T^{f_1}\right)^{\Aut(\alpha)},
$$
so that by Corollary~\ref{0.corr}, \eqref{to.kill} is of the form
$\phi^cL^\hdot\rho_!E_\idot$ for some $E_\idot \in
\D(\C_{f_1},\Z\amod)$. This is equal to $0$: by definition of the
functor $\phi^c$, we have
$$
\Hom(\phi^cL^\hdot\rho_!E_\idot,M_\idot) \cong
\Hom(L^\hdot\rho_!E_\idot,M_\idot \otimes \ol{T}) \cong
\Hom(E_\idot,\rho^*\ol{T} \otimes M_\idot)
$$
for any $M_\idot \in \D(\Ab)$, and $\rho^*\ol{T} = 0$.
\endproof

\begin{corr}\label{tate.corr}
For any Galois-type category $\C$ and two objects $c,c' \in \C$, the
homology $\HH^\hdot(\Phi^c\Infl^{c'}(M_\idot))$ of the object
$\Phi^c\Infl^{c'}(M_\idot) \in \D(\Ab)$ is given by
\begin{equation}\label{vh.eq}
\HH^\hdot(\Phi^c\Infl^{c'}(M_\idot)) \cong \bigoplus_{f \in
  \C(c',c)}\vH^\hdot(\Aut(f),\{\Aut(\alpha)|\alpha \in
  \C_1(f)\},M_\idot)
\end{equation}
for any $M_\idot \in \D(\br{c'},\Ab)$.\endproof
\end{corr}

To write down a formula for $\wPhi^{c}_\idot \Infl^{c'}$, the gluing
functor in the semi-orthogonal decomposition of
Corollary~\ref{1.corr}, one has to incorporate the natural
$\Aut(c)$-action on $\Phi^c$ into \eqref{vh.eq}; we leave it to the
reader.

\subsection{Invertible objects.}\label{inve.subs}

Let us say that an object $M \in \D$ is a triangulated tensor
category $\D$ is {\em invertible} if the functor $M \otimes -:\D \to
\D$ given by multiplication by $M$ is an equivalence (for example, a
unit object $I \in \D$ is invertible, and so are all its shifts
$I[l]$, $l \in \Z$). For simplicity, for the moment we restrict our
attention to the case $\Ab = \Z\amod$. Keep the assumptions of the
previous subsection, and assume also that $\C$ has a terminal
object, so that $\DS(\C,\Z)$ is a tensor triangulated category. In
addition, assume that $[\C]$ is finite, so that $\DS_{fs}(\C,\Z) =
\DS(\C,\Z)$. As an application of Proposition~\ref{T.tate}, let us
prove the following criterion of invertibility for objects in
$\DT(\C,\Z) \cong \DS(\C,\Z)$.

\begin{prop}\label{inve.prop}
Assume given an object $M \in \D(\C^{opp},\Z)$ such that
\begin{enumerate}
\item the restriction $M(c) = j^{c*}M \in \D(\br{c},\Z)$ with
  respect to the inclusion $j^c:\br{c} \to \C$ is invertible for any
  $c \in \C$, and
\item for any map $f:c' \to c$ in $\C$, equip $M(c)$ with the
  trivial $\Aut(f)$-action, so that the map $M(f):M(c) \to M(c')$ is
  $\Aut(f)$-equivariant. Then the cone of this map lies in
\begin{equation}\label{inve.eq}
\D_{\{\Aut(\alpha)|\alpha \in \C_1(f)\}}(\Aut(f),\Z\amod) \subset
\D^b(\Aut(f),\Z\amod).
\end{equation}
\end{enumerate}
Then the induced object $q^{opp}_!M \in \DS(\C,\Z)$ is invertible.
\end{prop}

\proof{} Let $\Mm:\DT(\C,\Z) \to \DT(\C,\Z)$ be the functor given
by multiplication by $M$, $\Mm(E) = M \otimes E$. By
\eqref{prod.indu.1} of Lemma~\ref{prod.indu.como}, $\Mm$ preserves
supports: if $E \in \DS(\C,\Z)$ is supported in some subset $U
\subset [\C]$, then $\Mm(E)$ is also supported in $U$. Thus as in
the proof of Theorem~\ref{galois.comp}, it suffices to prove that
$$
\Mm:\DT_U(\C,\Z) \to \DT_U(\C,\Z)
$$
is an equivalence for any finite closed $U \subset [\C]$. Take a
maximal element $[c] \in U$ and let $U' = U \setminus \{[c]\}$. Then
$\Mm$ sends the pieces of the semiorthogonal decomposition of
Lemma~\ref{orth.coa} into themselves. By induction, $\Mm$ is an
autoequivalence of $\DT_{U'}(\C,\Z)$, and by the
condition~\thetag{i}, $\Mm$ is also an autoequivalence of the
orthogonal ${}^\perp\DT_{U'}(\C,\Z) \cong \D(\br{c},\Z)$. Again as
in the proof of Theorem~\ref{galois.comp}, by Lemma~\ref{orth.equi}
it suffices to prove that $\Mm$ sends the orthogonal
$\DT_{U'}(\C,\Z)^\perp$ into itself. But every object in this
orthogonal is of the form $\I(E_\idot)$ for some $h$-injective
complex $E_\idot$ of $\Aut(c)$-modules, where $\I$ is as in
\eqref{I.gal}, and by Corollary~\ref{tate.corr} and
Lemma~\ref{no.tate.prod}, the condition~\thetag{ii} then implies
that $\Mm(\I(E_\idot))$ is quasiisomorphic to $\I(E_\idot \otimes
M(c))$.
\endproof

\subsection{The case of Mackey functors.}\label{mack.tate}

We now turn to the case of Mackey functors -- we assume given a
group $G$, and take $\C = O_G$, the category of finite
$G$-orbits. First, we fulfil two earlier promises -- give a
construction of an adapted complex, and prove that inflation is
fully faithful.

\begin{lemma}
For any finite group $G$ with a family of subgroups $\{H_i\}$, there
exists a complex $P_\idot$ of $\Z[G]$-modules adapted to the family
$\{H_i\}$ in the sense of Definition~\ref{adapt}.
\end{lemma}

\proof{} Note that we may assume that together with any subgroup
$H_i \subset G$, the family $\{H_i\}$ contains all the subgroups $H
\subset H_i$ -- adding such subgroups to the family does not change
the conditions of Definition~\ref{adapt} (nor of
Definition~\ref{gen.tate}). Having assumed this, consider the full
subcategory $\C' \subset O_G$ spanned by the one-point orbit $[G/G]$
and all the orbits $[G/H_i]$, $H_i \in \{H_i\}$. For any subgroup
$H_1,H_2 \subset G$, every $G$-orbit in the product $G/H_1 \times
G/H_2$ of of the form $G/H'$ for some $H' \subset H_1$. Therefore
the category $\C'$ is of Galois type, so that
Proposition~\ref{T.tate} applies. Take as $f$ the unique map
$f:[G/\{e\}] \to [G/G]$; then $G = \Aut(f)$, $\Aut(\alpha) \in
\{H_i\}$ for any $\alpha \in \C'_1(f)$, and $\T^{\C'}_\idot(f)$ is
the required adapted complex.
\endproof

\begin{lemma}\label{infl.ff}
Assume given a normal subgroup $H \subset G$ with the quotient $N =
G/H$. Then the inflation functor
$$
\wh{\Infl}^N_G:\DM(N,\Ab) \to \DM(G,\Ab)
$$
is fully faithful.
\end{lemma}

\proof{} Let $U \subset [O_G]$ be the set of orbits $[G/H']$ such
$[G/H] \geq [G/H']$, that is, there exists a $G$-equivariant map
$G/H \to G/H'$, that is, $H' \subset G$ contains a conjugate of $H
\subset G$. Since $H \subset G$ is assumed to be normal, $H$ itself
must be contained in $H'$. Therefore the orbits $G/H'$, $[G/H'] \in
U$ are in one-to-one correspondence with orbits $N/(H'/H)$ of the
quotient group $N = G/H$, and the full subcategory in $O_G$ spanned
by objects with classes in $U$ is equivalent to the category
$O_N$. Moreover, it is easy to see that the embedding $O_N \subset
O_G$ is compatible with the graded $A_\infty$-coalgebras $\T_\idot$
-- for any map $f$ in $O_N$, we have
$$
\T^{O_N}_\idot(f) \cong \T^{O_G}_\idot(f).
$$
Thus $\DM(N,\Ab) \cong \DT_U(O_G,\Ab) \subset \DM(G,\Ab)$, and
$\wh{\Infl}^N_G$ is the full embedding $\DT_U(O_G,\Ab) \subset
\DT(O_G,\Ab) = \DM(G,\Ab)$.
\endproof

Now we make the following easy observation.

\begin{lemma}\label{redd}
Assume given a Galois-type category $\C$ with the corresponding
$\C$-graded $A_\infty$-coalgebra $\T^\C_\idot$ of
Definition~\ref{T.defn}, another $\C$-graded $A_\infty$-coalgebra
$\T_\idot'$, and an $A_\infty$-map $\nu:\T'_\idot \to \T^\C_\idot$
such that
\begin{enumerate}
\item if $f$ is invertible, the map $\nu:\T'_\idot(f) \to
  \T^\C_\idot(f)$ is an isomorphism,
\item for a non-invertible $f:c' \to c$ and any $\Z[\Aut(f)]$-module $V$,
  the map $\nu$ induces an isomorphism
\begin{equation}\label{T.T'}
\lim_{\overset{l}{\to}}H^\hdot(\Aut(f),V_\idot \otimes
F^l\T'_\idot(f)) \cong
\lim_{\overset{l}{\to}}H^\hdot(\Aut(f),V_\idot \otimes
F^l\T^\C_\idot(f)),
\end{equation}
where $F^\hdot\T^\C_\idot$, $F^\hdot\T'_\idot$ are the stupid
filtrations.
\end{enumerate}
Then the corestriction functor $\nu^*$ induces an equivalence
between the derived category $\DT_{fs}(\C,\Ab)$ of finitely
supported $A_\infty$-comodules over $\T^\C_\idot$ and the category
$\DT'_{fs}(\C,\Ab)$ of finitely supported $A_\infty$-comodules over
$\T'_\idot$.
\end{lemma}

\proof{} The same as Theorem~\ref{galois.comp} and
Proposition~\ref{inve.prop}: since $\C$ is lattice-like, both
categories have filtrations by support and the corresponding
semi-orthogonal decompositions, and the corestriction functor
$\nu^*$ is compatible with these decompositions. Then \thetag{i}
insures that $\nu^*$ is an equivalence of the associated graded
pieces, and \thetag{ii} insures that $\nu^*$ is compatible with the
gluing.
\endproof

The condition \thetag{ii} of this Lemma is satisfied, for example,
when $\T'_\idot(f)$ is adapted to $\{\Aut(\alpha)\}$ for any $f$, as
in Proposition~\ref{adapt}. However, it can also be satisfied for
other reasons: for instance, if the right-hand side of \eqref{T.T'}
is equal to $0$, the left-hand side may just be trivial,
$\T'_\idot(f)=0$. This is especially useful in the Mackey functors
case, for the following reason.

\begin{lemma}\label{no.tate}
Let $G$ be a finite group, and let $\{H_i\}$ be the family of all
proper subgroups $H \subset G$.
\begin{enumerate}
\item The generalized Tate cohomology $\vH^\hdot(G,\{H_i\},V)$ with
  coefficients in any $\Z[G]$-module $V$ is trivial unless $G$ is a
  $p$-group for some prime $p$, and in this case, it is annihilated
  by $p$.
\item If $G$ is a cyclic prime group, $G = \Z/p^n$ for some prime
  $p$ and integer $n \geq 2$, then $\vH^\hdot(G,\{H_i\},V)=0$ for
  any $\Z[G]$-module $V$.
\end{enumerate}
\end{lemma}

\proof{} For a proper subgroup $H \subset G$, we have the natural
induction and coinduction maps $\Z \to \Ind_G^H(\Z)$, $\Ind_G^H(\Z)
\to \Z$, and their composition is multiplication by the index
$|G/H|$ of the subgroup $H$. Therefore multiplication by $|G/H|$ is
trivial in the quotient category $\D(G,\{H_i\},\Z\amod)$, hence
annihilates generalized Tate cohomology. If $G$ is not a $p$-group,
then the greatest common denominator of the indices of its Sylow
subgroups is $1$, so that $1$ also annihilates Tate cohomology. If
$G$ is a $p$-group, then it contains a subgroup of index $p$. This
proves \thetag{i}.

For \thetag{ii}, note that for any $n$, the second cohomology group
$H^2(\Z/p^n\Z,\Z)$ contains a canonical periodicity element $u_n$
represented by Yoneda by the exact sequence
$$
\begin{CD}
0 @>>> \Z @>>> \Z[\Z/p^n\Z] @>{1-\sigma}>> \Z[\Z/p^n\Z] @>>> Z @>>>
0,
\end{CD}
$$
where $\sigma$ is the generator of the cyclic group $\Z/p^n\Z$. For
any positive integer $n' < n$, we have a natural quotient map
$q_{n,n'}:\Z/p^n\Z \to \Z/p^{n'}\Z$, and one easily checks that
$q_{n,n'}^*(u_{n'}) = p^{n-n'}u_n$. In the quotient category
$\D(G,\{H_i\},\Z\amod)$, all the objects $q_{n,n'}^*\Z[\Z/p^{n'}\Z]$
become trivial, so that all the maps $q_{n,n'}u_{n'}$ are
invertible. If $n \geq 2$, this implies that $p$ is
invertible. Since $p$ annihilates Tate cohomology by \thetag{i},
this is only possible if $\vH^\hdot(G,\{H_i\},-)=0$.
\endproof

\begin{corr}
Assume that $G$ is a finite group whose order is invertible in the
category $\Ab$. Then
$$
\DM(G,\Ab) \cong \bigoplus_{H \subset G}\D(N_H/H,\Ab),
$$
where the sum is over all the conjugacy classes of subgroups $H
\subset G$, and $\D(N_H/H,\Ab)$ is the derived category of the
representations of the quotient $N_H/H$ of the normalizer $N_H$ of
the subgroup $H \subset G$ by $H$ itself.
\end{corr}

\proof{} Since $G$ is finite, everything has finite support,
$\DS(O_G,\Ab) \cong \DS_{fs}(O_G,\Ab) \cong
\DT_{fs}(O_G,\Ab)$. Moreover, every map $f$ between two $G$-orbits
is isomorphic to the quotient map $[G/H] \to [G/H']$ for some
subgroups $H \subset H' \subset G$, the automorphisms group of an
object $[G/H]$ is precisely $N_H/H$, and for any map $f:[G/H] \to
[G/H']$, any proper subgroup in the automorphisms group $\Aut(f)$ is
of the form $\Aut(\alpha)$ for some $\alpha \in \C_1(f)$. Now take
$\T'_\idot(f) = \T^\C_\idot(f)$ for invertible $f$, $\T'_\idot(f) = 0$
otherwise, and apply Lemma~\ref{redd}. The condition \thetag{ii} is
satisfied, since the left-hand side of \eqref{T.T'} is isomorphic to
$\vH^\hdot(\Aut(f),\{\Aut(\alpha)\},-)$, and this is trivial by
Lemma~\ref{no.tate}.
\endproof

In the ordinary, non-derived Mackey functor case, this is a Theorem
of Thevenaz \cite{Th1}.

Now take $\Ab$ arbitrary, $G$ possibly infinite, $\C = O_G$. In this
case, we can still use Lemma~\ref{redd} and Lemma~\ref{no.tate} to
simplify the $A_\infty$-coalgebra $\T^\C_\idot$. Namely, any map
between two finite $G$-orbits is isomorphic to the quotient map
$f:[G/H] \to [G/H']$, where $H \subset H' \subset G$ are subgroups
of finite index. Denote by $\Ind(f)$ the index of $H$ in $H'$. Then
one can replace $\T^\C_\idot$ with the $\C$-graded
$A_\infty$-coalgebra $\T'_\idot$ defined as
$$
\T'_\idot(f) = 
\begin{cases}
\T^\C_\idot(f), &\quad \Ind(f) = p^n \text{ for some prime }p,\\
0, &\text{otherwise}.
\end{cases}
$$
Again, by Lemma~\ref{no.tate}~\thetag{i} this will not change the
category $\DT_\idot(\C,\Ab)$.

\subsection{Cyclic group.}

We can go even further in the important special case $G = \Z$, the
infinite cyclic group. For every integer $n$, choose a projective
resolution $P^n_\idot$ of the trivial $\Z[\Z/n\Z]$-module $\Z$. For
example, we can take the standard periodic resolution
\begin{equation}\label{perio}
\begin{CD}
@>{\id-\sigma}>> \Z[\Z/n\Z] @>{\id+\sigma+\dots+\sigma^{n-1}}>>
  \Z[\Z/n\Z] @>{\id-\sigma}>> \Z[\Z/n\Z],
\end{CD}
\end{equation}
where $\sigma \in \Z/n\Z$ is the generator. Let $\wt{P}^n_\idot$ be
the cone of the augmentation map $P^n_\idot \to \Z$ (that is,
$\wt{P}_0 = \Z$, and $\wt{P}^n_l = P^n_{l-1}$ for $l \geq 1$). For
any morphism $f:c' \to c$ in $O_G$ from $c' = [\Z/n'\Z]$ to $c
=[Z/n\Z]$, define the complex $\wt{T}_\idot(f)$ by
$$
\wt{\T}_\idot(f)=
\begin{cases}
\Z, &\quad \Ind(f) = 1,\\
\wt{P}^{n'}_\idot, &\quad \Ind(f) = \frac{n'}{n} = p \text{ is prime},\\
0, &\text{otherwise}.
\end{cases}
$$
To extend this to an $A_\infty$-coalgebra, we need to define
co-multiplication maps, and the only possible non-trivial
comultiplication maps are those that encode $\Aut(c') \times
\Aut(c)$-action on $\bigoplus_{f \in O_G(c',c)}\wt{\T}_\idot(f)$ in
the case $\Ind(f)=p$, $n' = np$. For any $c = [\Z/n\Z] \in O_G$, we
have $\Aut(c) = \Z/n\Z$. More generally, for any $c' = [\Z/n'\Z]$,
$c = [Z/n\Z]$, $n'$ divisible by $n$, we have $O_G(c',c)=\Z/n\Z$,
and the group $\Aut(c') \times \Aut(c)$ acts transitively on this
set; we have
$$
O_G(c',c) = (\Aut(c') \times \Aut(c))/\Aut(c'),
$$
where the embedding $\Aut(c') \subset \Aut(c') \times \Aut(c)$ is
the product of the identity map $\Aut(c') \to \Aut(c')$ and the
natural projection $\Aut(c') = \Z/n'\Z \to \Aut(c) = \Z/n\Z$. In the
particular case $n' = np$, we take the given $\Z/n'\Z$-action on
$\wt{P}^{n'}_\idot$, and induce the $\Aut(c') \times \Aut(c)$-action
on $\bigoplus_{f \in O_G(c',c)}\wt{\T}_\idot(f)$ via the identification
$$
\bigoplus_{f \in O_G(c',c)}\wt{\T}_\idot(f) = \bigoplus_{f \in
  O_G(c',c)}\wt{P}_\idot^{n'} \cong \Ind_{\Aut(c') \times
  \Aut(c)}^{\Aut(c')}\wt{P}^{n'}_\idot.
$$
All the higher maps $b_n$, $n \geq 3$ are set to be trivial, $b_n = 0$.

\begin{prop}
With the notation above, the category $\DS_{fs}(O_G,\Ab)$ is
equivalent to the derived category of finitely supported
$\Ab$-valued $A_\infty$-comodules over the $O_G$-graded
$A_\infty$-coalgebra $\wt{\T}_\idot$.
\end{prop}

\proof{} First, we construct a map of $O_G$-graded
$A_\infty$-coalgebras $\nu:\wt{\T}_\idot \to \T^{O_G}_\idot$. If
$\Ind(f)=1$, in other words $f$ is invertible, we have
$\wt{\T}_\idot(f) = \T^{O_G}_\idot(f) = \Z$ by definition. If
$\Ind(f)$ is not prime, $\wt{\T}_\idot(f) = 0$, so that there is
nothing to define. It remains to define an $\Aut(c') \times
\Aut(c)$-equivariant map
$$
\nu:\wt{\T}_\idot(f) \to \T^{O_G}_\idot(f)
$$
for any map $f:c' \to c$ of prime index $\Ind(f) = p$, $c =
[\Z/n\Z]$, $c' = \Z/pn\Z$. We have $\Aut(c') = \Z/np\Z$, and by
adjunction, defining $\nu$ is equivalent to defining an
$\Z/pn\Z$-equivariant map
$$
\overline{\nu}:\wt{P}^{np}_\idot \to \T_\idot(c',c) = \bigoplus_{f \in
  O_G(c',c)}\T^{O_G}_\idot(f).
$$
Both complexes are acyclic. In degree $0$, $\wt{P}^{np}_0 = \Z$,
$\T^{O_G}_0(c',c) = \Z[O_G(c',c)]$ is a free $\Z/np\Z$-module, and
we take as $\overline{\nu}$ the natural map which identifies $\Z$
with the subgroup of $\Z/np\Z$-invariants in $\T^{O_G}_0(c',c)$. In
higher degrees, $\wt{P}^{np}_{\idot+1}$ is a resolution of $\Z$, and
$\T^{O_G}_{\idot+1}(c',c)$ is a resolution of
$\T^{O_G}_0(c',c)$. But the resolution $\wt{P}^{np}_{\idot+1}$ is by
definition projective; therefore the given map
$\overline{\nu}:\wt{P}^{np}_0 \to \T^{O_G}_0(c',c)$ extends to a map
of resolutions $\overline{\nu}:\wt{P}^{np}_{\idot+1} \to
\T^{O_G}_{\idot+1}(c',c)$.

To finish the proof, it remains to check the conditions of
Lemma~\ref{redd}. The condition \thetag{i} is satisfied by
definition. As for \thetag{ii}, it is satisfied unless $\Ind(f)$ is
prime by virtue of Lemma~\ref{no.tate}~\thetag{ii}. And if
$\Ind(f)=p$ is prime, both $\T^{O_G}_\idot(f)$ and $\wt{\T}_\idot$
are complexes of $\Z[\Aut(f)]$-modules adapted to the family of
Lemma~\ref{no.tate} (which in this case consist of the unique proper
subgroup in $\Aut(f) = \Z/p\Z$, namely, the unit subgroup $\{e\}
\subset \Aut(f)$). This also implies \thetag{ii} of
Lemma~\ref{redd}.
\endproof

By Lemma~\ref{infl.ff}, this description of the category
$\DS_{fs}(O_\Z,\Ab)$ also yields a description of the category
$\DM(\Z/n\Z,\Ab)$ for every finite cyclic group $\Z/n\Z$ -- one just
restricts to the orbits $[\Z/l\Z]$ with $l$ a divisor of $n$ (for a
finite group, the ``finite support'' condition becomes vacuous). In
the simplest possible case $G = \Z/p\Z$, $p$ prime, this
boils down to the following.

\begin{corr}
Let $G = \Z/p\Z$. Then the category $\DM(G,\Ab)$ is obtained by
inverting quasiisomorphisms from the category of triples $\langle
V_\idot,W_\idot,\rho \rangle$ of a complex $W_\idot$ of objects in
$\Ab$, a complex $V_\idot$ of $G$-representations in $\Ab$, and
a map
$$
\rho:W_\idot \to \wh{C}^\hdot(G,V_\idot),
$$
where $\wh{C}^\hdot(G,V_\idot)$ is the standard $2$-periodic complex
which computes Tate cohomology $\vH^\hdot(G,V_\idot)$ by means of
the resolution \eqref{perio}.
\end{corr}

\proof{} Clear.\endproof

\section{Relation to stable homotopy.}\label{spectra}

To finish the paper, we fix a finite group $G$ and explain the
relation between the category of derived Mackey functors $\DM(G)$
and $G$-equivariant stable homotopy theory. We will only give a
skeleton exposition so as to present the general principles. In
particular, we will restrict our attention to finite CW spectra when
giving proofs. To keep the text accessible to a person without
knowledge of equivariant stable homotopy, we do recall some of the
basics; for further information, the reader should construct the
original references, which are \cite{br} for
Subsection~\ref{spa.subs} and \cite{may1} for
Subsection~\ref{spe.subs}.

\subsection{Spaces.}\label{spa.subs}

The homotopy category $\Hom(G)$ of $G$-spaces is defined in a
straightforward way: objects are topological spaces equipped with a
continuous action of the group $G$, maps between them are continuous
$G$-equivariant maps, a homotopy between two maps $f,f':X \to Y$ is
a continuous $G$-equivariant map $F:X \times I \to Y$ with $F=f$ on
$X \times \{0\}$ and $F=f'$ on $X \times \{1\}$, where $I$ is the
unit interval $[0,1]$ with the standard topology and trivial
$G$-action, morphisms in the homotopy category are homotopy classes
of maps. Note that in this category, for any subgroup $H \subset G$
the functor $X \mapsto X^H$ which sends a space $X$ to the space
$X^H$ of $H$-invariant points is well-defined. As usual, to get a
useful theory one passes to the category of pointed spaces, and
restricts one's attention to the full subcategory spanned by CW
complexes. Base points are assumed to be $G$-invariant. One also
needs some compatibility condition between the $G$-action and the CW
structure; a convenient notion is that of a $G$-CW {\em complex} ---
this is a CW complex $X$ equipped with a continuous $G$-action such
that for any $g \in G$, the fixed points subset $X^g \subset X$ is a
subcomplex (that is, a union of cells). From now on, we will by
abuse of notation denote the homotopy category of pointed $G$-CW
complexes by $\Hom(G)$.

For any subgroup $H \subset G$ and any $G$-CW complex $X$, the fixed
points subset $X^H \subset X$ is obviously a CW complex. The
geometric realization of a simplicial $G$-set is a $G$-CW complex;
therefore any $G$-set which is homotopy equivalent to a
$G$-equivariant CW complex is also equivalent to a $G$-CW
complex. For any two pointed $G$-spaces $X$, $Y$, their smash
product $X \wedge Y$ is also a pointed $G$-space, and if $X$ and $Y$
are $G$-CW complexes, then so is their product $X \wedge Y$.

Equivalently, the fixed points subset $X^H \subset X$ of a $G$-space
$X$ can be described as follows: take the orbit $[G/H]$, treat it as
a $G$-space with the discrete topology, and consider the space of
continuous maps $\Maps([G/H],X)$ with a natural topology on it. Then
we have
\begin{equation}\label{fixed.maps}
X^H \cong \Maps([G/H],X).
\end{equation}
This shows that for any $G$-space $X$, the correspondence $H \mapsto
X^H$ is a actually a functor from the category $O_G^{opp}$ opposite
to the category of finite $G$-orbits to the category of topological
spaces. If $X$ is a $G$-CW complex, then this is a functor from
$O_G^{opp}$ to the category of CW complexes and cellular maps
between; if the $G$-CW complex $X$ is finite, then so are all the
$X^H$.

Now, for any CW complex $X$, we have a natural cellular chain
complex $C_\idot(X,\Z)$, and this construction is functorial with
respect to cellular maps. Therefore for any $G$-CW complex, we have
a natural functor
$$
[G/H] \mapsto C_\idot(X^H,\Z)
$$
from the category $O_G^{opp}$ to the category of complexes of
abelian groups. Let us denote the corresponding complex of functors
in $\Fun(O_G^{opp},\Z\amod)$ by $C_\idot^G(X,\Z)$. This complex was
constructed and studied by Bredon \cite{br}; in his terminology,
functors from $O_G^{opp}$ to abelian groups are called {\em
coefficient systems}, so that $C_\idot^G(X,\Z)$ is a complex of
coefficient systems. We note that up to a quasiisomorphism, it does
not depend on the particular CW model for an object $X \in
\Hom(G)$. The embedding of the fixed point $\ppt \to X$ induces an
injective map $\Z \cong C_\idot^G(\ppt,\Z) \to C_\idot^G(X,\Z)$; let
$\wt{C}_\idot^G(X,\Z) = C_\idot^G(X,\Z)/\Z$ be the quotient. Again,
up to a quasiisomorphism, it only depends on the class of $X$ in
$\Hom(G)$.

\begin{defn}\label{naive.hom.defn}
The complex $\wt{C}_\idot^G(X,\Z)$ is called the {\em naive reduced
$G$-\-equ\-i\-va\-ri\-ant cellular chain complex} of the CW complex
$X$. The corresponding object $\wt{C}_\idot^G(X,\Z) \in
\D(O_G^{opp},\Z\amod)$ is called the {\em naive reduced
$G$-\-equ\-i\-va\-ri\-ant homology} of $X$ considered as an object
in $\Hom(G)$.
\end{defn}

\begin{lemma}\label{prod.hom}
For any two objects $X,Y \in \Hom(G)$, we have a natural
quasiisomorphism
$$
\wt{C}_\idot^G(X \wedge Y,\Z) \cong \wt{C}_\idot^G(X,\Z) \otimes
\wt{C}_\idot^G(Y,\Z),
$$
where the tensor product in the right-hand side is the pointwise
tensor product in the category $\D(O_G^{opp},\Z\amod)$.
\end{lemma}

\proof{} The corresponding statement ``without the group $G$'' is
completely standard. Since the tensor product in
$\D(O_G^{opp},\Z\amod)$ is pointwise, it suffices to show that $(X
\times Y)^H = X^H \times Y^H$; this is obvious from
\eqref{fixed.maps}.
\endproof

\subsection{Stabilization.}\label{spe.subs}

Recall that for the usual finite pointed CW complexes, the
Spanier-Whitehead stable homotopy category $\SW$ is defined as
follows: object are again finite pointed CW complexes, maps from $X$
to $Y$ are given by
\begin{equation}\label{sp.naive}
\Hom(X,Y) = \lim_{\overset{i}{\to}}[\Sigma^i X,\Sigma^i Y],
\end{equation}
where $[-,-]$ means the set of homotopy classes of maps, and
$\Sigma$ is the suspension functor, so that $\Sigma^iX = S^i \wedge
X$, where $S^i$ is the $i$-sphere.

By definition, the suspension functor $\Sigma$ is fully faithful on
the category $\SW$; in fact the definition is designed exactly to
achieve this. The triangulated {\em category of spectra}, denoted
$\Sthom$, is obtained by formally inverting the suspension, thus
making it not only fully faithful but an autoequivalence. The
passage from $\SW$ to $\Sthom$ involves a somewhat delicate limiting
procedure due to Boardman; I do not feel qualified to re-count it
here, and refer the reader to any of the standard textbooks
(e.g. \cite{adams}) for the definition of $\Sthom$ and further
discussion. We have a natural fully faithful embedding $\SW \subset
\Sthom$ which is compatible with the suspension and is usually
denoted by $\Sigma^\infty$. The category $\sthom$ of {\em finite CW
spectra} is given by
$$
\sthom = \bigcup_i \Sigma^{-i}(\Sigma^\infty(\SW)) \subset \Sthom.
$$
This is the smallest full triangulated subcategory in $\Sthom$
containing $\SW$. Sometimes, e.g. in \cite[1.1.5]{BBD}, it is this
category which is called ``the stable homotopy category''. It
can be equivalently described as the category of pairs $\langle X,n
\rangle$ of a finite CW complex $X$ and an integer $n$, with maps
from $\langle X,n \rangle$ to $\langle Y,m \rangle$ given by
$$
\lim_{\overset{l \geq \max(n,m)}{\to}}[\Sigma^{l-n}X,\Sigma^{l-m}Y].
$$
We should note here that the category $\sthom$ is certainly not
sufficient for many applications in topology, since almost all
spectra representing interesting generalized homology theories are
not finite CW spectra.

The cellular chain complex $\wt{C}^\hdot(-,\Z)$ is compatible with
suspensions, thus descends to the homology functor from $\sthom$ to
the derived category $\D(\Z\amod)$ which we denote by
$$
h:\sthom \to \D(\Z\amod);
$$
explicitly, $h(\langle X,n\rangle)$ is represented by the complex
$\wt{C}^\hdot(X)[n]$.

For $G$-spaces, one can repeat the Spanier-Whitehead construction
literally; this results in a ``naive $G$-equivariant
Spanier-Whitehead category'' $\SW^{naive}(G)$ -- objects are
$G$-spaces $X \in \Hom(G)$, maps are given by the same formula
\eqref{sp.naive} as in the non-equivariant case.

However, there is a more interesting and more natural
option. Namely, the $i$-sphere $S^i$ is the one-point
compactification of the $i$-dimensional real vector space $V$; in
the equivariant world, this vector space should be allowed to carry
a non-trivial $G$-action. To make sense of the direct limit in
\eqref{sp.naive}, one fixes a so-called {\em complete $G$-universe},
that is, an infinite dimensional representation $U$ of the group $G$
which has an invariant inner product and is ``large enough'' in the
sense that every finite-dimensional inner-product representation $V$
of the group $G$ occurs in $U$ countably many times. Then the limit
in \eqref{sp.naive} should be taken over all inner product
$G$-invariant subspaces $V \subset U$:
$$
[X,Y]_U = \lim_{\overset{V \subset U}{\to}}[\Sigma^V X,\Sigma^V Y],
$$
where $\Sigma^V$ is the $V$-suspension functor given by $\Sigma^V(X)
= S^V \wedge X$, and $S^V \in \Hom(G)$ is the one-point
compactification of $V$, with base point at infinity. This
construction works for an arbitrary {\em compact} topological group
$G$. However, when $G$ is a finite group with discrete topology,
there is an obvious preferred choice of the universe: we can take $U
= R^{\oplus \infty}$, the sum of a countable number of copies of the
regular representation $R = \RR[G]$. This results in the following
definition.

\begin{defn}\label{gen.sp}
The ``genuine'' {\em $G$-equivariant Spanier-Whitehead category}
$\SW(G)$ has finite pointed $G$-CW complexes as objects, and maps
from $X$ to $Y$ are given by
$$
\Hom(X,Y) = \lim_{\overset{i}{\to}}[\Sigma^{iR} X, \Sigma^{iR} Y],
$$
where $[-,-]$ stands for the set of morphisms in $\Hom(G)$,
$R=\RR[G]$ is the regular representation of the group $G$, and
$\Sigma^{iR}$ by abuse of notation denotes $(\Sigma^R)^i$.
\end{defn}

Among other representations of $G$, $R$ contains the trivial one, so
that passing to the limit in Definition~\ref{gen.sp} includes taking
the limit \eqref{sp.naive}, and we have a natural tautological
functor
\begin{equation}\label{change.of.universe}
i_!:\SW^{naive}(G) \to \SW(G)
\end{equation}
which is identical on objects. However, this is most certainly not
identical on morphisms, thus not an equivalence.

To pass to the spectra, one again inverts the suspension, either in
the genuine category $\SW(G)$ or in the naive category
$\SW^{naive}(G)$; this results in the so-called {\em category of
naive $G$-spectra} $\Sthom^{naive}(G)$ and the {\em category of
genuine $G$-spectra} $\Sthom(G)$. Both are triangulated
categories. For the precise definitions, I refer the reader to the
book \cite{may1}, the standard reference on the subject. The
suspension functor $\Sigma$ is an autoequivalence on
$\Sthom^{naive}(G)$, and all the suspension functors $\Sigma^V$, $V
\subset R^{\oplus \infty}$, are autoequivalences on $\Sthom(G)$. We
again have full embeddings $\Sigma^\infty_{naive}:\SW^{naive}
\subset \Sthom^{naive}(G)$, $\Sigma^\infty:\SW \subset \Sthom(G)$,
and we can define the triangulated categories of finite $G$-spectra
by
$$
\sthom^{naive}(G) = \bigcup_i
\Sigma^{-i}(\Sigma^\infty_{naive}(\SW(G))) \subset \Sthom^{naive}(G)
$$
and
$$
\sthom(G) = \bigcup_i \Sigma^{-iR}(\Sigma^\infty(\SW(G))) \subset
\Sthom(G),
$$
where $\Sigma^{-iR}$ stands for $(\Sigma^R)^{-i}$. As in the
non-equivariant case, one can equivalently describe these categories
as the categories of pairs $\langle X,n \rangle$, $X \in \Hom(G)$,
$n \in \Z$, where $\langle X,n \rangle$ corresponds to
$\Sigma^{-i}(X) \in \Sthom^{naive}(G)$ in the naive case, and to
$\Sigma^{-iR}(X) \in \Sthom(G)$ in the genuine case.

\subsection{Equivariant homology.}

The Bredon equivariant homology functor
$\wt{C}^G_\idot(-,\Z):\Hom(G) \to \D(O_G^{opp},\Ab)$ of
Definition~\ref{naive.hom.defn} obviously extends to the naive
category $\SW^{naive}(G)$, and then further to the category of
finite spectra, so that we have a natural functor
\begin{equation}\label{h.g.naive}
\sthom^{naive}(G) \to \D(O_G^{opp},\Ab)
\end{equation}
which we denote by $h^G_{naive}(-)$. However, a moment's reflection
shows that it does {\em not} extend to the genuine category
$\sthom(G)$. This is where the Mackey functors come in.

\begin{defn}\label{true.hom.defn}
The {\em $G$-equivariant homology} $h^G(\Sigma^\infty X,\Z)$
of a finite point\-ed $G$-CW complex $X \in \Hom(G)$ is the induced
derived Mackey functor
$$
h^G(X,\Z) = q^{opp}_!\wt{C}^G_\idot(X,\Z) \in \DM(G),
$$
where $q^{opp}_!:\D(O_G^{opp},\Z\amod) \to \DS(O_G,\Z\amod) =
\DM(G)$ is the induction functor of Subsection~\ref{der.indu}.
\end{defn}

\begin{lemma}\label{prod.hom.Mackey}
For any two finite pointed $G$-CW complexes $X$, $Y$, we have a natural
isomorphism
$$
h^G(X \wedge Y,\Z) \cong h^G(X,\Z) \otimes h^G(Y,\Z) \in \DM(G).
$$
\end{lemma}

\proof{} A combination of Lemma~\ref{prod.hom} and
\eqref{prod.indu.2} of Lemma~\ref{prod.indu.como}.
\endproof

\begin{prop}\label{sph.inve}
For any finite-dimensional representation $V$ of the group $G$, the
object
$$
h^G(S^V,\Z) \in \DM(G)
$$
is invertible in the sense of Subsection~\ref{inve.subs}.
\end{prop}

\proof{} Since $h^G(S^V,\Z) = q_!^{opp}\wt{C}_\idot^G(S^V,\Z)$, we
can use the criterion of Proposition~\ref{inve.prop}. To check the
condition \thetag{i}, note that for any subgroup $H \subset G$, we
have
$$
\wt{C}_\idot^G(S^V,\Z)([G/H]) \cong \wt{C}_\idot((S^V)^H,\Z),
$$
and since $(S^V)^H = S^{V^H}$ is a sphere, its reduced homology is
$\Z[\dim_{\RR}V^H]$. To check the condition \thetag{ii}, fix two
subgroups $H_1 \subset H_2 \subset G$, let $S_i = (S^V)^{H_i} =
S^{V^{H_i}}$, $i=1,2$, let $f:[G/H_1] \to [G/H_2]$ be the quotient
map, let $\wt{f}:S_2 \to S_1$ be the corresponding cellular
embedding, and let $W = \Aut(f) \subset \Aut[G/H_1]$. Since for any
proper subgroup $W' \subset W$, we have a factorization
$$
\begin{CD}
[G/H_1] @>>> [G/H_1]/W' @>>> [G/H_2]
\end{CD}
$$
of the map $f$ through some $G$-orbit $[G/H_1]/W'$, every proper
subgroup $W' \subset W$ appears as $\Aut(\alpha)$ in
\eqref{inve.eq}. The group $W$ itself appears only if the map
$[G/H_1]/W \to [G/H_2]$ is not an isomorphism, or in other words,
$S_2 \subset S_1^W$ is a proper inclusion. We have to show that
the relative homology object
$$
\wt{C}_\idot(S_1/S_2,\Z) \in \D^b(W,\Ab)
$$
lies in the subcategory $\D_{\{\Aut(\alpha)\}}(W,\Ab) \subset
\D^b(W,\Ab)$. If $S_2 \neq S_1^W$, $W$ appears in the family
$\{\Aut(\alpha)\rangle$, so that $\D_{\{\Aut(\alpha)\}}(W,\Ab) =
\D^b(W,\Ab)$ and there is nothing to prove; thus we may assume $S_2
= S_1^W$. But since by definition, the stabilizer of any non-trivial
cell in the quotient $S_1/S_2^W$ is a proper subgroup $W$, the claim
then immediately follows from Lemma~\ref{no.stab}.
\endproof

\begin{corr}\label{g.homol}
The $G$-equivariant homology of Definition~\ref{true.hom.defn}
extends to a well-defined functor
$$
h^G(-,\Z):\sthom(G) \to \DM(G).
$$
\end{corr}

\proof{} For any finite-dimensional representation $V$ of the group
$G$, let $\sigma^V:\DM(G) \to \DM(G)$ be the functor given by
$$
\sigma^V(M) = M \otimes h^G(S^V,\Z).
$$
Then by Lemma~\ref{prod.hom.Mackey}, we have $h^G(\Sigma^VX,\Z)
\cong \sigma^V(h^G(X,\Z))$ for any finite $G$-CW complex $X \in
\Hom(G)$, and by Proposition~\ref{sph.inve}, $\sigma^V$ is an
autoequivalence. By definition, every object $X' \in \sthom^G$ is of
the form $X' = (\Sigma^V)^{-1}(\Sigma^\infty(X))$ for some such $V$
and and $X$; the desired extension is then given by
$$
h^G(X',\Z) = (\sigma^V)^{-1}(h^G(X,\Z)),
$$
and by Lemma~\ref{prod.hom.Mackey}, this does not depend on the
choice of the identification $X' =
(\Sigma^V)^{-1}(\Sigma^\infty(X))$.
\endproof

\subsection{A dictionary.}

Corollary~\ref{g.homol} allows one to compare notions from the
$G$-equivariant stable homotopy theory with those of the theory of
Mackey functor. We finish the section, and indeed the whole paper,
with a short dictionary saying what should correspond to what. All
the material on equivariant stable homotopy is taken from
\cite{may1}. Personally, I also find very useful the brief
introduction to \cite{HM}.

First of all, the tautological functor \eqref{change.of.universe}
obviously extends to finite spectra, and by definition, we have a
commutative diagram of triangulated categories and triangulated
functors
\begin{equation}\label{comma}
\begin{CD}
\sthom^{naive}(G) @>{i_!}>> \sthom(G)\\
@V{h^G_{naive}}VV @VV{h^G}V\\
\D(O_G^{opp},\Z) @>{q^{opp}_!}>> \DM(G).
\end{CD}
\end{equation}
The smash product on the category $\Hom(G)$ extends to the
categories of finite spectra $\sthom^{naive}(G)$ and $\sthom(G)$, so
that they become symmetric tensor triangulated
categories. Unfortunately, the product does not combine well with
the stabilization procedure of \eqref{sp.naive}, so that this
extension is very non-trivial already in the non-equivariant case,
see e.g. \cite{adams}. Although in the last fifteen years, new and
more satisfactory approaches appeared (e.g. \cite{M} or \cite{sym}),
still, none of them can be recounted in a few pages. Nevertheless,
whatever construction one uses, our equivariant homology functor
$h^G$ ought to be compatible with the smash products.

\begin{lemma}
For any two objects $X,Y \in \sthom(G)$, we have
$$
h^G(X \wedge Y) \cong h^G(X) \otimes h^G(Y).
$$
\end{lemma}

\proof[Sketch of a proof.] To give an honest proof, we would need to
use an exact definition of the product of spectra, which would be
beyond the scope of the paper; instead, we show how the statement is
deduced from the standard properties of the smash product. Write
down $X = (\Sigma^V)^{-1}(X')$, $Y = (\Sigma^W)^{-1}(y')$, $x',Y'
\in \Hom(G)$, $V,W \subset R^{\oplus \infty}$. By
Lemma~\ref{prod.hom.Mackey}, we have $h^G(X' \wedge Y') \cong
h^G(X') \otimes h^G(Y')$ and $h^G(S^{V \oplus W}) \cong h^G(S^V
\wedge S^W) \cong h^G(S^V) \otimes h^G(S^W)$, so that
$$
\sigma^{V \oplus W} \cong \sigma^V \circ \sigma^W \cong \sigma^W
\circ \sigma^V,
$$
where $\sigma^V$ is as in the proof of Corollary~\ref{g.homol}. Then
\begin{align*}
h^G(X \wedge Y) &\cong h^G((\Sigma^V)^{-1}X' \wedge
(\Sigma^W)^{-1}Y') \cong h^G((\Sigma^{V \oplus W})^{-1}(X' \wedge
Y')) \\
&\cong (\sigma^{V \oplus W})^{-1}h^G(X') \otimes h^G(Y') \cong
h^G(X) \cong h^G(Y),
\end{align*}
as required.
\endproof

An analogous statement for $\sthom^{naive}(G)$ is also obviously
true. Moreover, the functor $i_!$ of \eqref{comma} is also tensor,
and $q_!^{opp}$ is tensor by \eqref{prod.indu.2} of
Lemma~\ref{prod.indu.como}, so that \eqref{comma} is actually a
diagram of tensor functors.

Next, for any subgroup $H \subset G$, the fixed points functor $X
\mapsto X^H$ from $\Hom(G)$ to the category of CW complexes
preserves products; since the fixed points subset of a sphere $S^V$
is also a sphere, the fixed points functor is compatible with
\eqref{sp.naive} and extends to the category $\sthom(G)$ of finite
$G$-spectra. The result is called the {\em geometric fixed points
functor} and denoted by
$$
\Phi^H:\sthom(G) \to \sthom.
$$
It is also compatible with the smash product. Under our the
equivariant homology functor $h^G$, it goes to the fixed points
functor on the category $\DM(G)$:

\begin{lemma}
For any subgroup $H \subset G$ and any $X \in \sthom(G)$, we have a
quasiisomorphism
\begin{equation}\label{phi.phi}
h \circ \Phi^H (X) \cong \Phi^{[G/H]} \circ h^G(X)
\end{equation}
\end{lemma}

\proof{} By \eqref{prod.indu.2} of Lemma~\ref{prod.indu.como}, we
have
$$
\Phi^{[G/H]}(\sigma^V(M)) \cong \Phi^{[G/H]}(M)[\dim V^H]
$$
for any $V \subset R^{\otimes \infty}$ and any $M \in
\DM(G)$. Therefore the statement is compatible with suspensions, and
it suffices to prove it for $X = \Sigma^\infty(X_0)$ for some finite
$G$-CW complex $X_0 \in \Hom(G)$. Then it immediately follows from
\eqref{prod.indu.1} of Lemma~\ref{prod.indu.como}.
\endproof

Consider now the full categories of spectra $\Sthom^{naive}(G)$,
$\Sthom(G)$. The tautological functor $i_!$ of
\eqref{change.of.universe} extends to all spectra, although it
becomes rather non-trivial in the process. So do the smash product
and the geometric fixed point functors. Just as in the case of the
finite spectra, the functors $\Phi^H$, $H \subset G$, and $i_!$ are
tensor functors. It is natural to expect that our equivariant
homology functor $h^G$ also extends to $\Sthom(G)$.

\begin{conj}\label{conj.1}
The equivariant homology functor $h^G$ of Corollary~\ref{g.homol}
and the naive homology functor $h^G_{naive}$ of \eqref{h.g.naive}
extend to the categories of spectra, so that \eqref{comma} extends
to a commutative diagram
\begin{equation}\label{comma.bis}
\begin{CD}
\Sthom^{naive}(G) @>{i_!}>> \Sthom(G)\\
@V{h^G_{naive}}VV @VV{h^G}V\\
\D(O_G^{opp},\Z\amod) @>{q^{opp}_!}>> \DM(G)
\end{CD}
\end{equation}
of tensor triangulated functors. Moreover, for any subgroup $H
\subset G$, \eqref{phi.phi} extends to an isomorphism
$$
\Phi^{[G/H]} \circ h^G(X) \cong h \circ \Phi^H(X)
$$
of functors from $\Sthom(G)$ to $\D(\Z\amod)$.
\end{conj}

Moreover, on the level of spectra, the tautological functor $i_!$
acquires a right-adjoint $i^*:\Sthom(G) \to \Sthom^{naive}(G)$, and
we can consider the diagram
$$
\begin{CD}
\Sthom^{naive}(G) @<{i^*}<< \Sthom(G)\\
@V{h^G_{naive}}VV @VV{h^G}V\\
\D(O_G^{opp},\Z\amod) @<{q^{opp*}}<< \DM(G),
\end{CD}
$$
where $q^{opp*}$ is the restriction with respect to the embedding
$q^{opp}:O_G^{opp} \to \Q\Gamma_G$ -- that is, the right-adjoint
functor to $q_!^{opp}$. By base change, we have a natural map
$$
h^G_{naive} \circ i^* \to q^{opp*} \circ h^G.
$$

\begin{conj}\label{conj.2}
The base change map $h^G_{naive} \circ i^* \to q^{opp*} \circ h^G$
is an isomorphism.
\end{conj}

Using $i^*$, one can define another fixed point functor associated
to a subgroup $H \subset G$. Indeed, we have an obvious fixed points
functor
$$
\sthom^{naive}(G) \to \sthom,
$$
and it extends to a functor from $\Sthom^{naive}(G)$ to $\Sthom$
(isomorphic to $\Phi^H \circ i_!$); composing it with $i^*$, we
obtain a triangulated functor
$$
\Sthom(G) \to \Sthom.
$$
This functor is called the Lewis-May fixed points functor. There is
no standard letter associated to it; usually it is denoted simply by
$X \mapsto X^H$.

If Conjecture~\ref{conj.1} and Conjecture~\ref{conj.2} are known,
then the following is immediate.

\begin{corr}
For any $X \in \Sthom(G)$, we have a functorial isomorphism
$$
\Psi^{[G/H]}(h^G(X)) \cong h(X^H).
$$
\end{corr}

\proof{} By Conjecture~\ref{conj.2} and Conjecture~\ref{conj.1}, we
have
$$
h(X^H) \cong h(\Phi^H(i_!i^*(X))) \cong
\Phi^{[G/H]}(q_!^{opp}q^{opp*}h^G(X)),
$$
and by Lemma~\ref{prod.indu.como}, this is isomorphic to
$$
q^{opp*}(h^G(X))([G/H]) \cong \Psi^{[G/H]}(h^G(X)),
$$
as required.
\endproof

Let me conclude the paper with the following remark. As we have
noted in Remark~\ref{loop}, the categories $\Q\C(-,-)$ used in
Section~\ref{mack.der} in our definition of derived Mackey functors
are naturally symmetric monoidal, so that the group completions
$\Omega B|\Q\C(-,-)|$ of their classifying spaces are infinite loop
spaces. Thus one can define a category ${\sf B}_G$ ``enriched in
spectra'' with objects $[G/H]$ and morphisms given by $\Omega B|\Q
O_G(-,-)|$. If one does this with enough precision, one can then
define the category of ``enriched functors'' from ${\sf B}_G$ to
$\Sthom$ in such a way that it becomes a triangulated category. As I
understand, this construction is well-known in topology, and it is
well-known that the resulting category is $\Sthom(G)$. It seems that
this has not been written down (probably because the technology
needed to make this precise appeared not long ago and is still
rather cumbersome, and the difficulties of working out all the
details outweigh the benefits of giving yet another description of a
well-studied object). Be it as it may, if this sketch is assumed to
work, then the comparison with the present paper becomes quite
transparent: all we do is replace spectra with complexes which
compute their homology, and replace ``enriched functors'' with
$A_\infty$-functors. What I don't know is whether our third
description of $\DM(G)$, namely the one given in
Section~\ref{galois.sec}, has any counterpart with spectra instead
of complexes. It seems it would be very useful, since for Mackey
functors, this description is by far the most effective. On the
other hand, on the level of complexes one has to work with an
$A_\infty$-coalgebra $T^\C_\idot$, not an $A_\infty$-algebra, and
various finiteness phenomena become crucially important; it is not
clear whether this can be made to work in $\Sthom$.

\bigskip

\noindent
{\sc
Steklov Math Institute\\
Moscow, USSR}

\bigskip

\noindent
{\em E-mail address\/}: {\tt kaledin@mccme.ru}

\end{document}